\documentclass[11pt, oneside]{amsart}
\usepackage[text={6.5in,9in},centering,letterpaper,dvips]{geometry}
\usepackage{graphicx}
\usepackage{amsfonts}
\usepackage[dvipsnames]{xcolor}
\usepackage{subcaption}
\usepackage{epsf}
\usepackage{amssymb}
\usepackage{amsmath}
\usepackage{amscd}
\usepackage{tikz-cd}
\usepackage{adjustbox}
\usepackage{pdfpages}
\usepackage{fancyhdr}
\usepackage{setspace}
\usepackage{soul} 
\usepackage[all]{xy}
\usepackage{verbatim}
\usepackage{enumerate}
\usepackage[colorlinks=true, urlcolor=NavyBlue, linkcolor=NavyBlue, citecolor=NavyBlue,backref=page]{hyperref}
\usepackage{etoolbox}
\usepackage{bm}

\renewcommand*{\backref}[1]{}
\renewcommand*{\backrefalt}[4]
{%
    \ifcase #1 (Not cited.)%
        \or        (Cited on page~#2.)
        \else      (Cited on pages~#2.)
    \fi
}

\makeatletter

\renewcommand{\p@subfigure}{\thefigure.} 
\makeatother

\makeatletter
\renewcommand{\tocsection}[3]{%
  \indentlabel{\@ifnotempty{#2}{\bfseries\ignorespaces#1 #2\quad}}\bfseries#3}
\renewcommand{\tocsubsection}[3]{%
  \indentlabel{\@ifnotempty{#2}{\ignorespaces#1 #2\quad}}#3}

\newcommand\@dotsep{4.5}
\def\@tocline#1#2#3#4#5#6#7{\relax
  \ifnum #1>\c@tocdepth 
  \else
    \par \addpenalty\@secpenalty\addvspace{#2}%
    \begingroup \hyphenpenalty\@M
    \@ifempty{#4}{%
      \@tempdima\csname r@tocindent\number#1\endcsname\relax
    }{%
      \@tempdima#4\relax
    }%
    \parindent\z@ \leftskip#3\relax \advance\leftskip\@tempdima\relax
    \rightskip\@pnumwidth plus1em \parfillskip-\@pnumwidth
    #5\leavevmode\hskip-\@tempdima{#6}\nobreak
    \leaders\hbox{$\m@th\mkern \@dotsep mu\hbox{.}\mkern \@dotsep mu$}\hfill
    \nobreak
    \hbox to\@pnumwidth{\@tocpagenum{\ifnum#1=1\bfseries\fi#7}}\par
    \nobreak
    \endgroup
  \fi}
\AtBeginDocument{%
\expandafter\renewcommand\csname r@tocindent0\endcsname{0pt}
}
\def\l@subsection{\@tocline{2}{1pt 
}{5pc 
}{}{}}
\makeatother

\theoremstyle{theorem}
\newtheorem{theorem}{Theorem}[section]
\newtheorem{proposition}[theorem]{Proposition}
\newtheorem{lemma}[theorem]{Lemma}
\newtheorem{question}[theorem]{Question}
\newtheorem{corollary}[theorem]{Corollary}
\newtheorem{conjecture}[theorem]{Conjecture}

\newtheorem*{namedthm}{\namedthmname}
\newcounter{namedthm}

\makeatletter 
\newenvironment{named}[1]
  {\def\namedthmname{#1}%
   \refstepcounter{namedthm}%
   \namedthm\def\@currentlabel{#1}}
  {\endnamedthm}
\makeatother

\makeatletter
\newcommand{\newreptheorem}[2]{%
\newtheorem*{rep@#1}{\rep@title}
\newenvironment{rep#1}[1]{%
\def\rep@title{#2 \ref{##1}}
\begin{rep@#1}}{\end{rep@#1}}}
\makeatother

\newreptheorem{theorem}{Theorem}
\newreptheorem{lemma}{Lemma}
\newreptheorem{question}{Question}
\newreptheorem{corollary}{Corollary}
\newreptheorem{proposition}{Proposition}
\newreptheorem{conjecture}{Conjecture}

\theoremstyle{definition}
\newtheorem{definition}[theorem]{Definition}
\newtheorem{remark}[theorem]{Remark}

\newtheorem{example}[theorem]{Example}

\makeatletter
\newcommand{\newreptheoremD}[2]{%
\newtheorem*{rep@#1}{\rep@title}
\newenvironment{rep#1}[1]{%
\def\rep@title{#2 \ref{##1}}
\begin{rep@#1}}{\end{rep@#1}}}
\makeatother

\newreptheoremD{definition}{Defintion}

\newcommand{\Z}{\mathbb{Z}}
\newcommand{\C}{\mathbb{C}}
\newcommand{\N}{\mathbb{N}}

\newcommand{\R}{\mathbb{R}}
\newcommand{\RP}{\mathbb{RP}}
\newcommand{\CP}{\mathbb{CP}}

\newcommand{\id}{\text{id}}
\newcommand{\Int}{\text{Int}}

\newcommand{\Bb}{\mathcal B}
\newcommand{\Cc}{\mathcal C}
\newcommand{\Dd}{\mathcal D}
\newcommand{\Ee}{\mathcal E}
\newcommand{\Ff}{\mathcal F}

\newcommand{\Ii}{\mathcal I}

\newcommand{\Ll}{\mathcal L}

\newcommand{\Ss}{\mathcal S}
\newcommand{\Tt}{\mathcal T}
\newcommand{\Vv}{\mathcal V}

\newcommand{\Zz}{\mathcal Z}
\newcommand{\Mod}{\text{Mod}}

\newcommand{\im}{\text{Im}}

\newcommand{\DD}{\mathfrak{D}}

\newcommand{\TT}{\mathfrak{T}}

\newcommand{\Diff}{\text{Diff}}
\newcommand{\Sing}{\text{Sing}}
\newcommand{\Fix}{\text{Fix}}

\newcommand{\Stab}{\text{Stab}}

\newcommand{\Diag}{\text{Diag}}

\newcommand{\bigast}{\mathop{\scalebox{1.5}{\raisebox{-0.2ex}{$\ast$}}}}%

\makeatletter
\def\@seccntformat#1{%
  \protect\textup{\protect\@secnumfont
    \ifnum\pdfstrcmp{subsection}{#1}=0 \bfseries\fi
    \csname the#1\endcsname
    \protect\@secnumpunct
  }%
}  
\makeatother

\begin{document}

\rhead{\thepage}
\lhead{\author}
\thispagestyle{empty}


\raggedbottom
\pagenumbering{arabic}
\setcounter{section}{0}


\title{Equivariant Trisections for Group Actions on Four-Manifolds}

\author{Jeffrey Meier}
\address{Department of Mathematics, Western Washington University 
Bellingham, WA 98229}
\email{jeffrey.meier@wwu.edu}
\urladdr{http://jeffreymeier.org} 

\author{Evan Scott}
\address{Department of Mathematics, CUNY Graduate Center
New York, NY 10016}
\email{escott@gradcenter.cuny.edu}
\urladdr{https://sites.google.com/view/evanscott}

\begin{abstract}
	Let $G$ be a finite group, and let $X$ be a smooth, orientable, connected, closed 4--dimensional $G$--manifold.
	Let $\Ss$ be a smooth, embedded, $G$--invariant surface in $X$.
	We introduce the concept of a $G$--equivariant trisection of $X$ and the notion of $G$--equivariant bridge trisected position for $\Ss$ and establish that any such $X$ admits a $G$--equivariant trisection such that $\Ss$ is in equivariant bridge trisected position.
	Our definitions are designed so that $G$--equivariant (bridge) trisections are determined by their spines; hence, the 4--dimensional equivariant topology of a $G$--manifold pair $(X,\Ss)$ can be reduced to the 2--dimensional data of a $G$--equivariant shadow diagram.
	
	As an application, we discuss how equivariant trisections can be used to study quotients of $G$--manifolds.
	We also describe many examples of equivariant trisections, paying special attention to branched covering actions, hyperelliptic involutions, and linear actions on familiar manifolds such as $S^4$, $S^2\times S^2$, and $\CP^2$.
	We show that equivariant trisections of genus at most one are geometric, and we give a partial classification for genus-two.
\end{abstract}

\maketitle

\vspace{-10mm}
\tableofcontents

\section{Introduction}\label{sec:intro}

Since trisections were introduced by Gay and Kirby as a novel approach to studying 4-manifolds \cite{GayKir_16_Trisecting-4-manifolds}, a question has lingered regarding how these objects interface with finite group actions on $4$--manifolds.
The present paper answers this question by introducing the notion of an \emph{equivariant trisection} that satisfies equivariant analogs of many of the fundamental results about trisections.
Most importantly, the equivariant topology of an equivariantly trisected manifold is determined entirely by the equivariant topology of the $3$--dimensional spine of the trisection.
This gives, for example, that a group action on a trisection diagram uniquely defines an equivariantly trisected group action on the trisected manifold.

The study of finite group actions on $4$--manifolds has a rich history in both the topological and smooth categories, but in the smooth category many fundamental questions remain unanswered; see~\cite{Che_10_Group-actions-on-4-manifolds:} and~\cite{Edm_18_A-survey-of-group-actions-on-4-manifolds} for surveys.
In dimension three, equivariant questions such as the Smith Conjecture served as important catalysts for development of powerful tools and ideas~\cite{MorBas_84_The-Smith-conjecture}, and we hope the present work will draw attention to (and offer new perspective on) the 4--dimensional setting.

Central among these questions are those concerning \emph{linearity} of cyclic actions.
For example, it is known by the failure of the Smith conjecture in dimension four~\cite{Gif_66_The-generalized-Smith-conjecture,Gor_74_On-the-higher-dimensional-Smith} that the fixed-point set of a smooth $\Z_p$--action on $B^4$ can be knotted, so such actions cannot be smoothly equivalent to an action by a subgroup of $SO(4)$.
It is not known, however, whether there exist smooth, nonlinear $\Z_p$--actions on $B^4$ fixing an isolated point or fixing an unknotted disk; see Section~\ref{sec:open_questions_further_directions} for a discussion of open problems, and see Section~\ref{sec:background} for background on finite group actions.

The falsity of the Smith conjecture in dimension four, contrasting its validity in dimension three, is an obstruction to a $4$--dimensional group action being determined by the group action on the $3$--dimensional boundary.
Since we want equivariant trisections to be determined by their spines, we must require something extra of the sectors of an equivariant trisection: these pieces must be \emph{linearly parted} equivariant 4--dimensional 1--handlebodies; see Definition~\ref{def:linearly_parted} and the ensuing discussion in Section~\ref{sec:definitions}, as well as~\cite[Section~3]{MeiSco_LP}.
With this in mind, we have our main definition.

\begin{repdefinition}{def:equiv_tri}
	Let $X$ be a 4--dimensional $G$--manifold.
	A \emph{$G$--equivariant trisection of $X$} is a decomposition
	$$X = X_1\cup X_2\cup X_3$$
	such that, for each $i\in\Z_3$:
	\begin{enumerate}
		\item $X_i$ is an invariant 4--dimensional 1--handlebody $\natural^{k_i}(S^1\times B^3)$, and the induced $G$--action on $X_i$ is linearly parted;
		\item $H_i = X_{i-1}\cap X_i$ is an invariant 3--dimensional 1--handlebody; and
		\item $\Sigma = X_1\cap X_2\cap X_3$ is an invariant, closed, genus--$g$ surface.
	\end{enumerate}	
	The union $H_1\cup H_2\cup H_3$ is called the \emph{spine} of the trisection.
\end{repdefinition}

Given a $G$--invariant surface $\Ss$ in $X$, the above definition can be extended to formalize the notion that $\Ss$ be in \emph{equivariant bridge trisected position} with respect to $\TT$ (see Definition~\ref{def:equiv_bridge}), giving an equivariant version of the theory of bridge trisections first introduced by the first author and Zupan~\cite{MeiZup_17_Bridge-trisections,MeiZup_18_Bridge-trisections}.
Our first main result establishes that equivariant (bridge) trisections always exist.
The proof, which proceeds by taking an equivariant triangulation of $X$ and explicitly decomposing each $4$--dimensional simplex, is given in Section~\ref{sec:existence}.

\begin{reptheorem}{thm:exist_main}
	Let $G$ be a finite group, let $X$ be a smooth, orientable, connected, closed, $4$--dimensional $G$--manifold, and let $\{\Ss_i\}_{i=1}^N$ be a (possibly empty) collection of embedded, closed, $G$--invariant surfaces in $X$.
	Then, $X$ admits a $G$--equivariant trisection for which each $\Ss_i$ is in equivariant bridge trisected position.
\end{reptheorem}

As mentioned above, an equivariant (bridge) trisection is determined by its spine.

\begin{reptheorem}{thm:spines_bridge} 
	Let $G$ be a finite group.
	Let $\mathfrak T$ and $\mathfrak T'$ be $G$--equivariant bridge trisections for $(X,\Ss)$ and $(X',\Ss')$, respectively.
	If the restrictions of the actions to the spines of the bridge trisections are $G$--diffeomorphic, then the bridge trisections are $G$--diffeomorphic.
\end{reptheorem}

These results rely on theorems from recent work of the authors \cite{MeiSco_LP} that gives equivariant analogs of the foundational theorem of Laudenbach and Po\'enaru~\cite{LauPoe_72_A-note-on-4-dimensional-handlebodies} and the corresponding theorem involving boundary-parallel disks first established in a special case by Livingston~\cite{Liv_82_Surfaces-bounding-the-unlink} and later generalized by the first author and Zupan~\cite[Lemma~8]{MeiZup_18_Bridge-trisections}. 

It is fully understood when the quotient $X/G$ of a smooth group action of $G$ on a smooth $4$--manifold $X$ is a \textit{topological} manifold, and in all such cases there exists at least one smooth structure on $X/G$ compatible with the quotient map.
However, in most cases, it is not known whether such a smooth structure on $X/G$ is unique; see \cite[Appendix]{HamHau_11_Conjugation-spaces-and-4-manifolds} for more on this problem.
Regardless, we show that there is a canonical smooth structure on $X/G$ associated to any equivariant trisection $\TT$ of $X$, expressed as the quotient trisection $\TT/G$ of the topological $4$--manifold $X/G$.
In fact we show slightly more.

\begin{reptheorem}{thm:quotient_trisection}
	Let $\TT$ be a $G$--equivariant trisection of $X$.
	Let $N\triangleleft G$ be normal.
	Then, the following are equivalent.
	\begin{enumerate}
		\item $X/N$ admits a smooth structure making it a $G/N$--manifold.
		\item $\TT/N$ is a $G/N$--equivariant trisection of $X/N$.
		\item $\Stab_N(x)$ acts as a rotation group on $\nu(x)$ for each $x\in X$.
	\end{enumerate}
\end{reptheorem}

We also establish an analogous quotient theorem for bridge trisections (Theorem~\ref{thm:quotient_bridge_trisection}), which requires the extra hypothesis that the induced action of $N$ on a bridge trisected surface $\Ss$ have no elements with $1$--dimensional fixed point set.
In particular, Theorem~\ref{thm:quotient_bridge_trisection} applies to the fixed-point set of a group action when that fixed-point set is a disjoint union of embedded surfaces, such as when $G$ is a prime-order cyclic group and $X/G$ is a smooth manifold.

We'd like to now mention a few results that make precise the idea that the equivariant topology of $X$ is controlled by the $G$--action on the spine of an equivairant trisection of $X$.
The most important example of this is the following characterization of equivariant (bridge) trisections in terms of \emph{invariant shadow diagrams}; see Subsection~\ref{subsec:spines}.

\begin{repproposition}{prop:diag_spine}
	Fix a finite group $G$.
	\begin{enumerate}
		\item A $G$--invariant shadow diagram determines a unique $G$--equivariant (bridge) trisection, up to $G$--diffeomorphism.
		\item $G$--diffeomorphic $G$--invariant shadow diagrams determine $G$--diffeomorphic $G$--equivariant (bridge) trisections.
		\item Every $G$--equivariant (bridge) trisection can be represented by a $G$--invariant shadow diagram.
	\end{enumerate}
\end{repproposition}

There are other surprising manifestations of this idea.
Corollary~\ref{cor:bridge_in_spine} shows that the fixed-point set of an action must be in bridge trisected position whenever it intersects the spine of the trisection in the spine of a bridge trisected surface.
Corollary~\ref{cor:quotient_diagram} establishes that the quotient $X/G$ is a manifold if and only if the quotient $\DD/G$ of any $G$--invariant shadow diagram $\DD$ for $X$ is a shadow diagram.

Section~\ref{sec:examples} is devoted to studying examples of equivariant trisections in a variety of important settings.
In Subsection~\ref{subsec:branched_coverings}, we discuss how regular branched coverings of bridge trisected surfaces are naturally equivariantly trisected (Proposition~\ref{prop:branched_cover}), and we prove a structural lemma for linearly parted branched covering actions on handlebodies  that, in some cases, allows us to conclude that the singular set of an equivariant trisection is `automatically' bridge trisected (Lemma~\ref{lem:branched_parted_sub}).
We then explore an interesting example of a $Q_8$--equivariant trisection of $S^2\times S^2$ whose quotient is a link of unknotted projective planes in $S^4$ before showing, in Corollary~\ref{cor:Q8}, that every $Q_8$--equivariant trisection is irreducible.
In Subsection~\ref{subsec:hyperelliptic_trisections}, we discusses trisections that are equivariant with respect to a hyperelliptic involution on the trisection diagram, and we show in Theorem~\ref{thm:hyperelliptic} that all such trisections are branched double covers of knotted surfaces in $S^4$.
We use this to give a simple proof (Theorem~\ref{thm:massey-kuiper}) that the quotient of $\CP^2$ by complex conjugation is $S^4$, recovering an old result of Massey, Kuiper, and Arnol'd~\cite{Mas_73_The-quotient-space-of-the-complex-projective, Kui_74_The-quotient-space-of-bf-CP2-by-complex-conjugation, Arn_88_The-branched-covering-bf-Crm-P2to}, and we give a similar analysis of involutions of $S^2\times S^2$.
Finally, in Subsection~\ref{subsec:PU(3)_trisections}, we consider finite subgroups of the linear action of $PU(3)$ on $\CP^2$, and we give equivariant trisections for the vast majority of all such actions on $\CP^2$.

In Subsections~\ref{subsec:genus_0_classification} and~\ref{subsec:genus_1_classification}, we classify equivariant trisections of genus-zero and genus-one, showing they are always well-adapted to the geometry of the underlying manifold, which is always $S^4$, $\CP^2$, or $S^1\times S^3$.

\begin{repcorollary}{cor:geometric}
	Every group action admitting a genus-one equivariant trisection is geometric.
\end{repcorollary}

In Subsection~\ref{subsec:genus_two}, we adapt an approach of Bruno Zimmermann for studying equivariant Heegaard splittings~\cite{Zim_96_Genus-actions-of-finite} to the setting of equivariant trisections and classify all genus-two trisections which are \emph{maximally symmetric}.
The resulting actions again turn out to be quite simple.

Finally, in Section~\ref{sec:open_questions_further_directions}, we outline interesting open questions regarding finite group actions on $4$--manifolds, as well as further directions in which to pursue the study of equivariant trisections.
The equivariant analogs to several fundamental theorems about trisections are still missing.
In particular, we have no results about uniqueness up to stabilization nor any notion of a set of equivariant bridge trisection moves that are sufficient to relate equivariantly isotopic invariant surfaces.

\subsection*{Acknowledgements}

The authors wish to thank David Gabai, Malcolm Gabbard, David Gay, Daniel Hartman, Maggie Miller, and Michael Young for helpful conversations about various aspects of this project.
The first author was supported by NSF grants DMS-2006029 and DMS-2405324.

\section{Background on finite group actions}
\label{sec:background}

In this section, we recall basic definitions about group actions.
We refer the reader to~\cite{Bre_72_Introduction-to-compact-transformation-groups} for a thorough treatment of basic equivariant topology with complete details.

Throughout the paper, $G$ denotes a finite group, and all group actions we consider are faithful, smooth, and properly discontinuous.
Let $M$ be a manifold, and let $G$ act on $M$.
Unless explicitly indicated otherwise, we assume that every element of $G$ acts on $M$ as an orientation-preserving diffeomorphism, i.e. the group acts \emph{orientation-preservingly} on $M$.
We encode this set-up by saying that $M$ is a \emph{$G$--manifold}.

To compare different actions of a given group $G$ on a given manifold $M$, we regard a group action as a homomorphism $\rho\colon G\to \Diff(M)$.
We sometimes write $(M, \rho)$ to denote the $G$--manifold $M$ with $G$ acting by $\rho$.
Strictly speaking, we identify a group action $\rho$ with its image in $\Diff(M)$, so that $\rho$ is identified with $\rho\circ\phi$ for any automorphism $\phi$ of $G$. 
This amounts to relabelling the elements in the image of $\rho$ and will not be regarded as a distinct action.
Since we assume all actions of $G$ are faithful, we assume $\rho$ is injective.
Two actions $\rho_1$ and $\rho_2$ of $G$ on $M$ are \emph{(smoothly) equivalent} if there exists $\psi\in \Diff(M)$ so that $\psi^{-1}\rho_1(g)\psi = \rho_2(g)$ for all $g\in G$, i.e. the two actions are conjugate in $\Diff(M)$\footnote{
Again, strictly speaking we ask that the images are conjugate subgroups by $\phi$, but we can always choose isomorphisms $\rho_1$ and $\rho_2$ of $G$ onto these two subgroups which satisfy the above equation on the nose.}.
This coincides with the notion that $(M, \rho_1)$ and $(M, \rho_2)$ are equivariantly diffeomorphic by a diffeomorphism $\psi$.

A submanifold $N\subseteq M$ is \emph{invariant} if $g\cdot N = N$ for all $g\in G$.
An invariant submanifold always admits an \emph{open invariant tubular neighborhood}, $\nu(N)\subseteq M$, which is equivariantly diffeomorphic to a (smooth) $G$--vector bundle over $N$; any such open invariant tubular neighborhood contains a closed invariant tubular neighborhood~\cite[Chapter~VI, Theorem~2.2]{Bre_72_Introduction-to-compact-transformation-groups}.
We will sometimes use $\nu(N)$ to denote a \textit{closed} invariant tubular neighborhood and will make this distinction explicit when it is not clear from context.
A submanifold $N\subseteq M$ is \emph{equivariant} if, for all $g\in G$, either $g\cdot N = N$ or $(g\cdot N)\cap N = \varnothing$.
Clearly all invariant submanifolds are equivariant.

Let $N\subseteq M$ be an equivariant submanifold. 
The \emph{stabilizer} of $N$ is the subgroup $\Stab_G(N)<G$ consisting of those $g\in G$ such that $g\cdot N = N$.
We will sometimes refer to the stabilizer of $N$ as the \emph{set-wise stabilizer} for emphasis.
The \emph{point-wise stabilizer} of $N$ is the set of elements $g\in G$ that fix each point in $N$.
Note that the point-wise stabilizer of $N$ is normal in the set-wise stabilizer of $N$; we denote the corresponding quotient by $G_N$.
The quotient group $G_N$ acts on $N$, and we refer to this action as the \emph{induced action} of $G$ on $N$.
We say that $G$ acts on $N$ with some property (say, orientation-preservingly) if $G_N$ acts on $N$ with that property.
Note that if $N$ is codimension-zero (such as when $N$ is a tubular neighborhood, as will often be the case) then $G_N = \Stab_G(N)$: if $g\in G$ point-wise fixes a codimension-zero submanifold of (connected) manifold $M$, then $g$ fixes all of $M$ and is hence trivial by faithfulness.
Therefore, for orientation-preserving $G$--actions, $G$ always acts on $\nu(N)$ orientation-preservingly, even when $G$ does not act orientation-preservingly on $N$.

The \emph{fixed-point set} of $G$ is the set of points $x\in M$ such that $g\cdot x = x$ for all $g\in G$ and is denoted $\Fix(G)$.
The fixed-point set of $G$ is always a disjoint union of closed submanifolds of $M$, though the dimension can vary across connected components.
The \emph{singular set} of $G$ is the set of points $x\in M$ such that $g\cdot x = x$ for some nontrivial $g\in G$ and is denoted $\Sing(G)$.
The singular set is decomposed into the fixed-point sets of the subgroups $H\leq G$.
Thus, the singular set is a union of submanifolds of $M$, though these submanifolds need not be disjoint or transverse in general.
In this paper, $M$ will have dimension at most four, and $G$ will usually act via orientation-preserving diffeomorphisms; in this case, the components of $\Sing(G)$ fixed by cyclic subgroups will either have dimension-zero or codimension-two.
If there is a component of $\Sing(G)$ which is codimension $1$, the (point-wise) stabilizer of that component acts orientation-reversingly.

A group action is \emph{free} if every element of $G$ is a fixed-point-free map $g\colon M\to M$.
A group action is \emph{semi-free} if the action is not free, and every point of $M$ with nontrivial stabilizer is fixed by $G$.
In other words, in a semi-free action, the singular set of $G$ coincides with the fixed-point set.
For example, all actions of $\Z_p$ with $p$ a prime are semi-free.
We say an action is \emph{pseudo-free} if the action is not free, and the fixed-point set of every element of $G$ is $0$--dimensional.
In other words, the singular set of a pseudo-free action is discrete.

\begin{remark}
\label{rmk:action_terms}
	The definitions given above appear frequently in the literature with minor alterations, and the terms can be confusing.
	We offer the following mnemonics.
	The term ``semi-free'' can be abbreviated SF, which might be understood to stand for ``singularity-free,'' reflecting the fact that a semi-free action has no \emph{proper} singular points, just fixed points.
	The term ``pseudo-free'' can be abbreviated PF, which might be understood to stand for ``point-fixing,'' reflecting the fact that each element in a pseudo-free action fixes a collection of points.
\end{remark}

\section{Definitions and basic notions}
\label{sec:definitions}

In this section, we define equivariant (bridge) trisections.
These definitions depend on some ancillary concepts that require further development and discussion:
In Subsection~\ref{subsec:linearly_parting}, we discuss linearly parted actions on handlebodies, and in Subsection~\ref{subsec:tangles}, we extend this discussion to the setting of trivial tangles and trivial disk tangles.
In Subsection~\ref{subsec:spines}, we discuss how equivariant (bridge) trisections are determined by their spines and by their diagrams.

First, we recall our main definition from the introduction.

\begin{definition}
\label{def:equiv_tri}
	Let $X$ be a 4--dimensional $G$--manifold.
	A \emph{$G$--equivariant trisection of $X$} is a decomposition
	$$X = X_1\cup X_2\cup X_3$$
	such that, for each $i\in\Z_3$:
	\begin{enumerate}
		\item $X_i$ is an invariant 4--dimensional 1--handlebody $\natural^{k_i}(S^1\times B^3)$, and the induced $G$--action on $X_i$ is linearly parted;
		\item $H_i = X_{i-1}\cap X_i$ is an invariant 3--dimensional 1--handlebody; and
		\item $\Sigma = X_1\cap X_2\cap X_3$ is an invariant, closed, genus--$g$ surface.
	\end{enumerate}
	To record the complexity of the trisection, we sometimes refer to it as a $(g;\bm k)$--trisection, where $\bm k = (k_1,k_2,k_3)$, or a $(g,k)$--trisection if $k=k_i$ for all $i\in\Z_3$; the parameter $g$ is called the \emph{genus} of the trisection.
	
	The union $H_1\cup H_2\cup H_3$ is called the \emph{spine} of the trisection.
\end{definition}

The notion of an equivariant bridge trisection adapts the equivariant trisection framework to incorporate invariant surfaces in $X$, generalizing the last definition.

\begin{definition}
\label{def:equiv_bridge}
	Let $X$ be a 4--dimensional $G$--manifold, and let $\TT$ be a $G$--equivariant trisection of $X$.
	Let $\Ss$ be a $G$--invariant, embedded surface in $X$.
	We say that $\Ss$ is in \emph{equivariant bridge trisected position} with respect to $\TT$ if, for each $i\in\Z_3$:
	\begin{enumerate}
		\item $\Dd_i = X_i\cap \Ss$ is an invariant $p_i$--patch disk-tangle, and $(X_i, \Dd_i)$ is linearly parted as a pair; and
		\item $\Tt_i = H_i\cap \Ss$ is an invariant $b$--bridge trivial tangle.
	\end{enumerate}
	In this case, we refer to $\TT$ as a \emph{$G$--equivariant bridge trisection} of the pair $(X,\Ss)$.
	To record the complexity of the bridge trisection, we refer to it as a $(b;\bm p)$--bridge trisection, where $\bm p = (p_1,p_2,p_3)$, or a $(b,p)$--bridge trisection if $p=p_i$ for all $i\in\Z_3$; the parameter $b$ is called the \emph{bridge number} of the trisection.
	
	The union $(H_1,\Tt_1)\cup(H_2,\Tt_2)\cup(H_3,\Tt_3)$ is called the \emph{spine} of the bridge trisection.
\end{definition}

We have the following notions of equivalence for equivariant (bridge) trisections.

\begin{definition}
\label{def:trisection_equivalence}
	Let $\TT$ and $\TT'$ be equivariant trisections for $G$--manifolds $X$ and $X'$.
	We say that $\TT$ and $\TT'$ are \emph{$G$--diffeomorphic} if there exists an equivariant diffeomorphism $\phi\colon X\to X'$ such that, up to permuting the indices, $\phi(X_i) = X_i'$ for each $i\in\Z_3$.
	
	Let $\TT$ and $\TT'$ be equivariant bridge trisections for $G$--manifold pairs $(X,\Ss)$ and $(X',\Ss')$.
	We say that $\TT$ and $\TT'$ are \emph{$G$--diffeomorphic} if there exists an equivariant diffeomorphism of pairs $\phi\colon (X,\Ss)\to (X',\Ss')$ such that, up to permuting the indices, $\phi(X_i,\Dd_i) = (X_i',\Dd_i')$ for each $i\in\Z_3$.
\end{definition}

Note that the genus-one trisections of $\CP^2$ and $\overline\CP^2$ become diffeomorphic after transposing a pair of indices, so when it is necessary to be careful about orientation, only cyclic permutations of the indices should be allowed.

\subsection{Linearly parted actions}
\label{subsec:linearly_parting}

A key feature of the theory of trisections is the topological simplicity of $1$--handlebodies in dimensions three and four.
This feature is imported into the present work by breaking up $4$--dimensional $G$--manifolds into simple pieces for which the induced actions can be understood as amalgamations of linear actions.
We now recall some basic facts about linear actions.
We refer to \cite[Section~2]{MeiSco_LP} for further details.

\begin{definition}
	A $G$--action on $B^n$ is \emph{linear} if it is smoothly equivalent to the action of an $n$--dimensional, orthogonal representation of $G$ on $B^n$.
	A $G$--action on $S^{n-1}$ is \emph{linear} if it is the boundary of a linear $G$--action on $B^n$.
\end{definition}

Note that any real representation $G\hookrightarrow GL_n(\R)$ is equivalent to an orthogonal representation $G\hookrightarrow O(n)$.
This follows from the fact that every inner product is equivalent to the standard one up to change-of-basis (by Gram-Schmidt), together with the fact that we can always produce a $G$--invariant inner product by averaging over $G$.
Now, taking the change-of-basis that takes our $G$--invariant inner product to the standard one, we conjugate our $G$--representation to one that preserves the standard inner product, i.e. an orthogonal representation.
See~\cite[Chapter~0]{Bre_72_Introduction-to-compact-transformation-groups} for terminology and details.
This justifies the orthogonality hypothesis in the above definition.

The next definition makes precise the language in Definition~\ref{def:equiv_tri}.

\begin{definition}
\label{def:linearly_parted}
	Let $H$ be an $n$--dimensional 1--handlebody.
	A \emph{ball-system} for $H$ is a disjoint union $\Bb$ of neatly embedded\footnote{
	A submanifold $N$ is \emph{neatly embedded} in $M$ if $\partial N = N\cap\partial M$ and $N$ is transverse to $\partial M$~\cite[Chapter~1.4]{Hir_94_Differential-topology}.
	}
	$(n-1)$--balls in $H$ such that $H\setminus\nu(\Bb)$ is a disjoint union of $n$--balls.
	We say that a $G$--action on $H$ is \emph{parted} if there exists a $G$--invariant ball-system $\Bb$ for $H$.
	We say that a $G$--action on $H$ is \emph{linearly parted} if it is parted by $\Bb$, and the action of $G$ on each $(n-1)$--ball of $\Bb$ and on each $n$--ball component of $H\setminus\Bb$ is linear.
\end{definition}

In other words, the $G$--manifold $H$ is linearly parted if and only if it admits a decomposition into equivariant $0$--handles and $1$--handles such that the induced action on each handle is linear.
For this reason, we refer $H$ as a \emph{$G$--equivariant} $1$--handlebody.
Such a decomposition gives rise to a $G$--invariant Morse function with critical points of only index--$0$ and index--$1$, and such a Morse function gives rise to an equivariant handle decomposition in this sense.
We refer the reader to~\cite{Was_69_Equivariant-differential-topology} for a complete development of equivariant Morse theory and handlebody theory; see also Subsection~\ref{subsec:handles}, where we describe how equivariant trisections give rise to equivariant handle-decompositions of $4$--manifolds.

If $H$ is a $3$--dimensional $1$--handlebody, then any group action is automatically linearly parted by deep theorems from $3$--manifold topology; see~\cite[Section~3.1]{MeiSco_LP} for a complete summary.
In this sense, the requirement that the action on our $4$--dimensional $1$--handlebodies be linearly parted forces the $4$--dimensional pieces into analogy with the $3$--dimensional pieces.
However, the most important reason for requiring the actions on the $4$--dimensional $1$--handlebodies to be linearly parted is that it allows for the equivariant versions of the Laudenbach-Po\'enaru theorem that are proved in~\cite{MeiSco_LP}.
A major consequence of this is that equivariant (bridge) trisections are determined by their spines and can be represented by diagrams; see Subsection~\ref{subsec:spines}.
 
In what follows we will make use of the Equivariant Loop Theorem of Meeks and Yau~\cite{MeeSimYau_82_Embedded-minimal,MeeYau_79_The-classical-Plateau,MeeYau_80_Topology-of-three-dimensional}, as stated by Edmonds~\cite{Edm_86_A-topological-proof}, who gave a topological proof.

\stepcounter{theorem}
\begin{named}{Equivariant Loop Theorem~\thetheorem}
\label{ELT}
	Let $G$ act on a $3$--manifold $M$.
	Let $C\subset\partial M$ be a simple closed curve such that
	\begin{enumerate}
		\item $C$ is null-homotopic in $M$,
		\item $C$ is $G$--equivariant, and
		\item $C$ is transverse to the singular set $E$ of the action of $G$ on $\partial M$.
	\end{enumerate}
	Then there is a $G$--equivariant, neatly embedded disk $D\subset M$ with $\partial D=C$.
\end{named}

\subsection{Trivial tangles and disk-tangles}
\label{subsec:tangles}

We now turn our attention to explaining some of the language in Definition~\ref{def:equiv_bridge}: namely, the notions of equivariantly boundary-parallel trivial tangles and trivial disk tangles, as well as what it means for an action on such an object to be linearly parted as a pair.

\begin{definition}
	A \emph{$b$--bridge tangle} in a $3$--dimensional $1$--handlebody $H$ is a disjoint collection $\Tt$ of $b$ neatly embedded arcs in $H$.
	A tangle $\Tt$ is \emph{trivial} (or \emph{boundary-parallel}) if $\Tt$ can be isotoped relative to $\partial T$ to lie in $\partial H$.
	Equivalently, if $\Tt = \cup_{i=1}^b\tau_i$, there are pairwise disjoint disks $\Delta = \cup_{i=1}^b\delta_i$ such that $\partial \delta_i =\tau_i\cup a_i$, where $a_i$ is an arc in $\partial H$.
	The disks of $\Delta$ are collectively called \emph{bridge-disks} for the tangle $\Tt$.
\end{definition}

A tangle being boundary-parallel is a notion of unknottedness for a collection of arcs in a handlebody.
Adapting this definition to the equivariant case presents a major concern: consider the $\Z_n$--action on $B^3$ generated by a rotation.
The fixed arc of the action is clearly unknotted, yet cannot be equivariantly isotoped into the boundary.
To resolve this problem, we give the following definition.

\begin{definition}
	Let $G$ act on a $3$--dimensional $1$--handlebody $H$.
	An invariant tangle $\Tt = \cup_i\tau_i$ in $H$ is \emph{equivariantly boundary-parallel} if it has a collection $\Delta = \cup_i\delta_i$ of bridge-disks that is $G$--invariant and \emph{almost-disjoint}: the collection of disks is disjoint except possibly intersecting along strands of $\Tt$.
\end{definition}

For an example, the total preimage of a bridge-disk under a cyclic branched covering of a handlebody along a tangle is an almost-disjoint collection of bridge-disks with respect to the action of the deck group.
A single representative of a disjoint orbit of submanifolds is called an equivariant submanifold; similarly we say that a single bridge-disk from an almost-disjoint collection is \emph{almost-equivariant}, which is equivalent to being disjoint from (or equal to) each disk of its orbit except along the strands of $\Tt$.

A tangle being equivariantly boundary-parallel is a priori stronger than being invariant and boundary-parallel, but we have the following lemma which shows that (here in dimension $3$) these conditions are the same.

\begin{lemma}
\label{lem:equiv_parallel}
	Let $G$ act on a $3$--dimensional $1$--handlebody $H$.
	A $G$--invariant tangle $\Tt\subseteq H$ is equivariantly boundary-parallel if and only if it is boundary-parallel.
\end{lemma}

\begin{proof}
	The forward direction is clear, so we proceed with the reverse direction.
	The proof follows essentially from Edmonds' formulation of the~\ref{ELT} and relies critically on the fact that Edmonds accounts for orientation-reversing actions (cf.~\cite[Lemma~3.3]{Edm_86_A-topological-proof}).

	Let $Y = H\cup_{\partial H} \overline{H}$ be the double of $H$ across its boundary, which carries a natural doubled $G$--action.
	The double $Y$ also carries a natural orientation-reversing $\Z_2$--action (given by swapping the factors of the double), whose fixed set is $\partial H$.
	These actions commute and together give a $(G\times \Z_2)$--action on $Y$.

	For each strand $\tau$ of $\Tt$, let $c_\tau = \tau\cup \overline\tau$, and let $C$ denote the union of these $c_\tau$.
	Given any collection of bridge-disks for $\Tt$, the double of these disks across their arcs of intersection with $\partial H$ is a collection of spanning disks for $C$.
	It follows that $C$ is an unlink.
	
	Consider the induced action on $\nu(c_\tau)$.
	This is an action of $G_{\nu(\tau)}\times \Z_2$ on $\nu(c_\tau)$ induced by the double, where $G_{\nu(\tau)}$ acts on $\nu(\tau)$ and $\nu(\overline\tau)$ (both copies of $D^2\times I$) as a subgroup of an orientation-preserving and product-preserving dihedral action, and $\Z_2$ acts orientation-reversingly by swapping $\nu(\tau)$ and $\nu(\overline\tau)$ by the natural reflection given by the double.
	The group $G_{\nu(\tau)}$ is either cyclic or dihedral, with the later case occurring when an element of $G$ interchanges the points of $\partial\tau$.
	
	If $G_{\nu(\tau)}$ is cyclic, then the singular set of the $(G_{\nu(\tau)}\times \Z_2)$--action on $\partial(\nu(c_\tau))$ is the two meridians $\partial(\nu(c_\tau))\cap\partial H$; these meridians are fixed by the $\Z_2$--action, and rotated by the $G_{\nu(\tau)}$--action.
	If $G_{\nu(\tau)}$ contains an element interchanging the points of $\partial\tau$, then this singular set also includes two meridional rings of evenly spaced points; these are the points fixed by the conjugates of the reflection in $G_{\nu(\tau)}$ when it is dihedral.
	
	In any event, we can choose an equivariant longitude curve $\ell_\tau\subset \partial(\nu(c_\tau))$ that is disjoint from these singular points and transverse to the $\Z_2$--singular meridians.
	Then, there is an almost-equivariant annulus $A_\tau$  describing an isotopy from $\ell_\tau$ to $c_\tau$ in the sense that the orbit annuli $G_\tau\cdot A_\tau$ are disjoint except for all intersecting simultaneously at $c_\tau$.
	Choose one such curve $\ell_\tau$ for each orbit of strands in $\Tt$.

	Now consider $Y' = Y\setminus \nu(C)$, where an open $(G\times\Z_2)$--invariant tubular neighborhood has been removed.
	Since $C$ is an unlink, the $\ell_\tau$ chosen above can be assumed to be null-homotopic in $Y'$.
	Each longitude is transverse to the singular set on $\partial Y'$ by construction, hence it meets the hypotheses of the~\ref{ELT}.
	Pick one $\ell_\tau$, and apply the theorem to obtain a $(G\times\Z_2)$--equivariant spanning disk $D_\tau$ for $\ell_\tau$.
	Remove an invariant tubular neighborhood of the orbit of $D_\tau$, and repeat with some remaining $\ell_\tau$.
	Continuing, we obtain an equivariant collection on spanning disks for the $\ell_\tau$.
	It remains to show that each $D_\tau\cup A_\tau$ is the double of an equivariant bridge-disk $\delta_\tau$ for $\tau$, at which point the union of orbits of these disks across all $\tau$ is an invariant almost-disjoint bridge-disk system, as desired.

	To start, we claim that $D_\tau$ is transverse to $\partial H$.
	First note that $D_\tau$ is invariant under the $\Z_2$--action, since $D_\tau$ is equivariant and intersects $\partial H$, the fixed-point set of the $\Z_2$--action.
	Let $\nu(D_\tau)$ be a $\Z_2$--invariant tubular neighborhood of $D_\tau$.
	Then, $\nu(D_\tau)$ is equivariantly diffeomorphic to $D_\tau\times [0,1]$, equipped with a product-preserving action of $\Z_2$.
	Since the action of $\Z_2$ on $D_\tau$ is orientation-reversing, and the action of $\Z_2$ is orientation-reversing on $\nu(D_\tau)$, the action on the $[0,1]$ factor of the product must be orientation-preserving and hence trivial.
	Let $F$ be the fixed-point set of the $\Z_2$--action on $D_\tau$, which is a neatly embedded, $1$--dimensional submanifold of $D_\tau$.
	Since the action of $\Z_2$ on $\nu(D_\tau) = D_\tau\times[0,1]$ is trivial on the second factor, the fixed-point set of $\Z_2$ in $\nu(D_\tau)$, i.e. the intersection $\partial H\cap \nu(D_\tau)$, must be $F\times[0,1]$, hence $\partial H$ intersects $D_\tau$ transversely.

	Since $D_\tau$ is transverse to $\partial H$, we have $D_\tau\cap \partial H$ contains a single arc $a_\tau$, which connects the two points of $\ell_\tau\cap\partial H$, and some number of simple closed curves.
	However, an innermost such simple closed curve would bound a disk to one side of $\partial H$ that would be reflected across $\partial H$ by the $Z_2$--action, implying that $D_\tau$ contains a 2--sphere component, a contradiction.
	It follows that $D_\tau\cap\partial H = a_\tau$, and $D_\tau$ is the double of a bridge-disk $\delta_\tau$ whose boundary is $\tau\cup a_\tau$, as claimed.
\end{proof}

Our final concern in adapting boundary-parallelism to the equivariant setting involves how the linearly parting disk-system should interact with equivariantly boundary-parallel tangles, and is addressed by the following definition.

\begin{definition}
\label{def:linearly_parted_pair}
	Let $G$ act on a $3$--dimensional $1$--handlebody $H$, and let $\mathcal{T}$ be an equivariantly boundary-parallel tangle in $H$.
	We say that the $G$--pair $(H, \mathcal{T})$ is \emph{linearly parted as a pair} if there exists a linearly parting disk-system $\Bb$ for $H$ and an invariant, almost-disjoint collection $\Delta$ of bridge-disks for $\Tt$ such that $\Bb\cap\Delta = \varnothing$ and each $3$--ball component of $H\setminus\Bb$ contains at most one strand of $\Tt$.
\end{definition}

Again, this is \emph{a priori} stronger than simply being equivariantly boundary-parallel in $H$, but again, we have a result showing showing this is not the case.

\begin{corollary}
\label{cor:one_arc}
	If a tangle $\Tt\subset H$ is equivariantly boundary-parallel, then $(H,\Tt)$ is linearly parted as a pair.
\end{corollary}

\begin{proof}
	Let $\Delta$ be an invariant, almost-disjoint collection of bridge-disks for $\Tt$.
	Let
	$$D = \partial(\nu(\Delta))\setminus\partial H.$$ 
	Then $D$ is an invariant, disjoint collection of neatly embedded disks in $H$ with the property that $H\setminus\nu(D)$ is the disjoint union of some handlebodies (collectively denoted $H'$), together with some 3--balls, each of which contains precisely one strand of $\Tt$.
	Note that each disk $d\in D$ contributes one disk to $\partial H'$ by construction, which in abuse of notation we will also call $d$.

	A consequence of the~\ref{ELT} and the Linearization Theorem (see~\cite[Theorem~3.5]{MeiSco_LP}) is that every finite group action on a $3$--dimensional $1$--handlebody is linearly parted; see~\cite[Corollary~3.6]{MeiSco_LP}.
	So, the induced action on the handlebody $H'$ has a linearly parting disk-system $\Bb'$. 
	
	If the $G$--action on some $b\in\Bb'$ is orientation-reversing on $b$, then replace $b$ with $\partial(\nu (b))\setminus\partial H'$ and extend this replacement equivariantly along the orbit of $b$. 
	It is clear that $\Bb'$ is still a linearly parting ball-system after this replacement, since we have just subdivided one $1$--handle in the handle decomposition into two $1$--handles and a $0$--handle.
	Thus we can assume the induced action of $G$ on each $b\in\Bb'$ is orientation-preserving.

	Suppose the boundary of a disk $b\in \Bb'$ intersects a disk $d\in D$.
	Since $d\subseteq \partial H'$ and $G$ is orientation-preserving on $H'$, the stabilizer $G_d$ is trivial or cyclic, acting by rotation.
	If the action is trivial, then $b$ can be $G_d$--equivariantly isotoped off of $d$, hence the orbit of this isotopy gives an equivariant isotopy of the orbit of $b$ off the orbit of $d$.
	If the action is cyclic, then let $c\in d$ be the fixed-point center of $G_d$.
	Since $b$ is equivariant, and $G$ acts on $b$ orientation-preservingly, $\partial b$ does not intersect $c$.
	Thus since $d$ is a $G_d$--equivariant cone, $\partial b$ can be moved via a straight-line equivariant isotopy away from the cone point and off $d$.
	(Equivariant cones are discussed in detail in~\cite[Section~2.1]{MeiSco_LP}.)

	After completing these isotopies, the disk-system $\Bb'$ can be combined with $D$ to give the desired linearly parting disk-system $\Bb$ for $H$, since $\Bb\cap\Delta = \varnothing$ and the disks of $\Bb$ isolate the strands of $\Tt$, by construction.
\end{proof}

We now extend the notion of equivariant boundary-parallelism to the setting of trivial disk-tangles in $4$--dimensional $1$--handlebodies.

\begin{definition}
	Let $X$ be a 4--dimensional 1--handlebody.
	A \emph{$p$--patch disk-tangle} $\Dd$ is a disjoint collection of neatly embedded disks in $X$.
	A disk-tangle $\Dd$ is \emph{trivial} (or \emph{boundary-parallel}) if $\Dd$ can be isotoped relative to $\partial\Dd$ to lie in $\partial X$.
	Equivalently, if $\Dd = \{P_i\}_{i=1}^p$, there are pairwise disjoint $3$--balls $\Delta = \cup_{i=1}^p\delta_i$ such that $\partial\delta_i = P_i\cup D_i$, where $D_i$ is a spanning disk for $\partial P_i$ in $\partial X$.
	The balls of $\Delta$ are collectively called \emph{bridge-balls} for the disk-tangle $\Dd$.
\end{definition}

\begin{definition}
	Let $G$ act on  a $4$--dimensional $1$--handlebody $X$.
	An invariant trivial disk-tangle $\Dd = \cup_iP_i$  in $X$ is \emph{equivariantly boundary-parallel} if there is a collection $\Delta = \cup_i\delta_i$ of bridge-balls that is $G$--invariant and \emph{almost-disjoint}: each ball $\delta_i$ is disjoint except for possible intersections on their boundaries along patches $P_i$.
\end{definition}

We now have the following definition, analogous to Definition~\ref{def:linearly_parted_pair} above.

\begin{definition}
	Let $G$ act on a $4$--dimensional $1$--handlebody $X$, and let $\Dd$ be an equivariantly boundary-parallel disk-tangle in $X$.
	We say that the $G$--pair $(X, \Dd)$ is \emph{linearly parted as a pair} if there exists a linearly parting ball-system $\Bb$ for $X$ and a collection $\Delta$ of invariant, almost-disjoint bridge-balls for $\Dd$ such that $\Bb\cap\Delta = \varnothing$ and each $4$--ball component of $X\setminus\Bb$ contains at most one disk of $\Dd$.
\end{definition}

In contrast to the $3$--dimensional setting, it is not clear that every invariant disk-tangle should be equivariantly boundary-parallel, nor is it clear that every action on a disk-tangle should be linearly parted as a pair; see~\cite[Section~5]{MeiSco_LP}.

The following lemma is a nice criterion for equivariant boundary-parallelism, and is important for our existence proof in Section~\ref{sec:existence}.

\begin{lemma}
\label{lem:boundary_parallel_in_ball}
	Let $G$ act linearly on $B^4$.
	Let $D\subset B^4$ be a $G$--invariant, neatly embedded disk.
	If $D$ has a single critical point with respect to the radial Morse function on $B^4$, then $D$ is equivariantly boundary-parallel in $B^4$.
	Moreover, there exists a $G$--almost-equivariant bridge-ball for $D$ that misses any prescribed invariant collection of disjoint 3--balls in $\partial B^4$.
\end{lemma}

\begin{proof}
	Let $h\colon B^4\to[0,1]$ be the radial Morse function.
	For any $t\in[0,1]$, let $Y_t = h^{-1}(t)$.
	For any interval $t\in I\subseteq[0,1]$ and any compact subset $W\subseteq Y_t$, let $W_I$ denote the cylinder obtained by flowing $W\cap Y_t$ along the gradient of $h$ through the interval $I$.
	(We allow one-point intervals.)

	Let $c\in D$ be the critical point, which must be a minimum.
	Then, $D\cap h^{-1}(h(c)+\epsilon)$ is an unknot in the 3--sphere $Y_{h(c)+\epsilon}$.
	Since $h$ has no critical values greater than $h(c)$, the pairs $(Y_1,Y_1\cap D)$ and $(Y_{h(c)+\epsilon},Y_{h(c)+\epsilon}\cap D)$ are $G$--equivariantly diffeomorphic~\cite[Theorem~4.3]{Was_69_Equivariant-differential-topology}.
	It follows that $\partial D$ is a $G$--invariant unknot.
	By the~\ref{ELT}, $\partial D$ has a $G$--almost-equivariant Seifert disk $F$.
	Note that $F$ can be chosen to miss any prescribed invariant collection of disjoint 3--balls in $\partial B^4$.
	
	Let $\delta = F_{[h(c),1]}$.
	Then $\delta$ is a $G$--almost-equivariant bridge-ball for $D$, as desired.
\end{proof}

\subsection{Spines and diagrams}
\label{subsec:spines}

An important feature of equivariant (bridge) trisections is that they are uniquely determined by their spines.
The key ingredients in the proofs of these facts are the equivariant versions of the Laudenbach-Po\'enaru theorems proved in~\cite{MeiSco_LP}.

\begin{theorem}
\label{thm:spines}
	Let $\mathfrak T$ and $\mathfrak T'$ be $G$--equivariant trisections of $X$ and $X'$, respectively.
	If the restrictions of the actions to the spines of the trisections are $G$--diffeomorphic, then the trisections are $G$--diffeomorphic.
\end{theorem}

\begin{proof}
	Let $\psi\colon H_1\cup H_2\cup H_3 \to H_1'\cup H_2'\cup H_3'$ be the given $G$--equivariant diffeomorphism of the spines of the trisections $\TT$ and $\TT'$.
	Use this diffeomorphism to identify $Y_i = H_i\cup_\Sigma\overline H_{i+1}$ with $Y_i' = H_i'\cup_{\Sigma'}\overline H_{i+1}'$.
	Then $X_i$ and $X_i'$ are linearly parted fillings of $Y_i$ extending the same $G$--action.
	By~\cite[Theorem~4.1(2)]{MeiSco_LP}, the $G$--diffeomorphism $\psi$ extends to a $G$--diffeomorphism $\Psi_i\colon X_i\to X_i'$.
	Therefore, $\Psi = \Psi_1\cup\Psi_2\cup\Psi_3$ is a $G$--diffeomorphism of $\TT$ and $\TT'$.
\end{proof}

\begin{theorem}
\label{thm:spines_bridge}
	Let $\mathfrak T$ and $\mathfrak T'$ be $G$--equivariant bridge trisections for $(X,S)$ and $(X',S')$, respectively.
	If the restrictions of the actions to the spines of the bridge trisections are $G$--diffeomorphic, then the bridge trisections are $G$--diffeomorphic.
\end{theorem}

\begin{proof}
	Let $\psi\colon H_1\cup H_2\cup H_3 \to H_1'\cup H_2'\cup H_3'$ be the given $G$--equivariant diffeomorphism of the spines of the bridge trisections $\TT$ and $\TT'$ taking $\mathcal T_i = H_i\cap S$ to $\mathcal T_i' = H_i'\cap S'$.
	Use this diffeomorphism to identify $Y_i = H_i\cup_\Sigma\overline H_{i+1}$ with $Y_i' = H_i'\cup_{\Sigma'}\overline H_{i+1}'$ and $L_i = \mathcal T_i\cup \overline{\mathcal T}_{i+1}\subseteq Y_i$ with $L_i' = \mathcal T_i'\cup \overline{\mathcal T}'_{i+1}\subseteq Y_i'$.
	Then $(X_i, \Dd_i)$ and $(X_i', \Dd_i')$ are linearly parted fillings of $(Y_i, L_i)$ extending the same $G$--action.
	By~\cite[Theorem~5.9(2)]{MeiSco_LP}, the $G$--diffeomorphism $\psi$ extends to a $G$--diffeomorphism $\Psi_i\colon (X_i, \Dd_i)\to (X_i', \Dd_i')$.
	Therefore, $\Psi = \Psi_1\cup\Psi_2\cup\Psi_3$ is a $G$--diffeomorphism of the bridge trisections $\TT$ and $\TT'$.
\end{proof}

We conclude this subsection with a corollary indicating that the fixed-point set of an action must be in bridge trisected position whenever it intersects the spine of the trisection in the spine of a bridge trisected surface; cf. the examples of Subsection~\ref{subsec:marla_equivariant_trisections}.

\begin{corollary}
\label{cor:bridge_in_spine}
	Let $\mathfrak T$ be a $G$--equivariant trisection of $X$, and let $F_i$ be the fixed-point set of $G$ in $H_i$.
	If $F_1\cup F_2\cup F_3$ is the spine of an equivariant bridge trisection, then the $2$--dimensional part of the fixed-point set of $G$ is in equivariant bridge trisected position with respect to $\mathfrak T$.
\end{corollary}

\begin{proof}
	The $2$--dimensional part of the fixed point set in each $X_i$ is a neatly embedded invariant surface since the action on $X_i$ is linearly parted, and the boundary of this fixed-point set must be $F_i\cup \overline F_{i+1}$, an unlink.
	By~\cite[Corollary~5.6]{MeiSco_LP}, this surface must be trivial disk-tangle such that $(X_i,\Dd_i)$ is linearly parted as a pair as desired.
\end{proof}


The above results establish that equivariant (bridge) trisections are determined by their spines.
However, each piece of the spine can be encoded diagrammatically on the core surface.
This leads to the notion of equivariant shadow diagrams for equivariant (bridge) trisections, which we now develop.

\begin{definition}
\label{def:shadow}
	Let $G$ act on a surface $\Sigma$.	
	A tuple $(\Sigma;\alpha_1,\alpha_2,\alpha_3)$ is a \emph{$G$--invariant trisection diagram} if each $\alpha_i$ is a $G$--invariant cut-system for $\Sigma$ such that each triple $(\Sigma;\alpha_i,\alpha_{i+1})$ determines a Heegaard splitting of $\#^{k_i}(S^1\times S^2)$ for some $k_i\geq 0$.
	
	The arcs in a collection $\frak a$ of arcs in $\Sigma\setminus\nu(\alpha)$ are called \emph{$G$--invariant shadow arcs} with respect to $\alpha$ if
	\begin{enumerate}
		\item $\frak a$ is $G$--invariant,
		\item the arcs of $\frak a$ have disjoint interiors; and
		\item each connected component of $\Sigma\setminus\nu(\alpha)$ contains at most one connected component of $\frak a$.
	\end{enumerate}
	In abuse of notation, we write $\partial(\frak a)$ to mean the collection of points in $\Sigma$ occurring as a boundary point of some arc of $\frak a$, noting that shadow arcs in the same orbit can share endpoints.
	
	A tuple
	$$(\Sigma;\alpha_1,\alpha_2,\alpha_3,\frak a_1,\frak a_2,\frak a_3,\bm x)$$ is a \emph{$G$--invariant shadow diagram} if $(\Sigma;\alpha_1,\alpha_2,\alpha_3)$ is a $G$--invariant trisection diagram and $\frak a_i$ is a collection of $G$--invariant shadow arcs with respect to $\alpha_i$ with $\partial(\frak a_i) = \bm x$ such that
	$$(\Sigma;\alpha_i,\alpha_{i+1},\frak a_i,\frak a_{i+1})$$
	determines a generalized bridge splitting of a pair $\left(\#^{k_i}(S^1\times S^2),\sqcup_{p_i}U\right)$, where $\sqcup_{p_i}U$ is a $p_i$--component unlink, in the sense of~\cite[Section~2.A]{MeiZup_18_Bridge-trisections}.
\end{definition}

For ease of exposition, we will sometimes refer to a trisection diagram as a shadow diagram for which $\frak a_i =\varnothing$.
A $G$--invariant shadow diagram generally contains a great deal of redundant information.
For example, the $\Z_n$--action on a solid torus that rotates in the direction of the core curve would have a $\Z_n$--invariant cut-system consisting of $n$ parallel curves.
In practice, recording a subset of these curves whose orbit is the entire invariant set is sufficient; similar remarks hold for shadow arcs.

\begin{proposition}
\label{prop:diag_spine}
	Fix a finite group $G$.
	\begin{enumerate}
		\item A $G$--invariant shadow diagram determines a unique $G$--equivariant (bridge) trisection, up to $G$--diffeomorphism.
		\item $G$--diffeomorphic $G$--invariant shadow diagrams determine $G$--diffeomorphic $G$--equivariant (bridge) trisections.
		\item Every $G$--equivariant (bridge) trisection can be represented by a $G$--invariant shadow diagram.
	\end{enumerate}
\end{proposition}

\begin{proof}
	Consider a $G$--invariant shadow diagram
	$$\DD = (\Sigma;\alpha_1,\alpha_2,\alpha_3,\frak a_1,\frak a_2,\frak a_3,\bm x).$$
	First, we claim that the triple $(\Sigma;\alpha_i,\frak a_i)$ uniquely determines a $G$--pair $(H_i,\Tt_i)$ where $H_i$ is a $G$--equivariant 3--dimensional 1--handlebody with $\partial H = \Sigma$, $\Tt_i$ is a $G$--equivariantly boundary-parallel trivial tangle with $\partial\Tt_i = \bm x$, and $(H_i,\Tt_i)$ is linearly parted as a pair by $\Bb_i\cup\Delta_i$.
	We build $H_i$ by taking $\Sigma\times I$ and attaching 3--dimensional 2--handles along $\alpha_i\times\{1\}$, then capping off the resulting 2--sphere boundary components with 3--dimensional 3--handles.
	Since $\alpha_i$ is $G$--invariant, we can attach the 2--handles along a $G$--invariant tubular neighborhood of $\alpha_i$.
	
	We can extend the $G$--action on $\Sigma\times I$ (the product action) across each 2--handle in a unique way: Given a 2--handle $\frak h$, we have $\partial \frak h = A\cup d\cup d'$, where $A$ is the attaching region of the handle (an annulus) and $d$ and $d'$ are disks.
	The induced $G$--action on $A$ is predetermined, is a subgroup of a dihedral group, and has a unique extension across $d$ and $d'$.
	(The fact that every action on $S^1$ extends uniquely across $D^2$ is a simple application of the Smooth Palais' Theorem for lifting isotopies~\cite[Corollary~2.4]{Sch_80_Lifting-smooth-homotopies}.)
	The resulting action on $\partial h$ has a unique extension across $\frak h$, by the equivariant version of Munkres' theorem~\cite[Corollary~2.6]{DinLee_09_Equivariant_Ricci_flow}.
	The action is now determined on the attaching regions of the 3--handles, across which it extends uniquely for the same reason.
	This gives a unique extension of the $G$--action on $\Sigma$ to $H_i$.
	Let $\Bb_i$ denote the union of the cores of the 2--handles.
	So, $\Bb_i$ is a linearly parting disk-system for $H$ with $\partial\Bb_i =\alpha_i$.
	
	We remark in passing that if $\DD$ is a diagram for an equivariant trisection (i.e., $\frak a_i=\varnothing$), then the proof of part (1) is complete.
	
	The condition in Definition~\ref{def:shadow} that each component of $\Sigma\setminus\nu(\alpha_i)$ contain at most one connected component of $\frak a_i$ allows us to extend $\frak a_i$ to a well-defined collection of $G$--almost-equivariant bridge-disks $\Delta_i$, as follows.
	Let $a$ be a shadow arc of $\frak a_i$, and suppose that $a\times\{1\}$ (which we will henceforth call $a$) lies in $\partial Z$ for some 3--handle $Z$ of $H_i\setminus\Bb_i$.
	Consider the induced action $G_Z$ on $Z$, and note that $a$ is $G_Z$--equivariant.
	Since $a$ is the only shadow arc in $\partial Z$, the $G_Z$--action must fix $\partial a$ set-wise.
	Regarding the $G_Z$--action as an orthogonal action on $B^3$ (by the Linearization Theorem; cf.~\cite[Theorem~3.5]{MeiSco_LP}), we can assume $\partial a$ is a pair of antipodal points of $\partial Z$ and $a$ is a longitudinal arc.
	Since the action of $G$ on $\Sigma$ is orientation preserving, the same is true of the action of $G_Z$ on $Z$.
	It follows that $G_Z$ is a subgroup of a dihedral group.
	Let $t$ be a neatly embedded arc in $Z$ connecting the points of $\partial a$, and let $\delta$ be a $G_Z$--almost-equivariant bridge-disk for $t$ such that $\partial \delta = t\cup a$; this is easy to arrange by the linearity of the $G_Z$--action.
	
	Let $\Tt_i$ (respectively, $\Delta_i$) be the union of all arcs $t$ (respectively, disks $\delta$) constructed as above, ranging over the arcs $a$ of $\frak a$.
	By construction, $\Delta_i\cap\Bb_i = \varnothing$.
	Thus, the $G$--action on $(H_i,\Tt_i)$ is linearly parted as a pair, as desired.
	This gives a well-defined extension of the $G$--action on the diagram $\DD$ to the spine
	$$(H_1,\Tt_1)\cup(H_2,\Tt_2)\cup(H_3,\Tt_3).$$
	This spine determines a well-defined $G$--equivariant (bridge) trisection by Theorems~\ref{thm:spines} and~\ref{thm:spines_bridge}.	

	For part (2), suppose $\DD$ and $\DD'$ are $G$--equivariant shadow diagrams for  $G$--equivariant (bridge) trisections $\TT$ and $\TT'$.
	Let $\psi\colon\DD\to\DD'$ be a $G$--diffeomorphism so that $\psi(\alpha_i) = \alpha_i'$, and so on.
	Since every $G$--diffeomorphism of a circle extends to a $G$--diffeomorphism of disks, we can extend $\psi$ across the cut-disks $D_i$ bounded by $\alpha_i$.
	This determines $\psi$ on $\Sigma\cup D_1\cup D_2\cup D_2$.
	The cut-disks $D_i$ cut $H_i$ into 3--balls, each of which contains at most one arc of of $\Tt_i$.
	Since every $G$--diffeomorphism of $(S^2,\{N,S\})$ extends to a $G$--diffeomorphism of the trivial ball/arc-pair, we can extend $\psi$ to all of $(H_i,\Tt_i)$.
	Now that $\psi$ is defined on the spine of $\TT$, it extends to all of $\TT$ by Theorem~\ref{thm:spines}.

	For part (3), suppose that $(X,\Ss,\TT,G)$ is a $G$--equivariant (bridge) trisection (allowing for the possibility that $\Ss = \varnothing$).
	Let $(H_1,\Tt_1)\cup(H_2,\Tt_2)\cup(H_3,\Tt_3)$ be the spine of $\TT$, so, for each $i\in\Z_3$, we have that $H_i$ is a $G$--equivariant 3--dimensional 1--handlebody, $\Tt_i$ is a $G$--equivariantly boundary-parallel trivial tangle, and $(H_i,\Tt_i)$ is linearly parted as a pair by $\Bb_i\cup\Delta_i$, as in Definition~\ref{def:linearly_parted_pair}, so $\Bb_i\cap\Delta_i = \varnothing$.
	
	Let $\frak a_i = \Delta_i\cap\Sigma$, and let $\alpha_i = \partial\Bb_i$.
	Then, the tuple
	$$(\Sigma;\alpha_1,\alpha_2,\alpha_3,\frak a_1,\frak a_2,\frak a_3,\bm x)$$
	is a $G$--invariant shadow diagram for $\TT$.
	By part (1), this diagram can be used to recover $\TT$, establishing part (2).
\end{proof}

One useful application of this proposition is illustrated the examples of Subsection~\ref{subsec:marla_equivariant_trisections}, where $G$--actions on the underlying manifolds are ``discovered'' by noticing $G$--actions on a trisection diagrams.

\section{Existence of equivariant trisections}
\label{sec:existence}

In this section, we prove our main existence result.

\begin{theorem}
\label{thm:exist_main}
	Let $G$ be a finite group, let $X$ be a smooth, orientable, connected, closed, $4$--dimensional $G$--manifold, and let $\{\Ss_i\}_{i=1}^N$ be a (possibly empty) collection of embedded, closed, $G$--invariant surfaces in $X$.
	Then, $X$ admits a $G$--equivariant trisection for which each $\Ss_i$ is in equivariant bridge trisected position.
\end{theorem}

To prove this, we leverage a theorem of Illman~\cite{Ill_78_Smooth-equivariant-triangulations} that establishes the existence of equivariant triangulations.
The proof is completed in two stages.
First, an equivariant trisection of the underlying $G$--manifold $X$ is constructed from a given equivariant triangulation (Proposition~\ref{prop:exist_tri}).
Then, it is shown that, provided the triangulation restricts to a triangulation of a $G$--invariant surface $\Ss$, the resulting trisection bridge trisects $\Ss$ (Proposition~\ref{prop:exist_bridge}).

\subsection{Decompositions of equivariant triangulations}
\label{subsec:triangulation_to_trisection}

Given a finite group $G$, a simplicial complex $K$ is a \emph{$G$--simplicial complex} if $G$ acts on $K$ so that for all $g\in G$, the map $g\colon K\to K$ is a simplicial map.
Given a $G$--manifold $X$ and an equivariant triangulation $f\colon K\to X$, an invariant subset $W\subset X$ is said to be \emph{triangulated} by $f$ if there is a subcomplex $J\subset K$ such that $W = f(J)$.

\begin{theorem}[Illman,\cite{Ill_78_Smooth-equivariant-triangulations}]
\label{thm:Illman}
    Let $X$ be a smooth $G$--manifold for a finite group $G$.
    Let $S\subseteq X$ be an invariant surface.
    Then there exists a $G$--simplicial complex $K$ and a smooth equivariant triangulation $f\colon K\to X$ such that $S$ is triangulated by $f$.
\end{theorem}

Given this result, we construct an equivariant trisection on $X$ by decomposing each pentachoron (4--simplex) of an equivariant triangulation of $X$ into three pieces such that the union of these pieces across all the pentachora yields the desired trisection.
This key technical step is accomplished in the following proposition.

\begin{proposition}
\label{prop:pentachoron}
	Let $P$ be a pentachoron, and let the symmetric group $S_5$ act in the standard way on $P$.
	There exist $S_5$--invariant subsets $X_1, X_2, X_3\subset P$ such that
	$$P=X_1\cup X_2\cup X_3,$$
	and, for $i\in\Z_3$, we have:
	\begin{enumerate}
		\item $X_i$ is the neighborhood of a graph $\Gamma_i$ such that $\partial P\cap\Gamma_i$ is a graph $\partial\Gamma_i$;
		\item $H_i = X_i\cap X_{i-1}$ is $S_5$--invariant and is the neighborhood of a graph $\gamma_i$ such that $\partial P\cap\gamma_i$ is a graph $\partial\gamma_i$; and
		\item $\Sigma = X_1\cap X_2\cap X_3$ is an $S_5$--invariant, neatly embedded surface.
	\end{enumerate}
\end{proposition}

\begin{proof}

\begin{figure}[ht!]
     \centering
     \begin{subfigure}[b]{0.3\textwidth}
         \centering
         \includegraphics[width=\textwidth]{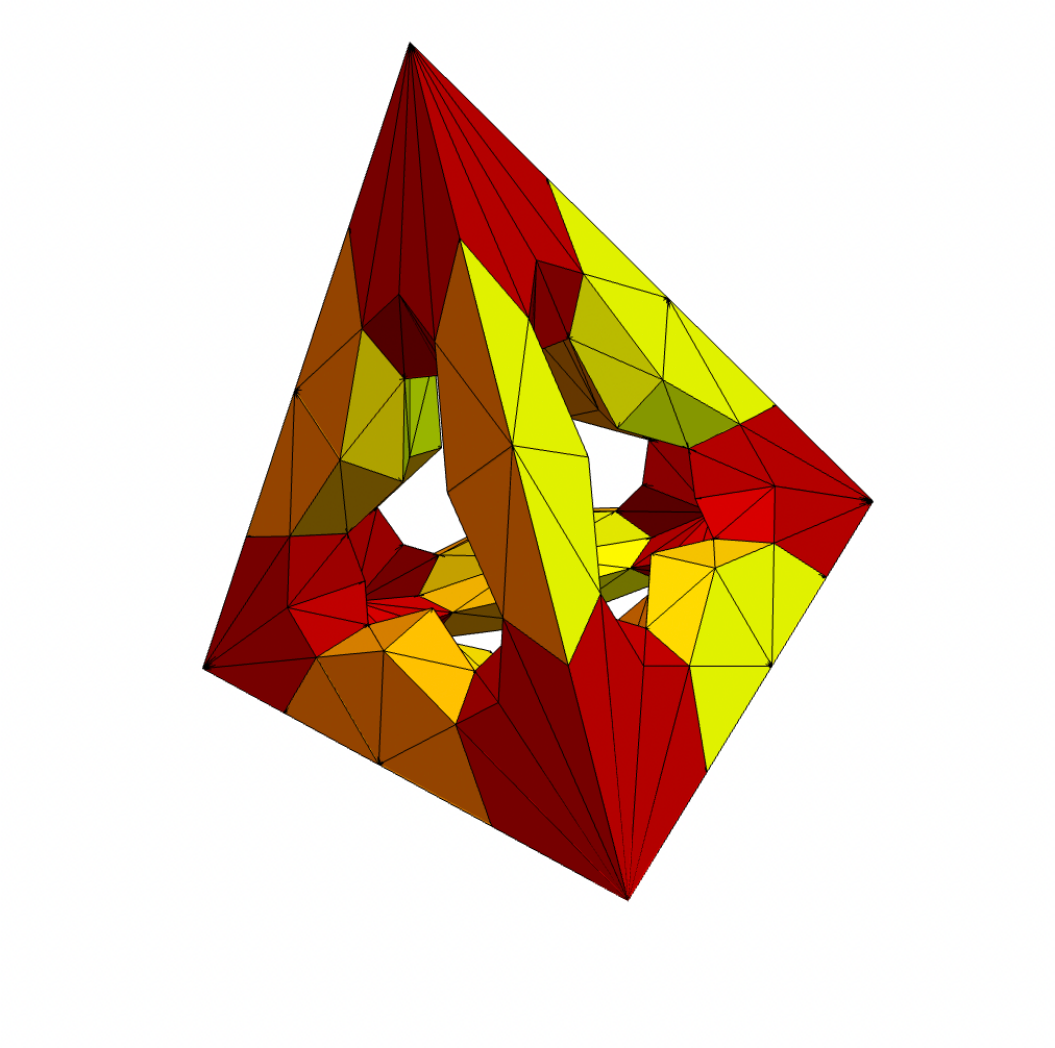}
         \caption{}
         \label{fig:vert_edge}
     \end{subfigure}
     \hfill
     \begin{subfigure}[b]{0.3\textwidth}
         \centering
         \includegraphics[width=\textwidth]{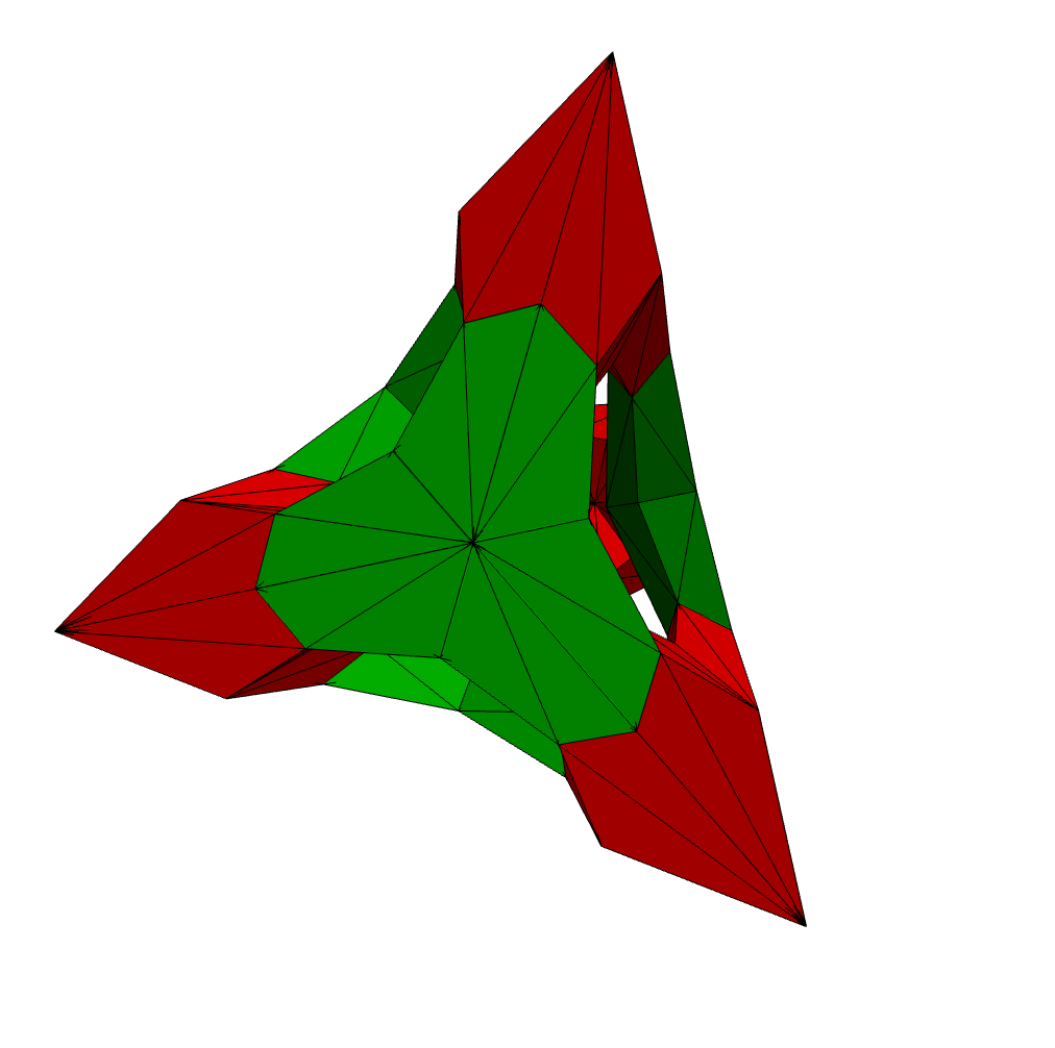}
         \caption{}
         \label{fig:vert_face}
     \end{subfigure}
     \hfill
     \begin{subfigure}[b]{0.3\textwidth}
         \centering
         \includegraphics[width=\textwidth]{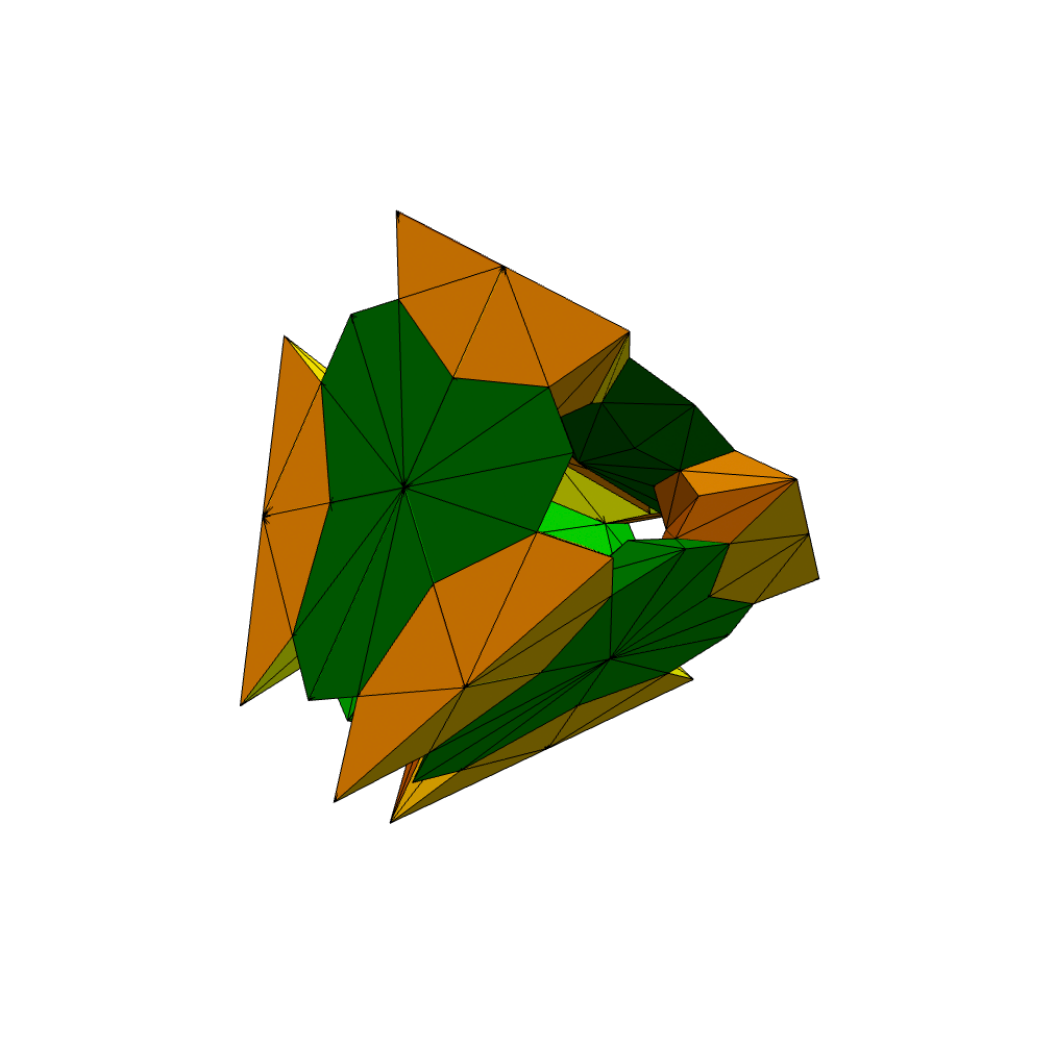}
         \caption{}
         \label{fig:face_edge}
     \end{subfigure}
          \begin{subfigure}[b]{0.3\textwidth}
         \centering
         \includegraphics[width=\textwidth]{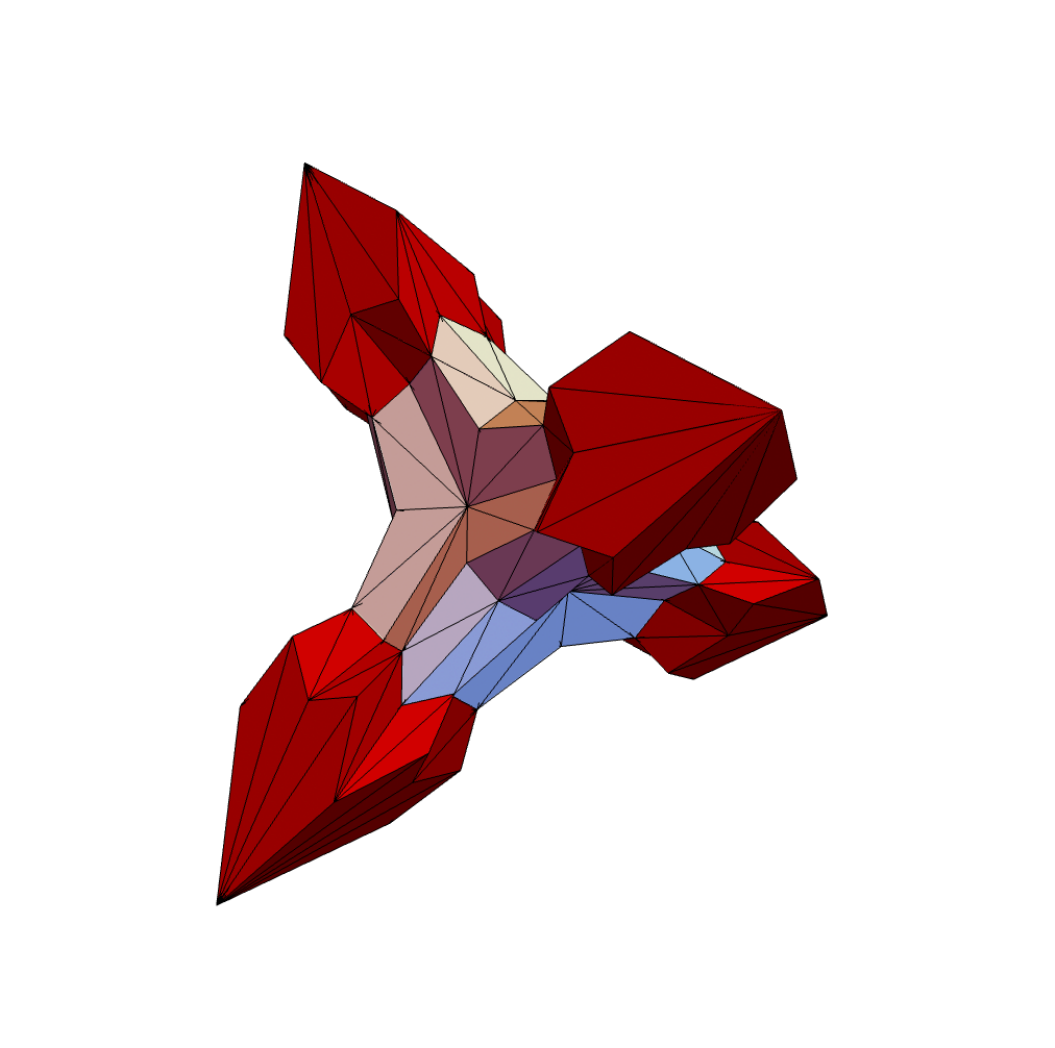}
         \caption{}
         \label{fig:core_vert}
     \end{subfigure}
     \hfill
     \begin{subfigure}[b]{0.3\textwidth}
         \centering
         \includegraphics[width=\textwidth]{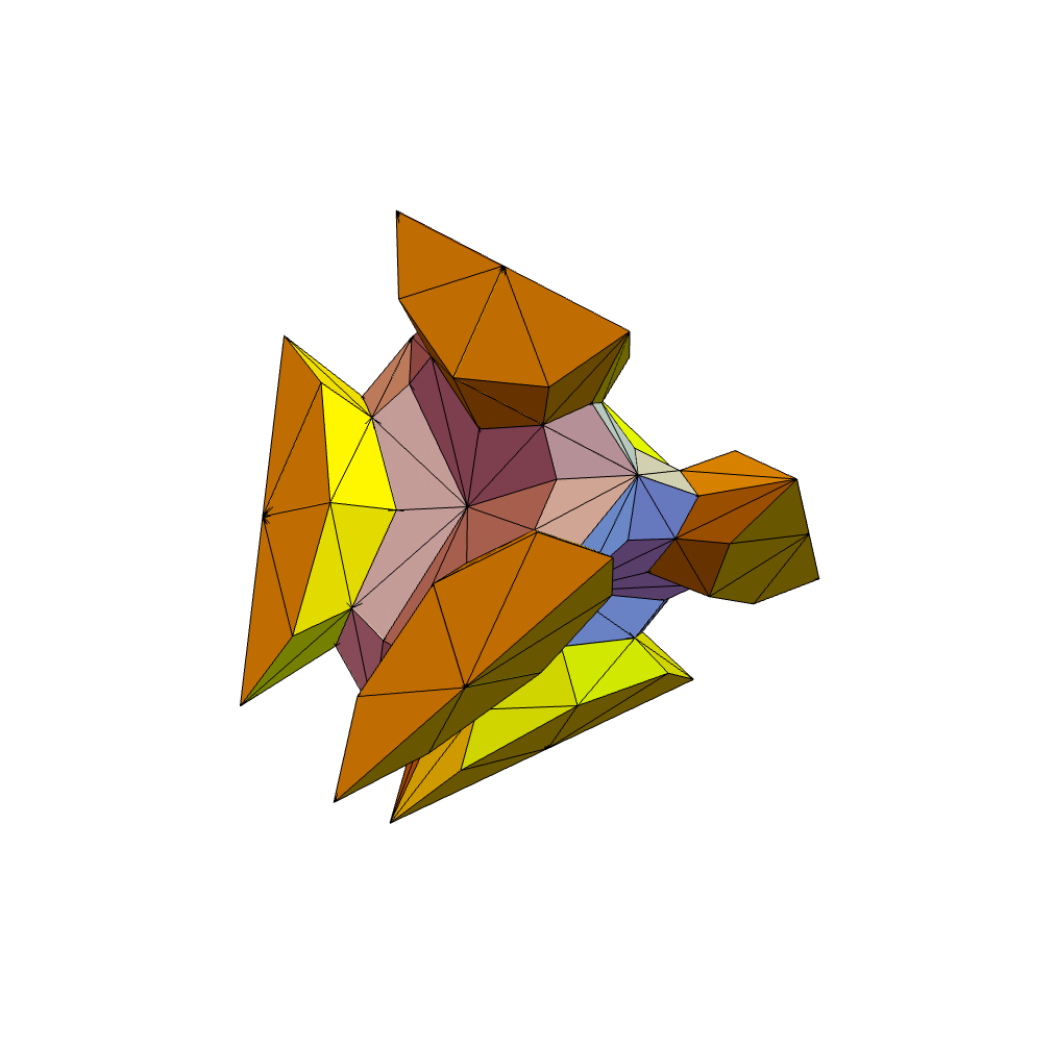}
         \caption{}
         \label{fig:core_edge}
     \end{subfigure}
     \hfill
     \begin{subfigure}[b]{0.3\textwidth}
         \centering
         \includegraphics[width=\textwidth]{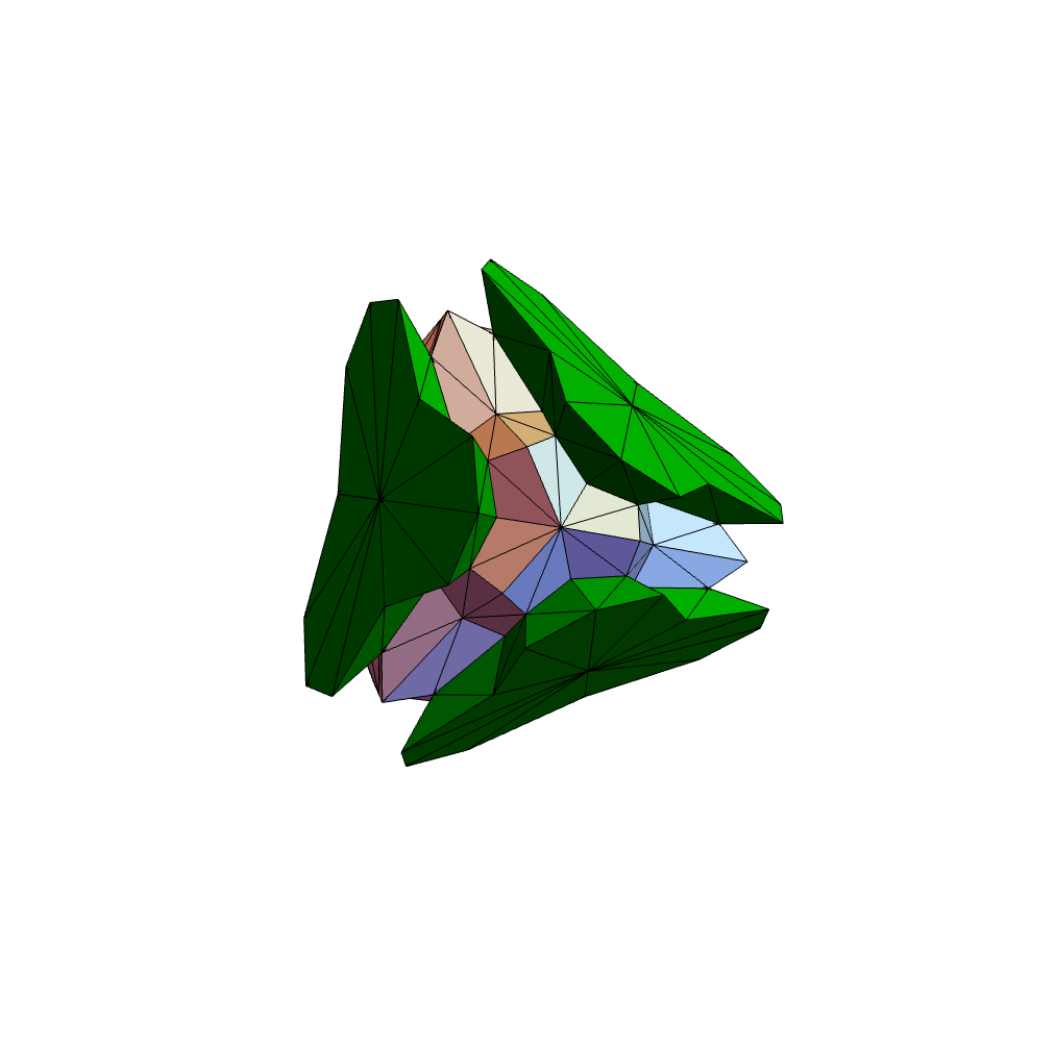}
         \caption{}
         \label{fig:core_face}
     \end{subfigure}
        \caption{Decomposition of tetrahedron into star neighborhoods using the double barycentric subdivision. These figures are derived from the wonderful Wolfram Demonstration by Aleksandr Berdnikov~\cite{Ber_24_Star-Neighborhoods-in-Double}.}
        \label{fig:tetrahedron}
\end{figure}

	Let $K$ denote the standard simplicial structure on $\partial P$, which consists of five tetrahedra $T$ that meet along ten triangular faces $F$, ten edges $E$, and five vertices $V$. 
	Let $K''$ denote the double barycentric subdivision of $K$.
	The restriction of $K''$ to a single tetrahedron of $\partial P$ is dissected in Figure~\ref{fig:tetrahedron}.
	Using $K''$, we decompose $\partial P$ into four pieces as follows:
	\begin{itemize}
		\item Let $\Vv$ denote the star neighborhood of the vertices $V$ of $K$; this is indicated in red in Figures~\ref{fig:vert_edge},~\ref{fig:vert_face}, and~\ref{fig:core_vert}.
		\item Let $\Ee$ denote the star neighborhood of the barycenters (midpoints) $B_E$ of the edges $E$ of $K$; this is indicated in yellow in Figures~\ref{fig:vert_edge},~\ref{fig:face_edge}, and~\ref{fig:core_edge}.
		\item Let $\Ff$ denote the star neighborhood of the barycenters $B_F$ of the triangular faces $F$ of $K$; this is indicated in green in Figures~\ref{fig:vert_face},~\ref{fig:face_edge}, and~\ref{fig:core_face}.
		\item Let $\Tt$ denote the star neighborhood of the barycenters $B_T$ of the tetrahedra $T$ of $K$; this is indicated in silver in Figures~\ref{fig:core_vert},~\ref{fig:core_edge}, and~\ref{fig:core_face}.
	\end{itemize}

	Regard $P$ as the cone on $\partial P$, with cone point $B_P$ the barycenter of $P$.
	Let $h\colon P\to[0,1]$ be the depth function: $h(B_P) = 1$ and $h(\partial P) = 0$.
	Let $Y_t = h^{-1}(t)$, so for each $t\in[0,1)$, $Y_t$ is a copy of $\partial P$.
	In what follows, for any compact subset $W\subset \partial P$ and any interval $I\subseteq[0,1]$, let $W_I$ denote the cylinder obtained by flowing $W$ along the gradient of $h$ through the interval $I$. (We allow one-point intervals.)
	We now define and analyze the pieces of the decomposition.
	See Figure~\ref{fig:depth_schem_X} below for a schematic of the decomposition for each $X_i$.

\begin{figure}[ht!]
	\centering
	\includegraphics[width = .9\linewidth]{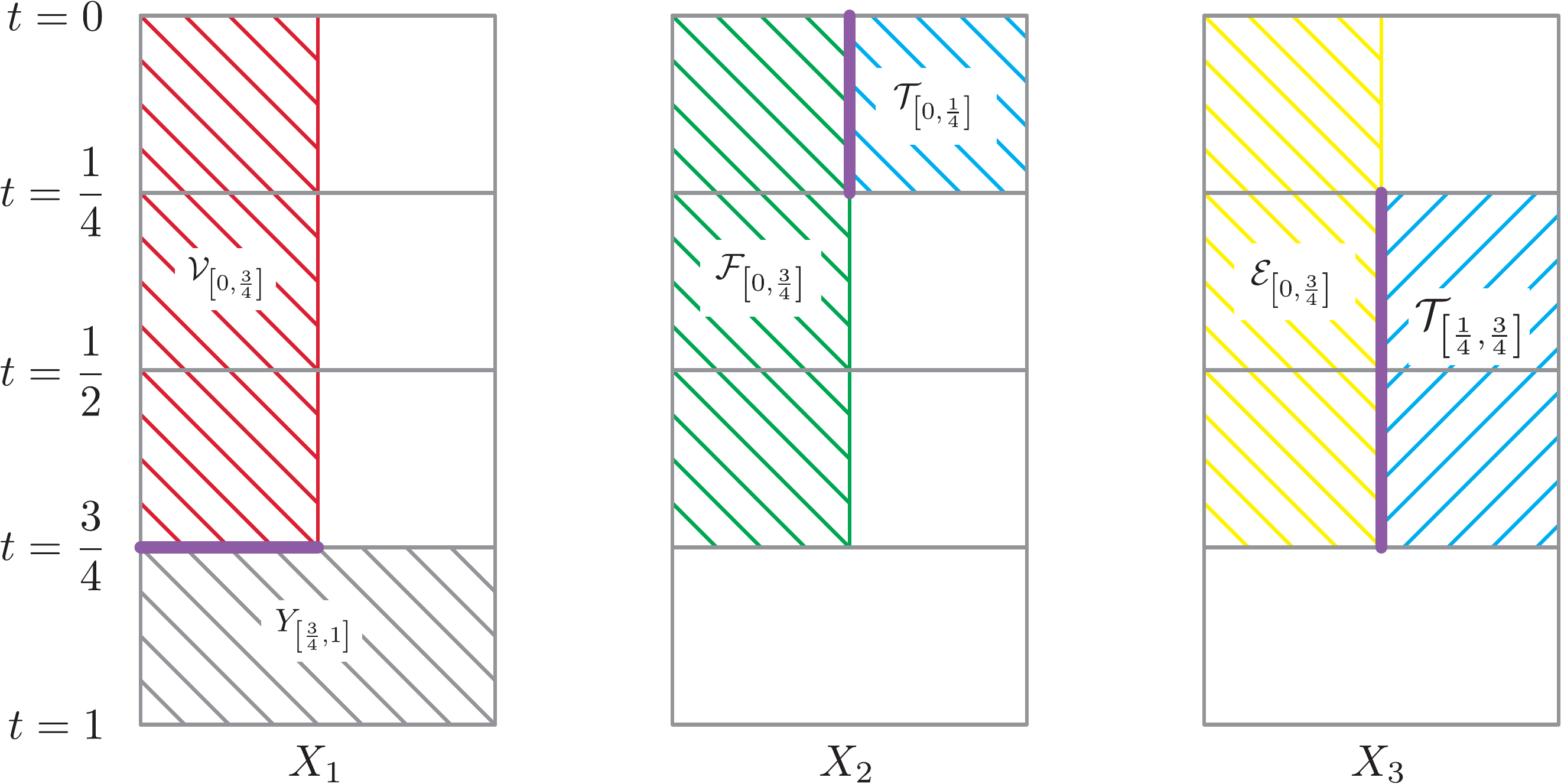}
	\caption{A schematic illustrating how the $X_i$ are decomposed based on the fundamental pieces $\Vv$, $\Ee$, $\Ff$, and $\Tt$ at each depth level.
	The purple segments indicate the location of the linearly parting ball-systems $\Bb_i$.}
	\label{fig:depth_schem_X}
\end{figure}
	
	Define
	$$X_1 = \Vv_{\left[0,\frac34\right]}\cup Y_{\left[\frac34,1\right]}
	\hspace{1cm} \text{and}\hspace{1cm}
	\Gamma_1 = V_{[0,1]}.$$
	Since $\Vv$ is the star neighborhood of $V$, it is clear that $X_1$ is a neighborhood of $\Gamma_1$ and is a 4--ball.
	Note that $\partial\Gamma_1 = V$, and the neighborhood $\Vv$ of $\partial\Gamma_1$ (restricted to a single tetrahedron) is shown in red throughout Figure~\ref{fig:tetrahedron}.
	
	Define
	$$X_2 = (\Ff\cup \Tt)_{\left[0,\frac14\right]}\cup \Ff_{\left[\frac14,\frac34\right]}
	\hspace{1cm} \text{and}\hspace{1cm}
	\Gamma_2 = \Gamma_2' \cup (B_F)_{\left[0,\frac12\right]},$$
	where $\Gamma_2'$ is the 1--skeleton of the dual triangulation $K^*$ of the triangulation $K$ of $\partial P$.
	The vertices of $\Gamma_2'$ are $B_F\cup B_T$, and the edges are radial in each tetrahedron of $K$.
	Note that $\partial\Gamma_2 = \Gamma_2'$ and $X_2$ is a neighborhood of $\Gamma_2$.
	The neighborhood of $\Gamma_2'$ (restricted to a single tetrahedron) is shown in Figure~\ref{fig:core_face}.
	
	Define
	$$X_3 = \Ee_{\left[0,\frac14\right]}\cup(\Ee\cup\Tt)_{\left[\frac14,\frac34\right]}
	\hspace{1cm} \text{and}\hspace{1cm}
	\Gamma_3 = (B_E)_{\left[0,\frac12\right]}\cup(\Gamma_3')_{\left\{\frac12\right\}},$$
	where $\Gamma_3'$ has vertices $B_E\cup B_T$ and edges the radial segments in the tetrahedra of $K$ connecting these vertices. (Figure~\ref{fig:core_vert} shows a thickening of the restriction of this graph to a single tetrahedron.)
	Note that $\partial\Gamma_3 = E$, and the neighborhood of $\partial\Gamma_3$ (restricted to a single tetrahedron) is shown in yellow throughout Figure~\ref{fig:tetrahedron}.
	It is clear that $X_3$ is a neighborhood of $\Gamma_3$.
	
	The pieces $X_i$ are clearly $S_5$--invariant and comprise all of $P$.
	This, together with the above discussion establishes part (1) of the theorem.
	We now turn our attention to analyzing the intersections $H_i = X_i\cap X_{i-1}$.
	See Figure~\ref{fig:depth_schem_H}.

\begin{figure}[ht!]
	\centering
	\includegraphics[width = .9\linewidth]{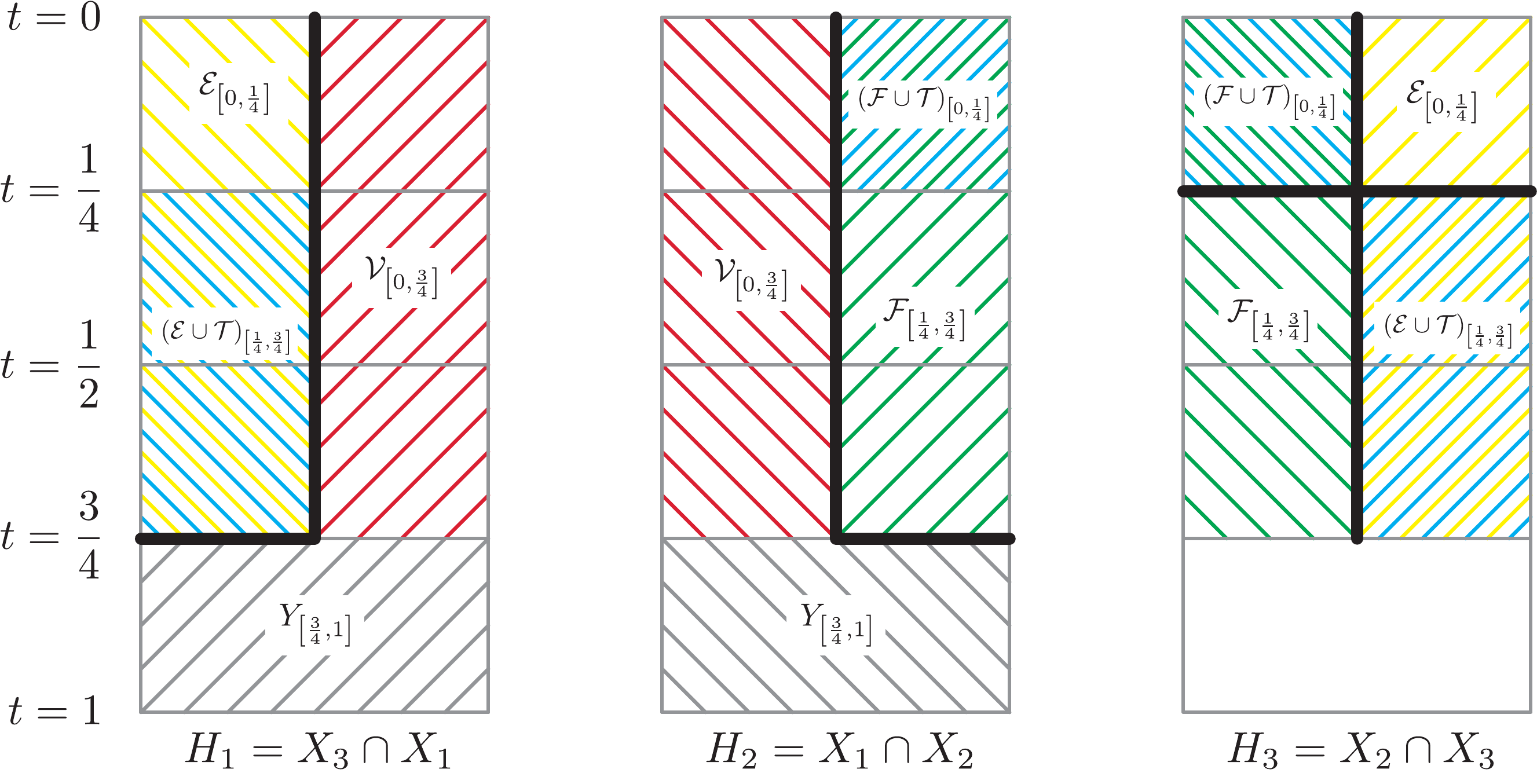}
	\caption{A schematic illustrating how the $H_i$ (represented as black segments) arise as the intersections of the sectors $X_i$ and $X_{i+1}$, which are decomposed based on the fundamental pieces $\Vv$, $\Ee$, $\Ff$, and $\Tt$ at each depth level.}
	\label{fig:depth_schem_H}
\end{figure}

	We have
	$$H_1 = (\Vv\cap\Ee)_{\left[0,\frac14\right]}\cup(\Vv\cap(\Ee\cup\Tt))_{\left[\frac14,\frac34\right]}\cup (\Ee\cup\Tt)_{\left\{\frac34\right\}}.$$

	To identify $H_1$ as the neighborhood of a graph $\gamma_1$, notice that $\Vv\cap\Ee$ is a collection of disks; let $B_{\Vv\cap\Ee}$ denote the barycenters of these disks.
	Also notice that $\Vv\cap(\Ee\cup\Tt)$ is the disjoint union of planar surfaces; let $\gamma_1'$ be a graph that contains $B_{\Vv\cap\Ee}\cup B_{\Vv\cap\Tt}$ as vertices and onto which this planar surface deformation retracts.
	Let $\gamma_1''$ be a graph that contains $B_\Ee\cup B_\Tt\cup B_{\Vv\cap\Tt}$ as vertices and onto which $\Ee\cup\Tt$ deformation retracts.
	Define 
	$$\gamma_1 = (B_{\Vv\cap\Ee})_{\left[0,\frac12\right]} \cup
	(\gamma_1')_{\left\{\frac12\right\}} \cup
	(B_{\Vv\cap\Tt})_{\left[\frac12,\frac34\right]} \cup
	(\gamma_1'')_{\left\{\frac34\right\}}.$$
	By construction, $H_1$ is a neighborhood of $\gamma_1$.
	Note that $\partial\gamma_1 = B_{\Vv\cap\Ee}$.
	
	We have
	$$H_2 = (\Vv\cap(\Ff\cup\Tt))_{\left[0,\frac14\right]}\cup(\Vv\cap\Ff)_{\left[\frac14,\frac34\right]}\cup\Ff_{\left\{\frac34\right\}}.$$
	
	To identify $H_2$ as the neighborhood of a graph $\gamma_2$, notice that $\Vv\cap(\Ff\cup\Tt)$ is a disjoint union of planar surfaces and $\Vv\cap\Ff$ is a disjoint union of disks; let $\gamma_2'$ denote a graph onto which this planar surface retracts such that the barycenters $B_{\Vv\cap\Ff}$ are included as vertices.
	Let $\gamma_2''$ denote the segments of $K'$ connecting $B_\Ff$ and $B_{\Vv\cap\Ff}$.
	Define
	$$\gamma_2 = \gamma_2' \cup
	(B_{\Vv\cap\Ff})_{\left[0,\frac34\right]} \cup
	(\gamma_2'')_{\left\{\frac34\right\}}.
	$$
	By construction, $H_2$ is a neighborhood of $\gamma_2$.
	Note that $\partial\gamma_2 = \gamma_2'$.
	
	We have
	$$H_3 = (\Ee\cap(\Ff\cup\Tt))_{\left[0,\frac14\right]}
	\cup \Tt_{\left\{\frac14\right\}}
	\cup(\Ff\cap(\Ee\cup\Tt))_{\left[\frac14,\frac34\right]}.$$

	To identify $H_3$ as the neighborhood of a graph $\gamma_3$, notice that $\Ee\cap(\Ff\cup\Tt)$ is a disjoint union of annuli.
	Let $\gamma_3'$ be a graph that includes $B_{\Ee\cap\Tt}$ as vertices and onto which this the collection of annuli deformation retracts.
	Similarly, note that $\Ff\cap(\Ee\cup\Tt)$ is a disjoint union of planar graphs.
	Choose a graph $\gamma_3''$ that includes $B_{\Ff\cap\Tt}$ as vertices and onto which the planar surface deformation retracts.
	Finally, let $\gamma_3'''$ be a graph that includes $B_{\Ee\cap\Tt}$ and $B_{\Ff\cap\Tt}$ as vertices and onto which $\Tt$ deformation retracts.
	Define
	$$\gamma_3 = \gamma_3' \cup 
	(B_{\Ee\cap\Tt})_{\left[0,\frac14\right]} \cup
	(\gamma_3''')_{\left\{\frac14\right\}} \cup
	(B_{\Ff\cap\Tt})_{\left[\frac14,\frac34\right]} \cup
	(\gamma_3'')_{\left\{\frac34\right\}}.$$
	By construction, $H_3$ is a neighborhood of $\gamma_3$.
	Note that $\partial\gamma_3 = \gamma_3'$.
	
	The $H_i$ are clearly $S_5$--invariant by construction, competing the proof of part (2).
	Finally, we consider $\Sigma = X_1\cap X_2\cap X_3$.
	We have
	$$\Sigma = (\Vv\cap\Ee\cap(\Ff\cup\Tt))_{\left[0,\frac14\right]} \cup
	(\Vv\cap(\Ff\cup\Tt)\cap(\Ee\cup\Tt))_{\left\{\frac14\right\}}$$
	$$ \cup
	(\Vv\cap\Ff\cap(\Ee\cup\Tt))_{\left[\frac14,\frac34\right]}\cup 
	(\Ff\cap(\Ee\cup\Tt))_{\left\{\frac34\right\}}.$$
	First, note that $\Vv\cap\Ee\cap(\Ff\cup\Tt) = \partial \Vv\cap\Ee$ is a disjoint union of circles, since $\Vv\cap\Ee$ is a disjoint union of disks.
	Similarly, $\Vv\cap\Ff\cap(\Ee\cup\Tt) = \partial\Vv\cap\Ff$ is a disjoint union of circles.
	The intersection of these two collections of circles is $\Vv\cap\Ee\cap\Ff$, which is a disjoint union of arcs.
	It follows that
	$$(\Vv\cap\Ee\cap(\Ff\cup\Tt))_{\left[0,\frac14\right]} \cup (\Vv\cap\Ff\cap(\Ee\cup\Tt))_{\left[\frac14,\frac34\right]}$$
	is the union of two collections of vertical annuli that meet along segments of their boundaries at $(\Vv\cap\Ee\cap\Ff)_{\left\{\frac14\right\}}$.
	Next,
	$$(\Vv\cap(\Ff\cup\Tt)\cap(\Ee\cup\Tt))_{\left\{\frac14\right\}} = ((\Vv\cap\Ee\cap\Ff)\cup(\Vv\cap\Tt))_{\left\{\frac14\right\}}$$
	is a collection of disks whose boundary decomposes into collections of arcs as
	$$(\Vv\cap\Tt\cap\Ee)_{\left\{\frac14\right\}}\cup(\Vv\cap\Tt\cap\Ff)_{\left\{\frac14\right\}}.$$
	These disks cap off the remaining boundary of the union of the vertical annuli in the level set $Y_\frac14$.
	We conclude by observing that $(\Ff\cap(\Ee\cup\Tt))_{\left\{\frac34\right\}}$ is a disjoint union of planar surfaces with boundary $(\Vv\cap\Ff\cap(\Ee\cup\Tt))_{\left\{\frac34\right\}}$, capping off the boundary of the union of vertical annuli in the level set $Y_\frac34$.
	
	By construction, $\Sigma$ is $S_5$--invariant, and
	$$\partial\Sigma = \Vv\cap\Ee\cap(\Ff\cup\Tt),$$
	which shows that $\Sigma$ is neatly embedded.
	This completes the proof of the proposition.
\end{proof}

We now show that given an equivariant triangulation of a $4$--manifold $X$, we can glue up these decompositions across all pentachora of the triangulation to give an equivariant trisection of $X$.

\begin{proposition}
\label{prop:exist_tri}
	Let $X$ be a closed, $4$--dimensional $G$--manifold.
	Let $K$ be a $G$--simplicial complex triangulating $X$.
	The triangulation $K$ gives rise to $G$--equivariant trisection of $X$ restricting on each pentachoron to the decomposition of Proposition~\ref{prop:pentachoron}.
\end{proposition}

\begin{proof}
	For each pentachoron $P$ of $K$, since the induced action on $P$ is determined by the induced action on the vertices of $P$, the induced action on $P$ is as a subgroup of $S_5$.
	Thus, we can invoke the $G$--invariant decomposition $P = X_1^P\cup X_2^P\cup X_3^P$ provided by Proposition~\ref{prop:pentachoron}.
	For $i\in\Z_3$, define
	$$X_i = \bigcup_{P\in K}X_i^P, \hspace{2cm}
	H_i = \bigcup_{P\in K}H_i^P, \hspace{1cm} \text{and}\hspace{1cm}
	\Sigma = \bigcup_{P\in K}\Sigma^P.$$
	Note that for every pentachoron $P$, the decomposition is identical on every tetrahedral face of $P$.
	It follows that the gluings described above are perfectly aligned at each face of each pentachoron.
	For example, it follows immediately that $\Sigma$ is a closed surface.
	 
	We claim that $X = X_1\cup X_2\cup X_3$ is a trisection of $X$.
	First, it follows that $H_i = X_{i-1}\cap X_i$ and $\Sigma = X_1\cap X_2\cap X_3$, since this was true in each pentachoron.
	Furthermore, it is immediate that the $X_i$, the $H_i$, and $\Sigma$ are $G$--invariant for the same reason.
	
	To see that each $X_i$ is a 4--dimensional 1--handlebody, recall that, in each pentachoron, $X_i^P$ is a 4--dimensional neighborhood of the graph $\Gamma_i^P$ and that $X_i^P\cap\partial P$ is a 3--dimensional neighborhood of a graph $\partial\Gamma_i^P$ in $\partial P$.
	If we let $\Gamma_i = \bigcup_{P\in T}\Gamma_i^P$, then we see that $X_i$ is a 4--dimensional neighborhood of $\Gamma_i$, hence a 4--dimensional 1--handlebody.
	An identical argument featuring the graphs $\gamma_i^P$ shows that each $H_i$ is a 3--dimensional handlebody.
	
	To see that the action of $G$ on $X_i$ is linearly parted, it suffices to find a collection $\Bb_i$ of neatly embedded 3--balls in $X_i$ such that $X_i\setminus\Bb_i$ is a collection $\Zz_i$ of 4--balls such that $G$ acts linearly on each 4--ball in $\Zz_i$.
	We can identify the $\Bb_i$,  the resulting 4--balls, and the action of $G$ as follows; cf. Figure~\ref{fig:depth_schem_X}, where the location of the $\Bb_i$ is indicated with purple arcs.
	
	First, we have
	$$\Bb_1 = \bigcup_{P\in K}\Vv^P_{\left\{\frac34\right\}}
	\hspace{5mm} \text{and} \hspace{5mm}
	\Zz_1 = \left(\bigcup_{P\in K}\Vv^P_{\left[0,\frac34\right]}\right)\bigsqcup\left(\bigcup_{P\in K}Y^P_{\left[\frac34,1\right]}\right).$$
	$$\Bb_2 = \bigcup_{P\in K}(\Ff\cap\Tt)^P_{\left[0,\frac14\right]}
	\hspace{5mm} \text{and} \hspace{5mm}
	\Zz_2 = \left(\bigcup_{P\in K}\Ff^P_{\left[0,\frac34\right]}\right)\bigsqcup\left(\bigcup_{P\in K}\Tt^P_{\left[0,\frac14\right]}\right).$$
	
	$$\Bb_3 = \bigcup_{P\in K}(\Ee\cap\Tt)^P_{\left[\frac12,\frac34\right]}
	\hspace{5mm} \text{and} \hspace{5mm}
	\Zz_3 = \left(\bigcup_{P\in K}\Ee^P_{\left[0,\frac34\right]}\right)\bigsqcup\left(\bigcup_{P\in K}\Tt^P_{\left[\frac12,\frac34\right]}\right).$$
	
	Among the $\Zz_i$, there are two different types of $4$--balls: those that are contained in a single pentachoron (the latter terms of $\Zz_1$ and $\Zz_3$) and those that are thickenings of cone neighborhoods of  barycenters.
	For the former, it suffices to note that a simplicial $G$--action on a pentachoron is clearly linear, and these subsets are star-like $4$--balls inside that linear action, hence the induced actions are linear.

	For the second type of $4$--balls, we use the example of the first term $\bigcup_{P\in K}\Vv^P_{\left[0,\frac34\right]}$ of $\Zz_1$:
	Let $v$ be a vertex of the triangulation $K$ of $X$. There may be many pentachora meeting at $v$, call this collection $S_v$.
	For each $P\in S_v$, $\Vv^P_{\left[0,\frac34\right]}$ consists of the star neighborhood of $v$ in the second barycentric subdivision of the $3$--skeleton of $P$ thickened in the depth direction.
	These thickened star neighborhoods are glued together around $v$ to give $N_v = \bigcup_{P\in S_v}\Vv^P_{\left[0,\frac34\right]}$. The $G$--action on $N_v$ is simplicial, swapping thickened star neighborhoods about.
	In particular, the action on $N_v$ is a cone, since $G$ respects the height function on each $P$ and the star neighborhoods of $v$ which were thickened were cones.
	Thus by~\cite[Corollary~2.4]{MeiSco_LP}, the $G$--action on $N_v$ is linear.
	We note that in fact $N_v$ is equivariantly isotopic to an honest star neighborhood of $v$ in the second barycentric subdivision of $K$ by an isotopy preserving the 3--skeleton of $K$.
	
	The other $4$--balls of this sort occur in $\Zz_2$, which has cone neighborhoods of the barycenters of faces and tetrahedra in $K$, and $\Zz_3$, which has cone neighborhoods of the barycenters of edges of $K$.
	Since all these $4$--balls are linear, the $G$--action on each $4$--dimensional $1$--handlebody $X_i$ is linearly parted and thus $X_1\cup X_2\cup X_3$ is an equivariant trisection of $X$.
\end{proof}

\subsection{Triangulated invariant surfaces}
\label{subsec:exist_bridge}
	
We now show that if $\Ss$ is a $G$--invariant surface triangulated by $K$, the above trisection is a bridge trisection of $\Ss$.

\begin{proposition}
\label{prop:exist_bridge}
	Let $X$, $G$, and $K$ be as above, and let $\Ss$ be an embedded, $G$--invariant surface triangulated by $K$.
	The surface $\Ss$ is in equivariant bridge trisected position with respect to the trisection given in Proposition~\ref{prop:exist_tri}.
\end{proposition} 

\begin{proof}
	The surface $\Ss$ is the union of a collection of triangular faces coming from the tetrahedra comprising the boundaries of the pentachora of $K$.
	These are precisely the triangular faces illustrated in Figure~\ref{fig:tetrahedron}.
	To establish that $\Ss$ is in bridge position with respect to the trisection described above, we need to establish (1) that the pairs $(X_i, \Dd_i)$, where $\Dd_i= X_i\cap \Ss$, are all linearly parted as pairs, and (2) each $\Tt_i = H_i\cap\Ss$ is a trivial tangle.	

	For (1), note that the disk components of $\Dd_1$, $\Dd_2$, and $\Dd_3$ are unions of the red, green, and yellow triangles respectively in the faces of the tetrahedra in Figure~\ref{fig:tetrahedron}.
	It is thus clear by definition that each $0$--handle of $X_i$ contains at most one disk of $\Dd_i$, so it suffices to show each $D\in \Dd_i$ is equivariantly boundary-parallel in its $0$--handle.
	We claim that, for the radial Morse function on the $0$--handle containing $D$, there is only one critical point, namely where $D$ intersects $\Gamma_i$.
	The graphs $\Gamma_1$, $\Gamma_2$, and $\Gamma_3$ intersect the faces of the tetrahedra shown in Figure~\ref{fig:tetrahedron} in the barycenters $B_\Vv$, $B_\Ff$, and $B_\Ee$, respectively.	
	These barycenters are the critical points of the induced Morse functions on $X_1$, $X_2$, and $X_3$, respectively.
	These barycenters are all vertices of the first barycentric subdivision $K'$ of $K$ and lie in the 2--skeleton $K^2$ of $K$.
	By definition, $D$ is a star neighborhood in $K^2\cap \Ss$ of one of these barycenters, so it contains exactly one critical point with respect to the radial Morse  function on its $0$--handle $Z$. 
	By Lemma~\ref{lem:boundary_parallel_in_ball}, $D$ has an almost-equivariant bridge-ball $\delta$ in $Z$.
	By the lemma, $\delta$ can be assumed to miss any prescribed invariant collection of 3--balls, so it is a bridge-ball for $D$ in $X_i$ missing the linearly parting ball-system $\Bb_i$ identified in Proposition~\ref{prop:exist_tri}.
	Repeating for each disk in $\Dd_i$ gives the required $G$--invariant, almost-disjoint bridge-balls $\Delta_i$ for the entire trivial disk-tangle.
	
	An identical approach involving the graphs $\gamma_i$ for which the handlebodies $H_i$ are neighborhoods works to verify that the tangles $\Tt_i$ are trivial tangles, hence equivariantly trivial by Lemma~\ref{lem:equiv_parallel}.
	We have already seen that the $\Dd_i$ are star neighborhood disks in the 2--skeleton of $\Ss$, $K^2\cap \Ss$.
	It follows that the $\Tt_i$ are arcs in the 2--skeleton.
	We will check that each such arc is a star neighborhood of a vertex of $\gamma_i\cap K^2$ intersecting $\gamma_i$ in a single point.
	The intersections of $\gamma_1$, $\gamma_2$, and $\gamma_3$ with $K^2$ are, respectively, $B_{\Vv\cap\Ee}$, $B_{\Vv\cap\Ff}$, and $B_{\Ee\cap\Ff}$.
	These barycenters are the midpoints of each arc in $\Tt_i$.
	Since each strand of $\Tt_i$ has a single critical point, $\Tt_i$ is trivial, as desired.
\end{proof}

We can now prove the main theorem of this section.

\begin{proof}[Proof of Theorem~\ref{thm:exist_main}]
	For each $i\in\{1,\ldots,N\}$, we can apply Theorem~\ref{thm:Illman} to the pair $(X,\Ss_i)$ to get an equivariant triangulation $f_i\colon K_i\to X$ for which $\Ss_i$ is triangulated by $f_i$.
	Then, we can apply a further result of Illman that shows that these equivariant triangulations have $G$--equivariant subdivisions $f\colon K_i'\to X$ that are $G$--isomorphic~\cite{Ill_78_Smooth-equivariant-triangulations}.
	Let $K=K'_1$, so $f\colon K\to X$ is a $G$--equivariant triangulation of $X$ such that each of the $\Ss_i$ is triangulated by $f$.
	By Proposition~\ref{prop:exist_tri}, this triangulation gives rise to a $G$--equivariant trisection for $X$.
	By Proposition~\ref{prop:exist_bridge}, each surface $\Ss_i$ is in equivariant bridge position with respect to this trisection, as desired.
\end{proof}

\begin{remark}
\label{rmk:sector1}
	If the action of $G$ on $X$ has isolated fixed points, then these fixed points will occur as vertices in any $G$--equivariant triangulation of $X$.
	As a consequence of the construction given above of the equivariant trisection, such fixed points will occur in $X_1$.
	On the other hand, there exist equivariant trisections in which isolated fixed points occur in more than one sector; see, for example, the trisections of Subection~\ref{subsec:PU(3)_trisections}.
	Similar analysis of the other sectors $X_2$ and $X_3$ shows that the possible stabilizers of $0$--handles in these sectors are limited by the combinatorics of the skeleta these sectors are built from.
	In this sense, the equivariant trisections constructed in the proof of Theorem~\ref{thm:exist_main} are quite special. 
\end{remark}

\subsection{Equivariant handle-decompositions}
\label{subsec:handles}

We conclude this section with a brief discussion of the connection between equivariant trisections and equivariant handle-decompositions.
See Wasserman~\cite[Section 4]{Was_69_Equivariant-differential-topology} for a complete, formal treatment, noting that the ``handle bundles over orbits'' of that paper, for finite $G$, coincide with orbits of handles as described below.

Let $M$ be a compact $G$--manifold.
A handle-decomposition of $M$ is \emph{$G$--equivariant} if

\begin{enumerate}
	\item the union of the handles of a given index is $G$--invariant, and
	\item the induced action on each handle $\frak h$ preserves the decomposition as a product of disks and acts linearly on each disk factor.
\end{enumerate}

\begin{corollary}
\label{cor:handle}
	Let $X$ be a closed, 4--dimensional $G$--manifold.
	Then, $X$ admits a $G$--equivariant handle-decomposition.
\end{corollary}

\begin{proof}
	By Theorem~\ref{thm:exist_main}, $X$ admits a $G$--equivariant trisection.
	The process of turning a trisection into a handle-decomposition in the non-equivariant setting is described in~\cite{MeiSchZup_16_Classification-of-trisections}.
	In the present setting, the linearly parting ball-system for $X_1$ gives the equivariant 0--handles and 1--handles.
	To get the 2--handles, let $\gamma$ be an invariant cut-system of curves for $H_3$; let $D\subset H_3$ be the invariant system of cut-disks bounded by $\gamma$.
	Let $\frak h = \overline{\nu(D)}$ be an equivariant thickening that intersects $\Sigma$ in an annular neighborhood of $\gamma$ and intersects $H_1\cup H_2$ in a solid torus neighborhood of $\gamma$.
	Then $\frak h$ is a collection of 4--dimensional 2--handles attached along $\gamma\subset Y_1 = \partial X_1$ with framing given by $\frak h\cap\Sigma$.
	Since $\frak h$ is a closed, equivariant tubular neighborhood of $D$, the induced action on each handle is linear and preserves the product structure.
	
	Removing the interior of $\frak h$ from $X_2\cup X_3$ gives $X_2'\cup_{H_3'}X_3'$, where $X_2'$ and $X_3'$ are just ``indented'' versions of $X_2$ and $X_3$, and $H_3'$ is a collection of 3--balls.
	Both $X_2'$ and $X_3'$ are still linearly parted 4--dimensional 1--handlebodies; let $\Bb_2'$ and $\Bb_3'$ be the corresponding linearly parting ball-systems.
	Since $H_3'$ is an invariant collection of 3--balls, $Z = X_2'\cup_{H_3'}X_3'$ is a $G$--invariant 4--dimensional 1--handlebody.
		
	To see that $Z$ is linearly parted as a $G$--equivariant handlebody, we apply the techniques of~\cite[Section~4]{MeiSco_LP} to obtain new linearly parting ball-systems $\Bb_2''$ and $\Bb_3''$ for $X_2'$ and $X_3'$ that are disjoint from $H_3'$.
	Then, $\Bb = H_3'\cup\Bb_2''\cup\Bb_3''$ is a linearly parting ball-system for the $G$--action on $Z$, as desired.
\end{proof}

For details regarding how to obtain a Kirby diagram from a trisection diagram, see~\cite{MeiSchZup_16_Classification-of-trisections,Kep_22_An-algorithm-taking-Kirby-diagrams}.
The algorithm given by Kepplinger may be difficult to apply in the equivariant setting.
There are necessary relationships between the framings of the 2--handles and the types of fixed points arising in the group action (see \cite[Section 1]{Edm_87_Construction-of-group-actions}) that may be interesting to investigate in future work.

Kepplinger also gives an algorithm for converting a Kirby diagram into a trisection diagram.
It should be possible to adapt this algorithm to one that is capable of converting an equivariant Kirby diagram to an equivariant trisection diagram.

\section{Quotients of equivariant trisections}
\label{sec:quotients}

In this section, we present necessary and sufficient conditions for the quotient of a $G$--equivariant trisection by a normal subgroup $N\triangleleft G$ to be a $G/N$--equivariant trisection.
The main theorem is the following.

\begin{theorem}
\label{thm:quotient_trisection}
	Let $\TT$ be a $G$--equivariant trisection of $X$.
	Let $N\triangleleft G$ be normal.
	Then, the following are equivalent.
	\begin{enumerate}
		\item $X/N$ admits a smooth structure making it a $G/N$--manifold.
		\item $\TT/N$ is a $G/N$--equivariant trisection of $X/N$.
		\item $\Stab_N(x)$ acts as a rotation group on $\nu(x)$ for each $x\in X$.
	\end{enumerate}
\end{theorem}

\begin{remark}
\label{rmk:conjugation}
	One consequence of this theorem is that the choice of the trisection $\TT$ for $X$ determines a smooth structure on $X/N$: specifically, the one determined by $\TT/N$.
	In general, it is not known whether $X/N$ admits a unique smooth structure compatible with the quotient map $X\to X/N$.
	Note that despite this, the above shows that the smooth structure on $X/N$ is an invariant of the $G$--equivariant trisection $\TT$; what is unknown is whether the smooth structure on $X/N$ is an invariant of the \emph{$G$--action} on $X$.
	
	Equivariant trisections may provide an avenue to study this uniqeness problem: In particular, if $X/N$ admits distinct smooth structures compatible with the quotient map, then it is possible the information is encoded in a pair of non-diffeomorphic equivariant trisections for the same action on $X$.
	For a positive result on this problem, Hambleton and Hausmann have shown that, for $\Z_2$ actions with codimension $2$ fixed-point set, there \emph{is} a unique smooth structure on the quotient compatible with the quotient map~\cite[Subsection 7.4]{HamHau_11_Conjugation-spaces-and-4-manifolds}.
\end{remark}

We have the following diagrammatic corollary of the above theorem.

\begin{corollary}
\label{cor:quotient_diagram}
	Let $\DD$ be a $G$--invariant trisection diagram for $X$.
	Let $N\triangleleft G$ be normal.
	Then $X/N$ is a $G/N$--manifold if and only if $\DD/N$ is a $G/N$--invariant trisection diagram.
\end{corollary}

We now introduce the relevant background on group actions before proving Theorem~\ref{thm:quotient_trisection} and Corollary~\ref{cor:quotient_diagram} in Subsection~\ref{subsec:quotient_proofs}.

\subsection{Quotient actions and manifolds}
\label{subsec:quotient_actions}

There is a body of literature, dating back to work of Miha\^{i}lova~\cite{Mih_78_Finite-imprimitive-groups} on groups generated by rotations, concerning which linear group actions on $B^n$ have quotient homeomorphic to $B^n$.
The problem was resolved completely in the topological category by Lange~\cite{Lan_19_When-is-the-underlying-space}.
The solution to the problem is complicated when $n\geq 5$ by the Double Suspension Theorem~\cite{Can_79_Shrinking-cell-like-decompositions},
but when $n\leq 4$ and $G$ acts orientation-preservingly, the situation is simpler.
Central is the concept of a group acting by rotations.

\begin{definition}
	A linear diffeomorphism $r$ of $B^n$ is a \emph{rotation} if its fixed-point set has codimension two.
	A finite group $G$ acting linearly on $B^n$ \emph{acts as a rotation group} if it is generated by (finite order) rotations.
\end{definition}

Note that a rotation $r$ fixes a hyperplane $\R^{n-2}$ in $\R^n$, so the codimension-two fixed-point set referenced above is a linear $B^{n-2}$ in $B^n$.

The key result that we make use of in this section is a theorem of Lange characterizing when quotients of $\R^n$ by linear group actions are topological manifolds~\cite{Lan_19_When-is-the-underlying-space}.
When $n\leq 4$ and the action is orientation-preserving, Lange's theorem holds in the smooth category as the following.

\stepcounter{theorem}
\begin{named}{Ball Quotient Theorem~\thetheorem}
\label{thm:ball_quotient}
	Let $G$ act linearly and orientation-preservingly on $B^n$ with $n\leq 4$.
	\begin{enumerate}
		\item $B^n/G$ is diffeomorphic to $B^n$ if and only if $G$ acts as a rotation group.
		\item If $N\vartriangleleft G$ acts as a rotation group, then the action of $G/N$ on the smooth ball $B^n/N$ is linear.
	\end{enumerate}	
\end{named}

\begin{proof}
	We first prove both directions of part (1).
	Suppose $B^n/G$ is diffeomorphic to $B^n$.
	Then $\R^n/G$ is diffeomorphic (hence, homeomorphic) to $\R^n$, so by Theorem~A of
	~\cite{Lan_19_When-is-the-underlying-space}, $G$ acts as a rotation group.
	(There are no reflections, since $G$ preserves orientation.)
	
	Now, suppose $G$ acts as a rotation group.
	Then $\R^n/G$ is homeomorphic to $\R^n$, by Theorem~A of ~\cite{Lan_19_When-is-the-underlying-space}.
	Since the action is radially invariant, the action of $G$ on $B^n$ is the cone of the action of $G$ on $S^{n-1}$, and $B^n/G$ is the cone on $S^{n-1}/G$, which is homeomorphic to $S^{n-1}$.
	Since $n\leq 4$, the solution to the 3--dimensional Poincar\'e Conjecture~\cite{Per_02_The-entropy-formula,Per_03_Finite-extinction,Per_03_Ricci-flow} implies that $S^{n-1}/G$ is diffeomorphic to $S^{n-1}$, so the cone $B^n/G$ is diffeomorphic to $B^n$.
	
	We now prove part (2).
	The conclusion that $B^n/N$ is diffeomorphic to $B^n$ follows from part (1), since $N$ acts as a rotation group.
	It remains to show that the action of $G/N$ on $B^n/N$ is linear; for this, we refer the reader to Section~6.3 of~\cite{Lan_16_Characterization-of-finite-groups}.
\end{proof}

With the Ball Quotient Theorem in hand, we proceed to the proofs of the main theorems.

\subsection{Proofs of Theorem~\ref{thm:quotient_trisection} and Corollary~\ref{cor:quotient_diagram}}
\label{subsec:quotient_proofs}

\begin{proof}[Proof of Theorem \ref{thm:quotient_trisection}]	
	First, note that part (2) implies part (1) trivially.
	We will show that part (1) implies part (3), implies part (2), to complete the proof.
	
	Suppose part (1) is true.
	Let $x\in X$, and let $B$ be an equivariant neighborhood of $x$ in $X$.
	Then $\Stab_N(x)$ acts smoothly on the 4--ball $B$, and the quotient $B/\Stab_N(x)$ is a 4--ball, since $X/N$ is a manifold.
	By the~\ref{thm:ball_quotient}(1), $\Stab_N(x)$ acts as a rotation group, as desired.
	
	Now suppose part (3) is true.
	Let $\Bb_i$ be a linearly parting ball-system for $X_i$, and let $\Zz_i = X_i\setminus\Bb_i$ denote the corresponding collection of 4--balls.
	Each $B\in \Bb_i$ is $N$--equivariant, and $\Stab_N(B) = \Stab_N(x_B)$, where $x_B$ is the center of $B$, since the action of $N$ is linear.
	Therefore, by the~\ref{thm:ball_quotient}(2), $B/\Stab_N(B)$ is a (neatly embedded) 3--ball, and the natural $G/N$--action is linear.
	Similarly, for each 4--ball $Z\subset \Zz_i$, the quotient $Z/\Stab_N(Z)$ is a 4--ball, and the natural $G/N$--action is linear.
	It follows that $X_i/N$ is a linearly parted 4--dimensional 1--handlebody with parting ball system given by the $G/N$--invariant ball-system $\Bb_i/N$, which cuts $X_i/N$ into the 4--balls $\Zz_i/N$.
	
	An identical argument featuring an invariant disk-system shows that $H_i/N$ is a $G/N$--invariant handlebody, which also implies that $\Sigma/N$ is $G/N$--invariant closed surface.
	It follows that $\TT/N$ is a $G/N$--equivariant trisection of $X/N$, as desired.
\end{proof} 

\begin{proof}[Proof of Corollary~\ref{cor:quotient_diagram}]
	Let $\DD$ be a $G$--invariant trisection diagram for $X$.
	The diagram $\DD$ uniquely determines a $G$--equivariant trisection $\TT$ of $X$, by Proposition~\ref{prop:diag_spine}.
	
	First, suppose that $X/N$ is a $G/N$--manifold.
	Then, by Theorem~\ref{thm:quotient_trisection}, $\TT/N$ is a $G/N$--equivariant trisection.
	The curves of $\DD$ bound invariant disk-systems in the $3$--dimensional handlebodies of the spine of $\TT$.
	As in the proof of Theorem~\ref{thm:quotient_trisection} above, the quotient of these disk-systems are disk-systems for the $3$--dimensional handlebodies of $\TT/N$.
	It follows that $\DD/N$ is a $G$--invariant trisection diagram for $\TT/N$ on $\Sigma/N$, as desired.

	Conversely, consider the quotient $\DD/N$ of $\DD$, which we assume to be a $G/N$--invariant trisection diagram.
	The quotient $H_i/N$ of each $3$--dimensional handlebody $H_i$ in the spine of $\TT$ is a handlebody, since the action of $N$ on $H_i$ is linearly parted by~\cite[Corollary 3.6]{MeiSco_LP}, and all linear actions on $3$--balls are actions by rotations by Euler's Theorem.
	The quotient $\DD/N$ consists of cut-curves for the $H_i/N$.
	This establishes that the spine determined by $\DD/N$ is the quotient of the spine determined by $\DD$.
	
	It remains to check that $X_i/N$ is a smooth manifold for each sector $X_i$ of $\TT$.
	As in the proof of Theorem~\ref{thm:quotient_trisection} above, we can see that for each sector $X_i$ of $\TT$, the quotient $X_i/N$ consists of the quotients of $0$--handles of $X_i$ attached along quotients $\Bb_i/N$ of parting $3$--balls $\Bb_i$ for $X_i$.
	The $\Bb_i/N$ are neatly embedded $3$--balls in $X_i/N$, by Euler's Theorem and the~\ref{thm:ball_quotient}(1).
	
	Let $Z$ be a $0$--handle of $X_i$.
	Let $N_Z$ be the stabilizer of $Z$, and let $R$ be the subgroup of $N_Z$ generated by all elements acting as rotations on $Z$. 
	
	\textbf{Claim 1:} $R$ is normal in $N_Z$.
	
	Consider $n^{-1}rn$ for a rotation $r$ and $n\in N_Z$. 
	Let $x$ be a point on the disk $F$ fixed by $r$. 
	Observing that $n^{-1}\cdot x$ is fixed by $n^{-1}rn$, we see that $s = n^{-1}rn$ is a rotation since it fixes a disk $n^{-1}\cdot F$, and thus $rn = ns$.
	For an element $w\in R$, we have that $w$ is a word in $N_Z$ consisting of elements acting as rotations, and thus for $n^{-1}wn$ by the above we can commute the rightmost $n$ leftward by replacing letters in $w$ with possibly novel rotations.
	This gives $n^{-1}wn = w'$, where $w'$ is a word in $N_Z$ consisting of rotations and thus $w'\in R$, proving the claim.
	
	\textbf{Claim 2:} $N_Z/R$ is trivial.
	
	Note that if $N_Z/R$ is trivial, then $N_Z$ was acting as rotations, and $Z/N_Z$ is a $4$--ball.
	Consider now the $N_Z/R$--action on $Z/R$, which is a linear action on a $4$--ball by the~\ref{thm:ball_quotient}(2).
	If an element of $N_Z$ fixes a point on $\partial Z/N$, then it acts as a rotation on this 3--sphere, hence as a rotation on the cone $Z/N$.
	It follows that each element of $N_Z/R$ acts freely on $\partial Z/R$.
	Thus $Z/N_Z = (Z/R)/(N_Z/R)$ is the cone on the spherical space form $(\partial Z/R)/(N_Z/R)$ with fundamental group $N_Z/R$.

	Returning to $X_i/N$, since $X_i/N$ is formed by connecting the quotients $Z/N$ of the 0--handles of $X_i$ with 4--dimensional 1--handles, we have that
	$$\pi_1(\partial(X_i/N)) \cong \left(\bigast_Z\pi_1(\partial(Z/N))\right)\ast F_m,$$
	where the product is over all 0--handles $Z$ in $X_i$ and $F_m$ is a free group.
	Since $\DD/N$ is a trisection diagram, $\partial(X_i/N)$ is $\#^{k_i}(S^1\times S^2)$, thus $\pi_1(\partial(X_i/N))$ is free.
	This implies that
	$\mathop{\huge\ast}_Z\pi_1(\partial(Z/N))$ is torsion-free, but since it is a free product of finite groups,
	$\pi_1(\partial(Z/N))\cong 1$ for all $Z$. By the above, $\pi_1(\partial Z/N) = N_Z/R$, so $N_Z/R\cong 1$ for each 0--handle $Z$.
\end{proof}

\begin{remark}
\label{rmk:orb_tri}
	Continuing the analysis above, we see that for any $G$--manifold $X$ with a $G$--equivariant trisection $\TT$ and any $N\triangleleft G$, the space $X/N$ admits a decomposition $\TT/N$ into pieces $X_i/N$, each of which has a parting system $\Bb/N$ (coming from the linearly parting ball-system $\Bb$ for $X_i$) for which $(X_i/N)\setminus (\Bb/N)$ is a collection of cones on $3$--dimensional spherical space forms.
	The quotient $\DD/N$ then encodes a triple of Heegaard splittings for $\partial(X_i/N)$, which are connected sums of $3$--dimensional spherical space forms and copies of $S^1\times S^2$.
	This gives a lot of structure to the orbifold $X/N$, and motivates notion of \emph{orbifold trisection}, which could warrant further study.
\end{remark}

\subsection{Quotients of bridge trisections}
\label{subsec:quotient_bridge_trisection}

Having determined above some necessary and sufficient conditions for the quotient of a $G$--equivariant trisection to be an equivariant trisection, we now address the scenario in which the set-up includes a $G$--invariant collection $\Ss$ of surfaces in bridge trisected position.
An important special case is when $\Ss$ is the codimension-two portion of $\Fix(N)$, for some $N\vartriangleleft G$, which is $G$--invariant by normality: If $s\in\Ss$ is fixed by $n\in N$, then, for any $g\in G$, we have
$$n\cdot(g\cdot s) = ng\cdot s = gm\cdot s = g\cdot(m\cdot s) = g\cdot s,$$
where $m = g^{-1}ng$.
So, $g\cdot s$ is fixed by $N$.

Our first task is to develop a version of the~\ref{thm:ball_quotient} for ball-disk pairs.
Note that if a $G$--invariant surface $\Ss$ has the property that an element $g\in G$ acts on $\Ss$ with $1$--dimensional fixed-point set, then $\Ss/G$ will have boundary, so a bridge trisection is out of the question.
We say that $G$ acts on $\Ss$ \emph{without reflections} if no element of $G$ has $1$--dimensional fixed-point set in $\Ss$.
This is equivalent to the restriction that, on any equivariant sub-disk of $\Ss$, the induced action of the stabilizer is by rotations (or trivial).

Note that this is not quite the same as asking that $G$ act orientation preservingly on $\Ss$, since we allow, for example, the antipodal involution of $S^2$.

\begin{theorem}
\label{thm:ball_quotient_redux}
	Let $G$ act linearly and orientation-preservingly on $B^n$ with $n\leq 4$.
	Let $D\subseteq B^n$ be a neatly embedded, $G$--invariant $(n-2)$--disk.
\begin{enumerate}
	\item $(B^n/G, D/G)$ is diffeomorphic to $(B^n, D^{n-2})$ if and only if $G$ acts as a rotation group on $B^n$, the induced action of $G$ on $D$ is by rotations, and $D$ is equivariantly boundary-parallel.
	\item If $N\vartriangleleft G$ acts as a rotation group on $B^n$, the induced action of $N$ on $D$ is by rotations, and $D$ is equivariantly boundary-parallel with respect to the action by $G$, then the action of $G/N$ on the smooth ball $B^n/N$ is linear and $D/N$ is equivariantly boundary-parallel with respect to the action of $G/N$.
\end{enumerate}
\end{theorem}

\begin{proof}
	First, suppose that $(B^n/G,D/G)$ is diffeomorphic to $(B^n,D^{n-2})$.
	By the~\ref{thm:ball_quotient}(1), $G$ acts as a rotation group on $B^n$.
	Since $D/G$ is neatly embedded, the action of $G$ on $D$ is by rotations.
	Since $D/G$ is boundary-parallel in $B^n/G$, we can choose a bridge-disk $\delta\cong B^{n-1}$ for $D^{n-2}$.
	The pre-image of $\delta$ under the quotient map is a $G$--invariant, almost-disjoint collection of bridge-balls for $D$ in $B^n$.
	
	Next, suppose that $G$ acts as a rotation group on $B^n$ and $D$ and that $D$ is equivariantly boundary-parallel.
	Let $\Delta$ be a $G$--invariant collection of almost-disjoint bridge-balls for $D$.
	We note that the induced action of $G$ on each bridge-disk in $\Delta$ is by rotations, since each such $(n-1)$--ball has $D$ as an invariant portion of its boundary, and the induced action of $G$ on $D$ is by rotations.
	By the~\ref{thm:ball_quotient}(1), applied thrice, the triple $(B^n/G, \Delta/G, D/G)$ is diffeomorphic to $(B^n,\delta,D^{n-2})$, where $\delta\cong B^{n-1}$ is a bridge-ball for $D/G$ in $B^n/G$.
	
	Finally, assume $N\vartriangleleft G$ acts as a rotation group on $B^n$ and $D$ and that $D\subset B^n$ is a $G$--equivariantly boundary-parallel $(n-2)$--disk.
	The action of $G/N$ on $B^n/N$ is linear by the~\ref{thm:ball_quotient}(2).
	Let $\Delta$ be a $G$--invariant collection of almost-disjoint bridge-balls for $D$.
	As above, the induced action of $N$ on each bridge-disk of $\Delta$ is by rotations.
	For each bridge-disk $\delta\in\Delta$, the triple $(B^n/N,\delta/N,D/N)$ is diffeomorphic to $(B^n,B^{n-1},B^{n-2})$, by the~\ref{thm:ball_quotient}(1), applied thrice.
	It follows that $\Delta/N$ is a $G/N$--invariant collection of almost-disjoint bridge-disks for $D/N$.
\end{proof}

And now we have the main quotient theorem as a corollary.

\begin{theorem}
\label{thm:quotient_bridge_trisection}
	Let $\TT$ be a $G$--equivariant trisection of $X$.
	Let $N\vartriangleleft G$ so that $X/N$ is a smooth manifold.
	Let $\Ss$ be a $G$--invariant surface-link in bridge trisected position with respect to $\TT$.
	If $N$ acts on $\Ss$ without reflections, then $\TT/N$ is a $G/N$--equivariant trisection of $X/N$ such that $\Ss/N$ is in bridge trisected position with respect to $\TT/N$.
\end{theorem}

\begin{proof}
	In light of Theorem~\ref{thm:quotient_trisection}, it suffices to show that the surface-link $\Ss/N$ is in bridge trisected position with respect to $\TT/N$.
	Since $\Ss$ is in bridge trisected position, we can choose a linearly parting ball-system $\Bb_i$ for $X_i$ such that each 0--handle $Z$ in the parting contains at most one patch $D$ of the patch-tangle $\Dd_i\subset \Ss$.
	Since $N$ acts on $\Ss$ without reflections, the induced action on each patch $D$ is by rotations.
	By Theorem~\ref{thm:ball_quotient_redux}(2), each such $(Z/N,D/N)$ is diffeomorphic to $(B^4,B^2)$, the action of $G/N$ on $Z/N$ is linear, and $D/N$ is $G/N$--equivariantly boundary-parallel.
	
	An identical argument can be applied to the $(H_i,\Tt_i)$ to see that the quotient is a $G/N$--equivariantly boundary-parallel trivial tangle.
	This shows that $\Ss/N$ is in $G/N$--equivariant bridge trisected position, as desired.
\end{proof}

\section{Examples}
\label{sec:examples}

We now turn our focus to analyzing various examples of equivariant trisections.
We discuss (regular) branched coverings in general, before giving an example of an irreducible $Q_8$--equivariant trisection of $\CP^2\#\overline\CP^2$ arising within a family of branched coverings of an interesting link of projective planes in $S^4$.
We then discuss hyperelliptic actions, give examples where the existence of group actions on $T^4$ and $T^2\times S^2$ can be inferred from a diagram, and analyze the actions on $\CP^2$ by most subgroups of $PU(3)$.

\subsection{Branched covering actions}
\label{subsec:branched_coverings}

Some of the simplest and most important examples of equivariant trisections arise via covering and branched covering constructions.
We adopt the formal definition for branched coverings given by Zuddas~\cite[Section~1.1]{Zud_08_Branched-coverings}, but given our interest in group actions, we will restrict our attention to \emph{regular} branched coverings.
We now recall a fundamental and well-known feature of bridge trisections in the language of equivariant trisections.

\begin{proposition}
\label{prop:branched_cover}
	Let $p\colon \widetilde X\to X$ be a regular branched covering with branch locus $\Ss\subset X$, singular locus $\widetilde\Ss = p^{-1}(\Ss)$, and deck transformation group $G$.
	Suppose that $\TT$ is a bridge trisection of the pair $(X,\Ss)$.
	Then, $p^{-1}(\TT)$ is a $G$--equivariant bridge trisection of $(\widetilde X,\widetilde\Ss)$.
\end{proposition}

\begin{proof}
	The fact that $\widetilde\TT = p^{-1}(\TT)$ is a bridge trisection of $(\widetilde X,\widetilde\Ss)$ was established in~\cite[Proposition~13]{LamMei_22_Bridge-trisections-in-rational}.
	It is immediate from the set-up that the pieces of $\widetilde\TT$ are $G$--invariant, so it remains to check that the action of $G$ on each sector is linearly parted.
	For this, consider a sector $(X_i,\Dd_i)$ of $\TT$.
	Choose a collection of neatly embedded 3--balls $\Bb_i\subset X_i$ such that each component of $X_i\setminus\Bb_i$ is a 4--ball that contains at most one component of $\Dd_i$.
	Let $\widetilde\Bb_i = p^{-1}(\Bb_i)$; we claim that the action of $G$ on $\widetilde X_i$ is linearly parted by $\widetilde\Bb_i$.
	Let $Z$ be a component of $\widetilde X_i\setminus\widetilde\Bb_i$.
	The restriction $p\vert_Z\colon Z\to p(Z)$ is a regular branched covering, and $p(Z)$ is a 4--ball component of $X_i\setminus\Bb_i$.
	It follows that $p\vert_Z$ is either a diffeomorphism, or a cyclic self-covering of the 4--ball, branched over a trivial disk.
	In either case, the induced action of $G$ on $Z$ is linear, and $p^{-1}(\Dd_i)$ is a collection of equivariantly boundary-parallel disks,  as desired.
\end{proof}

An action of a group $G$ on a compact, smooth manifold $X$ is called a \emph{branched covering action} if the quotient map $p\colon X\to X/G$ is a (regular) branched covering.
Trisections associated to branched coverings have been extensively discussed~\cite{BlaCahKju_24_Note-on-three-fold-branched-covers,CahKju_17_Singular-branched,CahMatRup_23_Algorithms-for-Computing-Invariants-of-Trisected,LamMei_22_Bridge-trisections-in-rational,LamMeiSta_21_Symplectic-4-manifolds-admit-Weinstein}.
However, these branched coverings are often either irregular or regular and cyclic.
In Example~\ref{ex:Q8-link} below, we will offer an extended example to illustrate the richness of equivariant trisections corresponding to more general regular branched coverings.

First, we give the following lemma and corollary which characterizes the possible singular (and hence, branching) sets of linearly parted branched covering actions on 1--handlebodies, which is of general use when discussing equivariant trisections associated to branched covering actions, and we will use when discussing branched covering actions in Theorem~\ref{thm:hyperelliptic} and Section~\ref{sec:classification}.

\begin{lemma}
\label{lem:branched_parted_sub}
	Let $(X,G)$ be an $n$--dimensional equivariant 1--handlebody with linearly parting $(n-1)$--ball-system $\Bb$ such that the quotient map $X \to X/G$ is a regular branched covering map.
	Let $S\subset X$ be the singular set.
	Then $S$ is an $(n-2)$--dimensional 1--handlebody neatly embedded in $X$ such that
	\begin{enumerate}
		\item The intersection $\mathcal B\cap S$ is a linearly parting $(n-3)$--ball-system for $S$, and
		\item Each disk component of $S\setminus \mathcal B$ is an equivariantly boundary-parallel $(n-2)$--disk in the $n$--ball of $X\setminus \mathcal B$ that contains it.
	\end{enumerate}
	In particular, $X/G$ is an $n$--dimensional $1$--handlebody, and the branching set $S/G$ is a disjoint union of $(n-2)$--dimensional 1--handlebodies.
\end{lemma} 

\begin{proof}
	We being by modifying $\Bb$ so that the induced action of $G$ on $\Bb$ is orientation-preserving:
	If there is some component $B$ of $\Bb$ for which some element of $G$ reverses orientation, replace $B$ with the two push-offs $\partial(\nu(B))\setminus\partial X$.
	In abuse of notation, we continue to denote the resulting ball-system $\Bb$.
	
	Let $Z$ be a 0--handle component of $X\setminus\Bb$.
	Let $G_Z$ be the stabilizer of $Z$, and let $\Bb_Z = Z\cap \Bb$.
	We will first show that $S_Z = Z\cap S$ is an equivariantly boundary-parallel $(n-2)$--disk on which the induced $G$--action is linear.
	Then, a similar argument will yield that, for each parting ball $B\in \Bb$, $B\cap S$ is an $(n-3)$--dimensional disk that is neatly embedded in $S$.
	From this, the claim about $S$ and parts (1) and (2) of the theorem follow.
	For the claims about $X/G$ and $S/G$, simply note that the quotient of a linearly parting ball-system is a parting ball-system for the quotient.

	Since the quotient map $Z\to Z/G_Z$ is a branched covering map, $Z/G_Z$ is a manifold, and since the action of $G_Z$ on $Z$ is linear, it follows that $Z/G_Z$ is an $n$-ball.
	By the forward direction of the~\ref{thm:ball_quotient}(1), which does not require the restriction to $n\leq 4$, we know that $G_Z$ acts on $Z$ as a rotation group.
		
	If $S_Z = Z\cap S$ is empty, then we are done.
	Otherwise, it must be that $G_Z$ is non-trivial and (since $Z\to Z/G_Z$ is a branched covering) $S_Z$ is a manifold.
	So, $G_Z$ contains a rotation, which fixes a disk $D_Z$, by definition.
	If $G_Z$ contained a rotation $r$ with $\Fix(r)\not=D_Z$, then $\Fix(r)\cup D_Z$ would be a non-manifold component of $S_Z$ (since they intersect at the origin of $Z$), a contradiction.
	Hence, $G_Z$ is cyclic and is generated by a rotation that fixes $D_Z$, so we have $S_Z = D_Z$, an equivariantly boundary-parallel $(n-2)$--disk in $Z$, as desired.
		
	Now consider some $(n-1)$--ball $B\subset \Bb$.
	By an argument similar to the one given above, the forward direction of the~\ref{thm:ball_quotient}(1)  implies that $B\cap S$ an $(n-3)$--disk that is neatly embedded in $S$.
\end{proof}

If $S$ is an $(n-2)$--dimensional invariant handlebody embedded in an $n$--dimensional handlebody $X$, then we can consider $S$ to be a ``sub-handlebody'' if condition (1) above holds. In this situation, we can apply~\cite[Lemma~3.7]{MeiSco_LP} simultaneously to a $0$--handle $Z$ and the $0$--handle $S_Z = S\cap Z$ of $S$.
Applying these techniques as was done in~\cite[Corollary~3.8]{MeiSco_LP}, we have the following.

\begin{corollary}
\label{cor:branched_auto_parted}
	Let $(X,G)$ be an $n$--dimensional linearly parted 1--handlebody such that $X/G\cong B^n$ and the quotient map $X\to B^n$ is a regular branched covering map.
	Let $S\subseteq X$ be the singular set.
	Then the branching set $S/G$ is an $(n-2)$--dimensional trivial disk--tangle in $B^n$, $S$ is an $(n-2)$--dimensional equivariantly boundary-parallel disk--tangle in $X$, and $(X, S)$ is linearly parted as a pair.
\end{corollary}

\begin{proof}
	Let $\mathcal B$ be a linearly parting $(n-1)$--ball-system for $X$.
	As in Lemma~\ref{lem:branched_parted_sub}, modify $\Bb$ so that each component has orientation-preserving stabilizer. 
	By Lemma~\ref{lem:branched_parted_sub}, $S$ is linearly parted by $\mathcal B\cap S$, and each disk component of $S\setminus\Bb$ is equivariantly boundary-parallel in its $0$--handle.

	Consider the quotient $\Bb^* = \mathcal B/G$, which is a parting ball-system for $B^n$.
	In each $0$--handle $Z^* = Z/G$ of $B^n\setminus\Bb^*$, either we see no points of the branching set $S^* = S/G$, or we see a boundary-parallel disk $S_Z/G_Z$, as in the proof of Lemma~\ref{lem:branched_parted_sub}.
	We now mirror the proof of~\cite[Corollary~3.8]{MeiSco_LP}.

	The graph $\Gamma^*$ associated to the handle decomposition corresponding to $\Bb^*$ must be a tree, since $B^n$ is contractible. 
	If this handle-decomposition consists of a single $0$--handle, we're done.
	Otherwise, $\Gamma^*$ has a leaf vertex $Z^*_2$ adjacent to a vertex $Z^*_1$ along and edge $B^*\in \mathcal B^*$.
	We now repeat the following inductive procedure to simplify the handle-structure on $B^n$.
	
	If $Z^*_2$ contains no disk of $S^*$, contract it into $Z^*_1$ along $B$.
	If $Z^*_2$ contains a disk $D_2^*\subseteq S^*$, then either $D_2^*$ intersects the cocore $B^*$ or it does not.
	
	If $D_2^*\cap B^* = \varnothing$, then set $(Z^*_2,D_2^*)$ aside as a copy of $(B^n, D^{n-2})$ occurring as a boundary-connect-summand of $(B^n, S^*)$.
	If $D_2^*$ meets $B^*$, then since $S\cap B$ was an $(n-3)$--disk neatly embedded in $S$, $D_2^*\cap B^*$ must be an $(n-3)$--disk in the boundary of $D^*$.
	In this case, $Z_1^*$ also contains a disk $D_1^*$ of $S^*$ and we can contract $(Z_2^*, D_2^*)$ into $(Z_1^*,D_1^*)$.

	At the end of this procedure, we will have decomposed $(B^n, S^*)$ into the boundary connected sum of some number of copies of $(B^n, D^{n-2})$, and so the branching set $S^*$ in $B^n$ is a boundary-parallel $(n-2)$--disk tangle.
	Let $\Delta^*$ denote a collection of bridge-balls for $S^*$, and let $(\Bb')^*$ denote a collection of neatly embedded $(n-1)$--balls disjoint from $\Delta^*$ that isolate the components of $S^*$ in $B^4$.
	Let $\Bb'$ and $\Delta$ denote the lifts of $(\Bb')^*$ and $\Delta^*$, respectively.
	Then, $\Bb'\cup\Delta$ exhibits that $(X,S)$ is linearly parted as a pair with respect to the $G$--action.
\end{proof}

This corollary is especially useful in the case of a cyclic group $\Z_n$ acting on a $3$--dimensional handlebody $X$.
In this setting, the stabilizer of each $0$--handle of a linear parting is cyclic and acts on each $0$--handle by rotation around an axis.
Thus, the quotient map $X\to X/\Z_n$ is a branched covering map.
This gives the following immediate corollary, which we'll use in Theorem~\ref{thm:genus_two}.

\begin{corollary}
\label{coro:cyclic_3d}
	Let $\Z_n$ act on a $3$--dimensional handlebody $X$ so that $X/\Z_n$ is a $3$--ball.
	The quotient map $X\to X/\Z_n$ is a branched covering, and the branched set (hence, also the singular set) is a trivial tangle.
\end{corollary}

We now give examples of non-trivial branched covering actions and their equivariant trisections.

\begin{example}
\label{ex:Q8-link}
	Throughout this example, let $N_e(S)$ denote the disk-bundle over the surface $S$ with Euler number $e$; when $\Ss$ is a surface in a 4--manifold with normal Euler number $e(\Ss)=e$, we identify $N_e(S)$ with a tubular neighborhood of $\Ss$.
	
	To begin, consider the unknotted projective plane $P_+$ in $S^4$ with normal Euler number 2; see~\cite{JosMeiMil_22_Bridge-trisections-and-classical} for complete details and remarks on conventions.
	It turns out that $S^4\setminus N_2(P_+)$ is diffeomorphic to $N_{-2}(P_-)$, where $P_-$ is the other unknotted projective plane; see~\cite[Chapter~6]{GomSti_99_4-manifolds-and-Kirby} for discussions of the handle-decompositions of $N_e(S)$ and $S^4\setminus N_e(\Ss)$.
	Let $Q_8$ denote the quaternion group of eight elements, which can be presented as
	$$Q_8 = \langle a, b\,|\,\bar baba, \bar abab\rangle,$$
	where $a\mapsto i$ and $b\mapsto j$ gives the standard description of the group.
	Let $Y_{Q_8}$ denote the circle-bundle over $\RP^2$ with Euler number 2, which has fundamental group $Q_8$ and coincides with $\partial N_2(P_+)$.
	Summarizing, we have 
	$$S^4 = N_2(P_+)\bigcup_{Y_{Q_8}}N_{-2}(P_-),$$
	and the link $\Ll = P_+\cup P_-$ satisfies $\pi_1(S^4\setminus\nu(\Ll))\cong Q_8$, since $S^4\setminus\nu(\Ll)\cong Y_{Q_8}\times I$.
	Note that this decomposition of $S^4$ is not a double.
	In~\cite{JosMeiMil_22_Bridge-trisections-and-classical}, the link $(S^4,\Ll)$ is discussed from the bridge trisection perspective: there, a tri-plane diagram is given in Figure~10, and the fundamental group is computed in Example~4.4.
	This tri-plane diagram is shown here in Figure~\ref{fig:Q8_covers}(\textsc{a}): the left, darker component is $P_+$, and the right, lighter component is $P_-$.
	A corresponding shadow diagram is shown in Figure~\ref{fig:Q8_covers}(\textsc{b}), and the $Q_8$--coloring of the meridians at the bridge sphere is indicated: each label is the element of $Q_8$ associated to a counter-clockwise meridian of that bridge point.
	
	The link $(S^4,\Ll)$ admits a regular branched cover for each epimorphism $\rho\colon Q_8\twoheadrightarrow G$ for a finite group $G$.
	For each such branched covering, we get a $G$--equivariant trisection of the covering manifold whose singular set is in bridge trisected position.
	We now analyze all such coverings; our conclusions are summarized in Figure~\ref{fig:Q8_table} and illustrated in Figure~\ref{fig:Q8_covers}.
	Figure~\ref{fig:Q8_covers}(\textsc{c}) shows the bridge sphere cut open along the shadows $\frak a_1$ of the tangle $\Tt_1$.
	This is the fundamental domain for the diagrams of the ensuing coverings.

\begin{figure}[ht!]
	\centering
	\includegraphics[width = .9\textwidth]{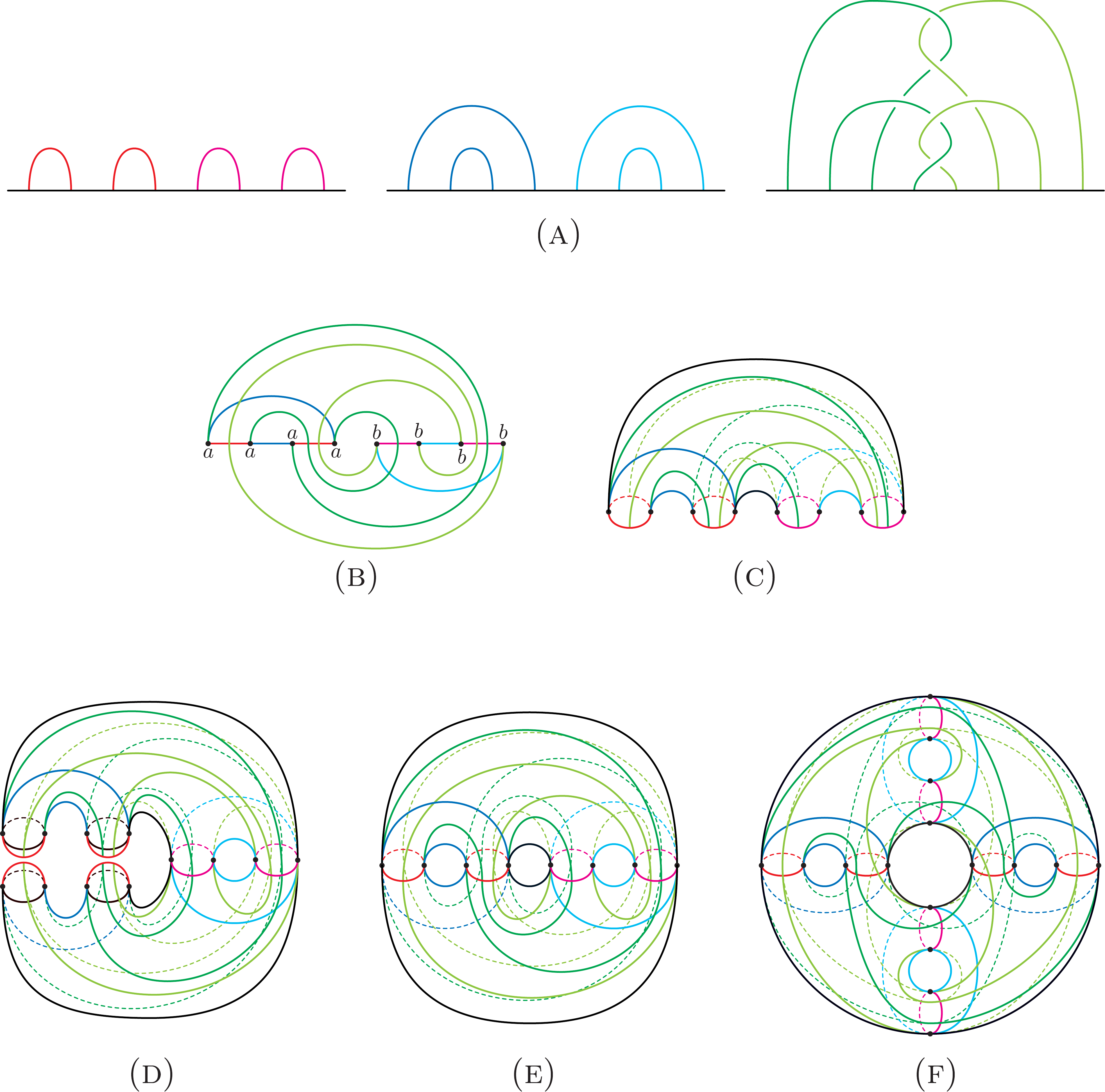}
	\caption{Diagrams for the link $\Ll$ and some of its branched covers}
	\label{fig:Q8_covers}
\end{figure}
	
	If $G\cong\Z_2$, there are three coverings, which are determined by the possible values of $\rho(a)$ and $\rho(b)$.
	If $\rho(a) = 0$ and $\rho(b)=1$, then the corresponding cover is branched along $P_-$.
	In this case, the covering space is $\CP^2$, and the lifts of $P_+$ and $P_-$ are $\Cc_2$ and $\RP^2$, respectively, where $\Cc_2$ denotes the complex curve of degree two and $\RP^2$ is the standard real projective plane.
	Note that the decomposition of $S^4$ as two disk-bundles lifts: we get
	$$\CP^2 = N_4(\Cc_2)\bigcup_{L(4,1)}N_{-1}(\RP^2),$$
	where $L(4,1)$ is the lens space that double covers $Y_{Q_8}$.
	The $\Z_2$--equivariant (bridge) trisection of the covering space is shown in Figure~\ref{fig:Q8_covers}(\textsc{d}).
	Note that this is the standard genus-one trisection of $\CP^2$; the deck transformation is the hyperelliptic involution.
	The lift $\Cc_2$ is in $4$--bridge position, as indicated by the darker-shaded shadow arcs in the diagram.
	The lift $\RP^2$ is in $2$--bridge position; the shadow arcs are contained in the lighter-shaded cut curves.
	Reversing orientation everywhere results in the description of the covering space corresponding to $\rho(a)=1$ and $\rho(b)=0$.
	
	When $\rho(a) = \rho(b) = 1$, the covering space is $\Ss^*(\RP^3)$, the twisted-spin of $\RP^3$, which can also be described as
	$$\Ss^*(\RP^3) = N_1(\widetilde P_+)\bigcup_{L(4,1)}N_{-1}(\widetilde P_-),$$
	where $\widetilde P_\pm$ is the lift of $P_\pm$.
	This decomposition is a double, since every automorphism of $L(4,1)$ extends across $N_{\pm 1}(\RP^2)$~\cite[Corollary~2.2]{FinSte_97_Rational-blowdowns-of-smooth-4-manifolds}.
	However, it is not an equivariant double, since the quotient is not a double.
	This could have implications for an equivariant version of the rational blowdown.
	The $\Z_2$--equivariant  $(3,1)$--trisection for $(\Ss^*(\RP^3),\widetilde P_+\cup\widetilde P_-)$ is shown in Figure~\ref{fig:Q8_covers}(\textsc{e}).
	Each lift $\widetilde P_\pm$ is in $(2,1)$--bridge trisected position, with the shadow arcs contained in the cut curves.
	The action is given by extending the hyperelliptic involution of the central surface.
	
	Next, consider the representation $\rho\colon Q_8\to \Z_2\oplus\Z_2$ given by $\rho(a) = (1,0)$ and $\rho(b) = (0,1)$.
	The covering space is the 2--fold universal covering of $\Ss^*(\RP^3)$ and is described as
	$$S^2\times S^2 = N_2(\Cc_{(1,1)})\bigcup_{\RP^3}N_{-2}(\Cc_{(1,-1)}),$$
	where $\Cc_{a,b}$ is the complex curve representing $(a,b)$ in $H_2(S^2\times S^2)\cong\Z\oplus\Z$.
	Again, the decomposition is the double.
	The $(\Z_2\oplus\Z_2)$--equivariant $(5,1)$--trisection is shown in Figure~\ref{fig:Q8_covers}(\textsc{f}).
	Each curve $\Cc_{(1,\pm1)}$ is in $(4,2)$--bridge trisected position, with the shadow arcs contained in the cut curves.
	(To see that the singular set is as claimed, note that the 2--fold universal covering map $S^2\times S^2\to \Ss^*(\RP^3)$ is given by $(x,y)\mapsto(-x,-y)$ and the involutions with quotient space $\CP^2$ and $\overline\CP^2$ are, respectively, $(x,y)\mapsto(y,x)$ and $(x,y)\mapsto(-y,-x)$.)
	
	Finally, consider the identity epimorphism $\rho\colon Q_8\to Q_8$.
	This covering space is the 2--fold cover of $S^2\times S^2$, branched along $\Cc_{(1,1)}\cup\Cc_{(1,-1)}$ and since the double structure is preserved, is described as
	$$\CP^2\#\overline\CP^2 = N_1(\CP^1)\bigcup_{S^3}N_{-1}(\overline\CP^1).$$
	The corresponding $(17,5)$--trisection is not pictured, but each projective line is in $(4,2)$--bridge trisected position.
\end{example}

\begin{figure}[ht!]
\centering
\adjustbox{scale=.7,center}{%
\begin{tikzpicture}[mybox/.style={draw, inner sep=5pt}]
\node[mybox] (box){%
\begin{tikzcd}[every arrow/.append style={dash}]
G = \im(\rho) & & (X,\Ss) & & \ker(\rho) = \pi(\Ss) \\
\hline \\
1\ar[d]&
&
\begin{array}{c} \displaystyle N_2(P_+)\bigcup_{Y_{Q_8}}N_{-2}(P_-) \\ (S^4,P_+\cup P_-) \end{array}\ar[d]\ar[dl]\ar[dr]
&& Q_8\ar[d]
\\
\Z_2\ar[d]
&
\begin{array}{c} \displaystyle N_4(\Cc_2)\bigcup_{L(4,1)}N_{-1}(\RP^2) \\ (\CP^2,\Cc_2\cup \RP^2) \end{array} \ar[dr]
&
\begin{array}{c} \displaystyle N_1(\widetilde P_+)\bigcup_{L(4,1)}N_{-1}(\widetilde P_-) \\ (\Ss^*(\RP^3),\widetilde P_+\cup\widetilde P_-)) \end{array}
 \ar[d]
&
\begin{array}{c} \displaystyle N_1(\RP^2)\bigcup_{L(4,1)}N_{-4}(\Cc_{-2}) \\ (\overline\CP^2,\overline\RP^2\cup \Cc_{-2}) \end{array}
\ar[dl]
&\Z_4\ar[d]
\\
\Z_2\oplus\Z_2\ar[d]&
&
\begin{array}{c} \displaystyle N_2(\Cc_{(1,1)})\bigcup_{\RP^3}N_{-2}(\Cc_{(1,-1)}) \\ (S^2\times S^2,\Cc_{(1,1)}\cup \Cc_{(1,-1)}) \end{array}
\ar[d]
&&\Z_2\ar[d]
\\
Q_8&
&
\begin{array}{c} \displaystyle N_1(\CP^1)\bigcup_{S^3}N_{-1}(\overline\CP^1) \\ (\CP^2\#\overline\CP^2,\CP^1\cup\overline\CP^1) \end{array}
&&1
\end{tikzcd}
};
\end{tikzpicture}
}
\caption{The branched coverings of the link $\Ll\subset S^4$ of projective planes corresponding to the representations $\rho\colon\pi_1(S^4\setminus\nu(\Ll))\to Q_8$}
\label{fig:Q8_table}
\end{figure}
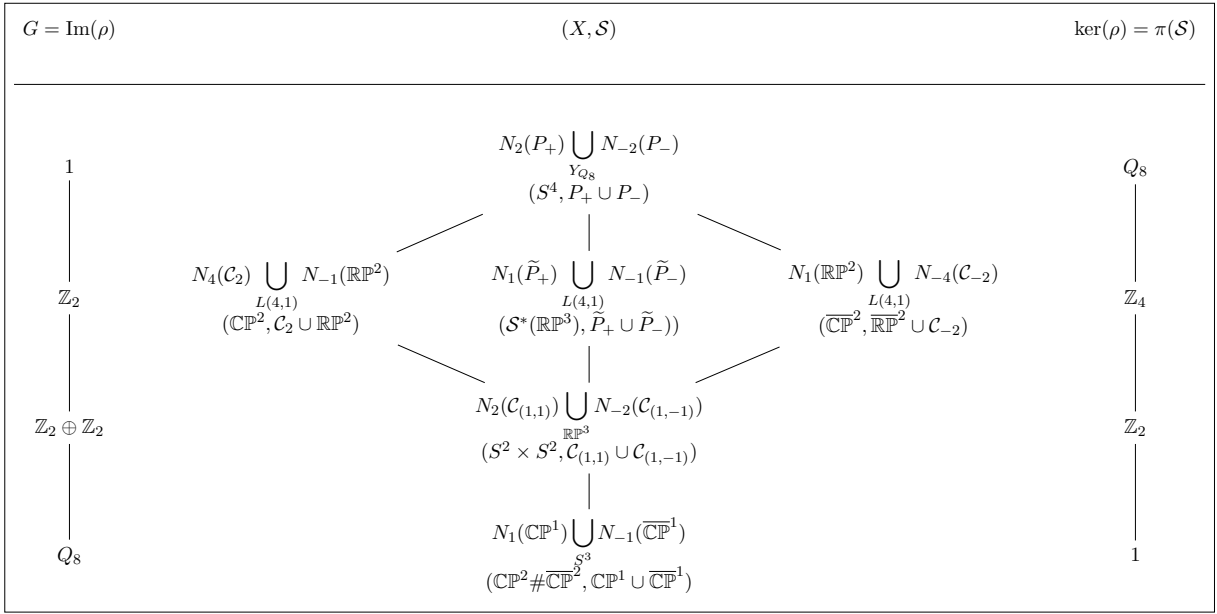

\begin{remark}
	The $Q_8$ action on $\CP^2\#\overline\CP^2$ described above is an equivariant connected sum, as indicated by the decomposition
	$$\CP^2\#\overline\CP^2 = N_1(\CP^1)\bigcup_{S^3}N_{-1}(\overline\CP^1).$$
	Considering the action on just $N_1(\CP^1)$, we see that the $Q_8$ action respects the bundle structure in the sense that the base $\CP^1$ is left invariant and fibers are permuted: for any point $b\in\CP^1$ and any $g\in Q_8$, we have $g\cdot D^2_b = D^2_{g\cdot b}$.

	In Subsection~\ref{subsec:PU(3)_trisections} below, we analyze and trisect some actions on $\CP^2$ involving finite subgroups of the subgroup $U(2)<PU(3)$, which we identify with $\{[1\oplus A]\,:\,A\in U(2)\}$.
	These actions also leave $\CP^1\subseteq \CP^2$ invariant and permute fibers, thus also respect the bundle structure on $N_1(\CP^1)$.
	One of these actions is a $Q_8$--action.

	Surprisingly, these two $Q_8$--actions on $N_1(\CP^1)$ are not the same.
	The action described above and obtained via branched covering has elements that act on the base $\CP^1$ by the antipodal map, whereas the $Q_8$--action described below and obtained from $U(2)<PU(3)$ does not.
	This phenomenon arises from the fact that $Q_8$ has two representations on $\R^4 = \C^2$: a two-dimensional complex one, realizing $Q_8$ as a subgroup of $U(2)$, and a four-dimensional real representation, which is not complex, arising from the identification $\R^4 = \mathbb H$, which we obtained from the branched covering above.
	This real representation no longer commutes with complex multiplication by elements $u\in \C^\times$, but for any $z\in \C^2$ and $g\in Q_8$, under the real representation, $g\cdot(uz)$ is either $u(g\cdot z)$ or $\overline{u}(g\cdot z)$.
	Since $\overline{u} \in \C^\times$, the representation still descends to a well-defined action on $\C^2/\C^\times = \CP^2$, which we see above.

	One last remark about this pair of actions is that the quotient of the preserved bundle structure gives a disk bundle over a surface whose boundary circle bundle is $S^3/Q_8 = Y_{Q_8}$.
	This gives two different presentations of $Y_{Q_8}$: first, as a Seifert fiber space $S^2(2,2,4)$, and second, as the circle bundle over $\RP^2$ with normal Euler number $2$.
	Capping off each circle-fiber with disks, the latter naturally extends to $N_2(\RP^2)$, as discussed above, while the former extends to a non-manifold object that is a sort of a disk-bundle over the orbifold $S^2(2,2,4)$.
\end{remark}

We conclude with an interesting observation about $Q_8$--equivariant trisections of $\CP^2\#\overline\CP^2$.
Recall that a trisection is \emph{reducible} if it is the connected sum of two trisections, neither of which is genus-zero; see~\cite{MeiSchZup_16_Classification-of-trisections} for more details.
We say that a $G$--equivariant trisection is \emph{reducible} if there is some such connected sum decomposition that is $G$--invariant.
A trisection (equivariant or not) is \emph{irreducible} if it is not reducible.

\begin{corollary}
\label{cor:Q8}
	Every $Q_8$--equivariant trisection is irreducible.
\end{corollary}

The proof is immediate from the two two propositions that follow.

\begin{proof}
	Suppose that $\TT = \TT^1\#\TT^2$ is a reducible, $Q_8$--equivariant trisection, with invariant reducing 3--sphere $P$.
	By Proposition~\ref{prop:reducing_curve}, the induced action of $Q_8$ on $P$ leaves a genus-one Heegaard splitting invariant.
	By Proposition~\ref{prop:Q8_S3}, the induced action of $Q_8$ on $P$ is free.
	Since there is no free $Q_8$--action on the torus (see Lemma~\ref{lem:genus_1_maximal_handlebody_group}), this is a contradiction. 
\end{proof}

\begin{proposition}
\label{prop:reducing_curve}
	Suppose that $\TT = \TT^1\#\TT^2$ is a reducible, $G$--equivariant trisection, with invariant reducing sphere $P$.
	Then, the induced $G$--action on $P$ respects a genus-one Heegaard splitting of $P$.
\end{proposition}

\begin{proof}
	By definition, $P\cap\Sigma = \alpha$ is a single curve that bounds a disk $D_i = P\cap H_i$ for each $i\in\Z_3$.
	The disks $D_i$ are pages of the genus-zero open book decomposition of $P$, and $\alpha$ is the binding.
	Since $G$ leaves $\alpha$ invariant, $G$ leaves the genus-one Heegaard splitting obtained by thickening $\alpha$ invariant.
\end{proof}

\begin{proposition}
\label{prop:Q8_S3}
	Every action of $Q_8$ on $S^3$ is free.
\end{proposition}

\begin{proof}
	Let $Q_8$ act on $S^3$.
	Suppose for a contradiction that $\Stab(x)$ is non-trivial for some $x\in S^3$.
	Every non-trivial subgroup of $Q_8$ is cyclic; so the fixed-point set of $\Stab(x)$ is an unknot $K$, by the Smith Conjecture.
	Since every element $g\in Q_8$ has some power $g^n = -1$ and $-1$ acts freely, every element of $Q_8$ acts freely on $S^3\setminus\nu(K)$; hence $Q_8$ acts freely on on $S^3\setminus\nu(K)$, a solid torus.
	However, as before, no such action exists, by Lemma~\ref{lem:genus_1_maximal_handlebody_group}.	 
\end{proof}

\subsection{Hyperelliptic trisections}
\label{subsec:hyperelliptic_trisections}

Given a genus--$g$ surface $\Sigma$, recall that an involution $\tau\colon\Sigma\to\Sigma$ is called \emph{hyperelliptic} if it induces $-\id$ on $H_1(\Sigma)$.
Equivalently, $\tau$ fixes $2g+2$ points and can be described as rotation of a $(4g+2)$--gon representation of $\Sigma$ through $\pi$ radians; any two hyperelliptic involutions are conjugate in $\Mod(\Sigma)$~\cite{FarMar_12_A-primer-on-mapping-class-groups}.
A $\Z_2$--equivariant trisection is called \emph{hyperelliptic} if the induced $\Z_2$--action on $\Sigma$ is generated by a hyperelliptic involution.
A thorough discussion of a related class of involutions of 4--manifolds, called conjugation manifolds, is given by Hambleton and Hausmann~\cite{HamHau_11_Conjugation-spaces-and-4-manifolds}.
Here, we pursue a trisection-theoretic perspective, starting with the following consequence of Corollary~\ref{cor:branched_auto_parted}.

\begin{theorem}
\label{thm:hyperelliptic}
	Let $G\cong\Z_2$, and let $X$ be a $G$--manifold.
	If $X$ admits a $G$--equivariant hyperelliptic trisection $\TT$, then $X/G\cong S^4$ and $\Fix(G)$ is in equivariant bridge trisected position with respect to $\TT$.
	Conversely, if $X/G\cong S^4$, then $X$ admits a $G$--equivariant hyperelliptic trisection.
	In either case, the quotient map $\TT\to\TT/G$ is a 2--fold branched covering.
\end{theorem}

\begin{proof}
	First, suppose that $X$ admits a hyperelliptic trisection $\TT$.
	Pick any diagram $\DD$ for $\TT$, and consider the quotient $\DD/G$.
	This is a sphere, as well as a trisection surface for the quotient $X/G$.
	It follows that $X/G\cong S^4$.
	
	The $G$--action is not free, since $S^4$ is simply connected, and cannot fix an isolated point since $S^4$ is a manifold.
	Furthermore, since $G$ has no nontrivial subgroups, the action is semifree, and thus we have that $\Sing(G) = \Fix(G)$ is a surface-link.
	It follows that the quotient map $X\to X/G$ is a 2--fold branched covering.
	
	By Corollary~\ref{cor:branched_auto_parted}, applied to each sector and to each handlebody, $\Fix(G) = \Sing(G)$ is in equivariant bridge trisected position, as desired.
	
	For the converse, suppose $X$ is a $G$--manifold with $G\cong\Z_2$ and $X/G\cong S^4$.
	By the above argument, the quotient map $X\to X/G$ is a 2--fold branched covering.
	Choose a genus-zero bridge trisection $\TT^*$ for $(S^4,\Fix(G)/G)$, and let $\TT$ be the pre-image of this trisection under the quotient map, which is a trisection by Proposition~\ref{prop:branched_cover}.
	Then, $\TT$ is a hyperelliptic $G$--equivariant trisection on $X$, and the quotient map $\TT\to\TT/G$ is a 2--fold branched covering.
\end{proof}

As a first application of hyperelliptic trisections, we prove the following classical theorem, which was proved independently by Massey~\cite{Mas_73_The-quotient-space-of-the-complex-projective}, Kuiper~\cite{Kui_74_The-quotient-space-of-bf-CP2-by-complex-conjugation}, and Arnol'd~\cite{Arn_88_The-branched-covering-bf-Crm-P2to}; see~\cite{HamHau_11_Conjugation-spaces-and-4-manifolds} for the more general context.
We give a proof making use of the genus-one hyperelliptic involution of $\CP^2$.

\begin{theorem}
\label{thm:massey-kuiper}
	Let $\iota\colon\CP^2\to\CP^2$ be the complex conjugation given by
	$$\iota([z_1\colon z_2\colon z_3]) = [\overline z_1\colon \overline z_2\colon \overline z_3].$$
	The quotient space $\mathbb{CP}^2/\langle \iota\rangle$ is diffeomorphic to $S^4$.
\end{theorem} 

\begin{proof}
	Consider the moment map $\mu\colon \mathbb{CP}^2\to \mathbb{R}^2$ given by 
	$$\mu([z_1\colon z_2\colon z_3]) = \left(\dfrac{|z_1|}{|z_1|+|z_2|+|z_3|}, \dfrac{|z_2|}{|z_1|+|z_2|+|z_3|}\right).$$
	The image of $\mu$ is the convex hull of the positive unit vectors in $\R^2$, and the pre-image $(x,y)\in\im(\mu)$ is either a torus, a circle, or a point, depending on whether $(x,y)$ is in the interior of the triangle, in the interior of an edge, or a vertex, respectively.
	Regardless, $\mu^{-1}(x,y)$ is $\langle\iota\rangle$--invariant.
	For example, for torus fibers, we have
	$$\mu^{-1}(x,y) = \mu^{-1}\left(\frac{a}{a+b+1}, \frac{b}{a+b+1}\right) = [ae^{i\theta}\colon be^{i\psi}\colon 1],$$
	for some fixed positive real $a,b$, where $\theta, \psi\in S^1$. 
	The action of $\iota$ on this fiber negates both $\theta$ and $\psi$, which is the hyperelliptic involution on the torus.
	
	Cutting $\im(\mu)$ into three pieces, we get a genus-one trisection of $\CP^2$; see Figure~\cite[Figure~4]{GayKir_16_Trisecting-4-manifolds}.
	This gives a hyperelliptic $\langle\iota\rangle$--invariant trisection.
	Furthermore, the $\langle\iota\rangle$--action on each $B^4$ sector of this trisection is the equivariant cone on its boundary action, since each sector of $\im(\mu)$ is the cone on the arc along its boundary whose pre-image is $H_i\cup\overline{H_j}$.
	This invariant trisection is, therefore, actually a hyperelliptic $\langle\iota\rangle$--equivariant trisection.
	By Theorem~\ref{thm:hyperelliptic}, $\CP^2/\langle\iota\rangle\cong S^4$, as desired.
\end{proof}

Next, we analyze the components of a $(\Z_2\oplus\Z_2)$--action on $S^2\times S^2$ using the genus-two trisection.

\begin{example}
\label{ex:s2xs2}
	Viewing $S^2$ as the unit 2--sphere in $\R^3$ with cylindrical coordinates $(\theta,z)$, define $\mu\colon S^2\times S^2\to \mathbb R^2$ by
	$$\mu(\theta_1,z_1,\theta_2,z_2) = (z_1, z_2).$$
	The image of $\mu$ is the unit square $[-1,1]\times[-1,1]$, and the pre-image of each point in the image is a torus, a circle, or a single point, as was the case with the moment map for $\CP^2$ above.
	See Figure~\ref{fig:s2xs2_moment}.
	In particular, for $(z_1,z_2)$ in the interior of the square, the fiber $F_{(z_1,z_2)}$ is parametrized by $(\theta_1,\theta_2)$.
	The map $\mu$ is equivariant with respect to the component-wise $(S^1\times S^1)$--action defined by rotating each factor around its $z$--axis along $\theta_i$.

\begin{figure}[ht!]
	\centering
	\includegraphics[width = 0.9\linewidth]{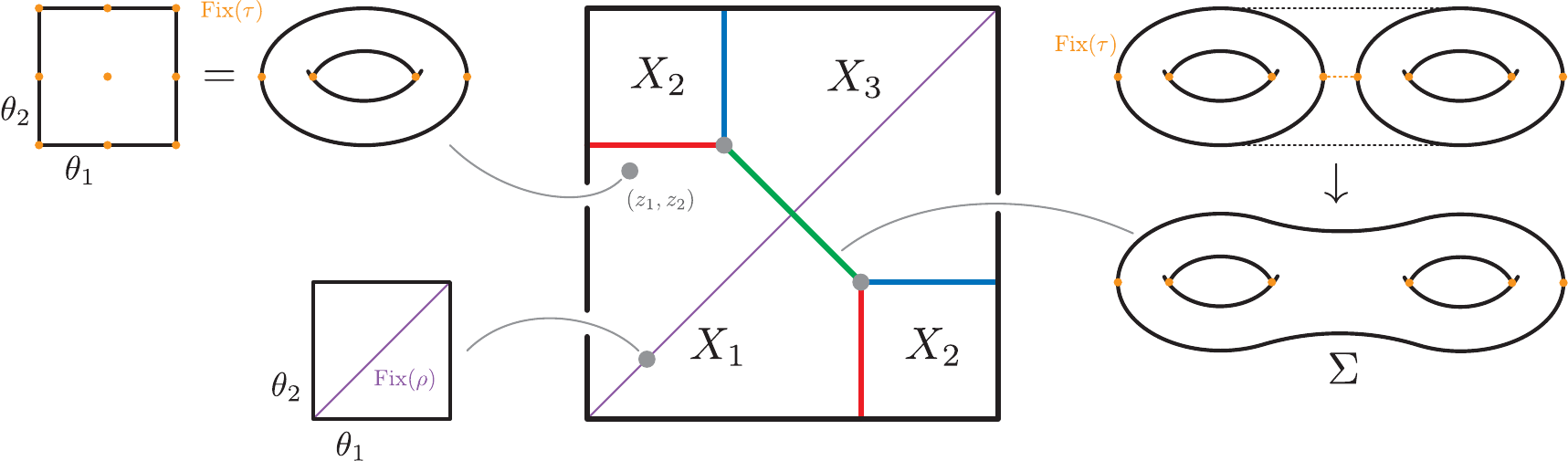}
	\caption{The genus-two trisection of $S^2\times S^2$ coming from the moment map}
	\label{fig:s2xs2_moment}
\end{figure}

	Let $\tau\colon S^2\times S^2\to S^2\times S^2$ be given by
	$$\tau(\theta_1,z_1,\theta_2,z_2) = (-\theta_1,z_1,-\theta_2,z_2),$$
	so $\tau$ is reflection of each $S^2$ across the $xz$--plane.
	Note that $\mu\circ\tau = \mu$.
	For $(z_1,z_2)$ in the interior of the unit square, $\Fix(\tau)\cap F_{(z_1,z_2)}$ is the four points where $\theta_i\in\{0,\pi\}$.
	The four points collapse to pairs of points over the interiors of the edges of the square and to single points over the vertices.
	Since $\tau$ is the product of a reflection in each factor of $S^2\times S^2$, $\Fix(\tau)$ is the product of the fixed-point sets in each factor; that is, $\Fix(\tau) = S^1\times S^1$.
	
	Let $\rho\colon S^2\times S^2\to S^2\times S^2$ be given by 
	$$\rho(\theta_1,z_1,\theta_2,z_2) = (\theta_2,z_2,\theta_1,z_1),$$
	so $\rho$ swaps the $S^2$--factors.
	Note that $\mu\circ\rho = \rho^*\circ\mu$, where $\rho^*(z_1,z_2) = (z_2,z_1)$ is reflection across the diagonal of the unit square.
	So, $\Fix(\rho)$ is the diagonal 2--sphere defined by $z_1=z_2$ and $\theta_1=\theta_2$.
	The fixed-sets of $\tau$ and $\rho$ are indicated in gray in Figure~\ref{fig:s2xs2_moment}.

	Let $G = \langle\tau,\rho\rangle$.
	We now construct a $G$--equivariant trisection $\TT$ of $S^2\times S^2$ that will have the property that it is hyperelliptic (in $\tau$).
	Following Gay and Kirby~\cite[Section~2]{GayKir_16_Trisecting-4-manifolds}, consider the decomposition of the unit square shown in Figure~\ref{fig:s2xs2_moment}.
	The pre-image of the short, negatively sloped arc $A$ in the center is
	$$\bigcup_{(z_1,z_2)\in A} F_{(z_1,z_2)}\cong T^2\times A.$$
	The sectors $X_i$ are the pre-images of the labeled regions of the unit square, with the modification that we borrow a neighborhood of $(0,0)\times A$ from $T^2\times A$ and donate it to $X_2$.
	This has the effect of turning $T^2\times\partial A$ into $\Sigma$, a genus-two surface, as indicated on the right side of Figure~\ref{fig:s2xs2_moment}.
		
	The arc $(0,0)\times A$ is invariant under $\rho$ and fixed by $\tau$, so it is clear the corresponding trisection (see~\cite{GayKir_16_Trisecting-4-manifolds} for details) is $G$--invariant.
	It remains to show that the sectors of $\TT$ are linearly parted and that the action of $\tau$ is induced by a hyperelliptic involution of $\Sigma$.
	The action of $G$ on $X_1$ and $X_3$ is the cone of the action of $G$ on $Y_1$ and $Y_3$, respectively, so the actions are linear on these 4--balls.
	For $X_2$, let $B$ be the neighborhood of $(0,0)\in F_{(0,0)}$ that was donated to $X_2$.
	Then, $X_2\setminus B$ is two 4--balls: the pre-images of the regions of the unit square labeled with $X_2$.
	These are interchanged by $\rho$, and the action of $\tau$ on each is the cone of the action of $\tau$ on its boundary, so $B$ is a linearly parting ball-system for the $G$--action on $X_2$.
	Since $\tau$ fixed four points on each of the two tori that were summed to get $\Sigma$, and since the summation occurred at two of these fixed points, we see that $\tau$ fixes six points on $\Sigma$.
	It follows that $\tau$ is hyperelliptic, as desired.
		
	The diagram shown in Figure~\ref{fig:s2xs2_involutions}(\textsc{a}) corresponds to the trisection $\TT$~\cite{GayKir_16_Trisecting-4-manifolds}.
	The separating curve separating the left half of the surface from the right half is rotated by $\tau$ and reflected by $\rho$.
	It follows that these involutions act on the surface as rotation through $\pi$ radians according to the indicated axes.
	Let $\eta =\tau\circ\rho$.

\begin{figure}[ht!]
	\centering
	\includegraphics[width = 0.8\linewidth]{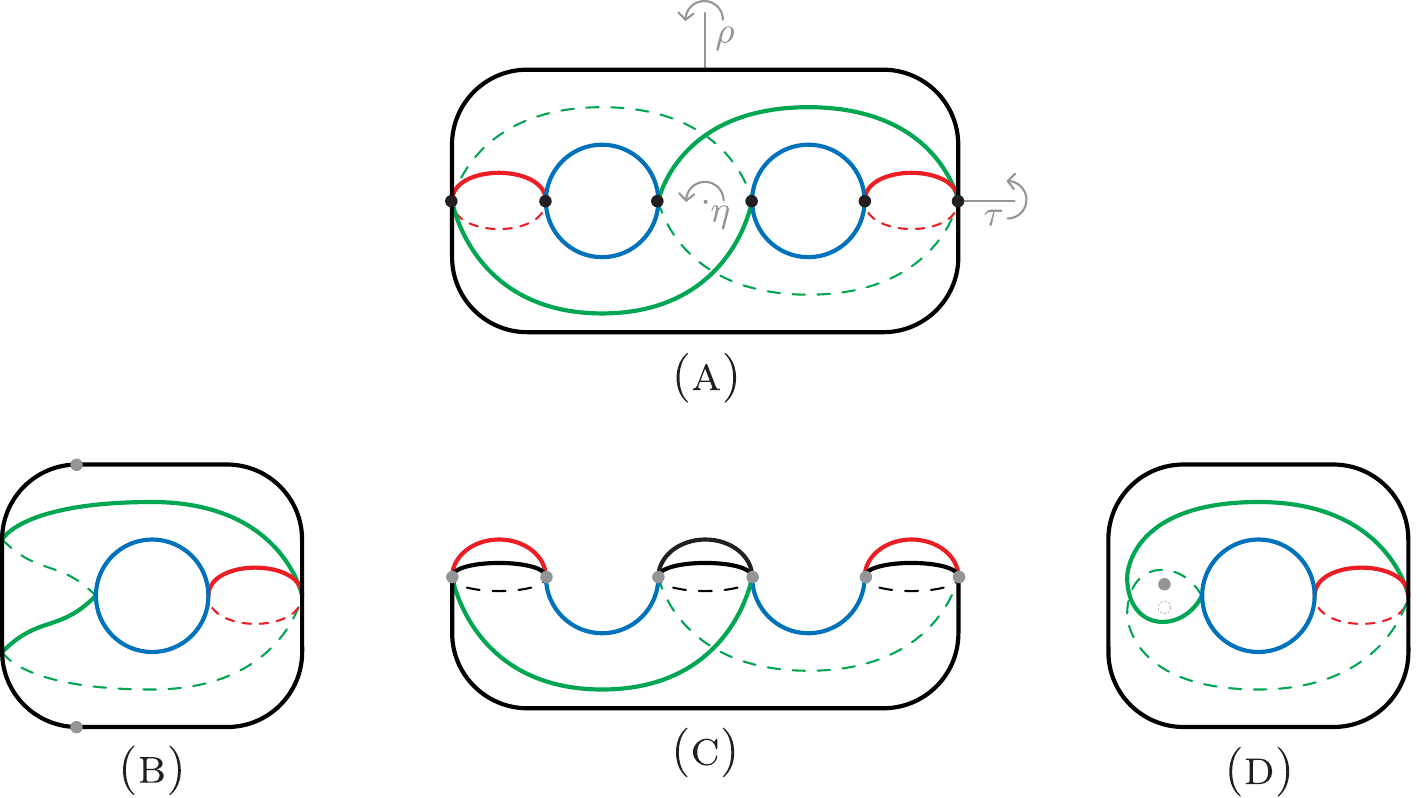}
	\caption{Quotients of involutions of the genus-two trisection of $S^2\times S^2$}
	\label{fig:s2xs2_involutions}
\end{figure}
	
	By Theorem~\ref{thm:hyperelliptic}, we know that $S^2\times S^2/\langle\tau\rangle\cong S^4$, the quotient map is a branched covering, and the fixed-point set is in bridge trisected position.
	More generally, by Theorem~\ref{thm:quotient_trisection}, the quotient of $\TT$ by any one of $\langle\rho\rangle$, $\langle\tau\rangle$, or $\langle\eta\rangle$ will be a trisection.
	Since the genus is low, we can identify the corresponding quotient of $S^2\times S^2$:
	\begin{enumerate}
		\item $S^2\times S^2/\langle\rho\rangle\cong \CP^2$ and $\Fix(\rho)/\langle\rho\rangle$ is the conic $\Cc_2\subset\CP^2$. A doubly pointed trisection diagram for $(\CP^2,\Cc_2)$ is shown in Figure~\ref{fig:s2xs2_involutions}(\textsc{b}).
		\item $S^2\times S^2/\langle\tau\rangle\cong S^4$ and $\Fix(\tau)/\langle\tau\rangle$ is the unknotted torus. A 3--bridge shadow diagram for $(S^4,T^2)$ is (almost) shown in Figure~\ref{fig:s2xs2_involutions}(\textsc{c}); the third shadow arc of each color is uniquely determined.
		\item $S^2\times S^2/\langle\tau\circ\rho\rangle\cong \overline\CP^2$ and $\Fix(\tau\circ\rho)/\langle\tau\circ\rho\rangle$ is the mirror conic $\overline\Cc_2\subset\overline\CP^2$. A doubly pointed trisection diagram for $(\overline\CP^2,\overline\Cc_2)$ is shown in Figure~\ref{fig:s2xs2_involutions}(\textsc{d}).
	\end{enumerate}
	See~\cite{GayMei_22_Doubly-pointed-trisection-diagrams} for more details on doubly pointed trisection diagrams.
	Note that the quotient maps (1) and (3) were encountered above in Example~\ref{ex:Q8-link}.
	All three of these quotient maps are branched coverings.
	Each quotient manifold inherits a quotient $\Z_2$--action.
	For (1) and (3), this is the hyperelliptic involution of the torus, and for (2), this is rotation of the sphere through $\pi$ radians; these actions respect the quotient diagrams in Figure~\ref{fig:s2xs2_involutions}.
\end{example}

\subsection{Inferring group actions on $T^2\times S^2$ and $T^4$ via trisection diagrams}
\label{subsec:marla_equivariant_trisections}

\begin{figure}[ht!]
	\centering
	\includegraphics[width = 0.8\linewidth]{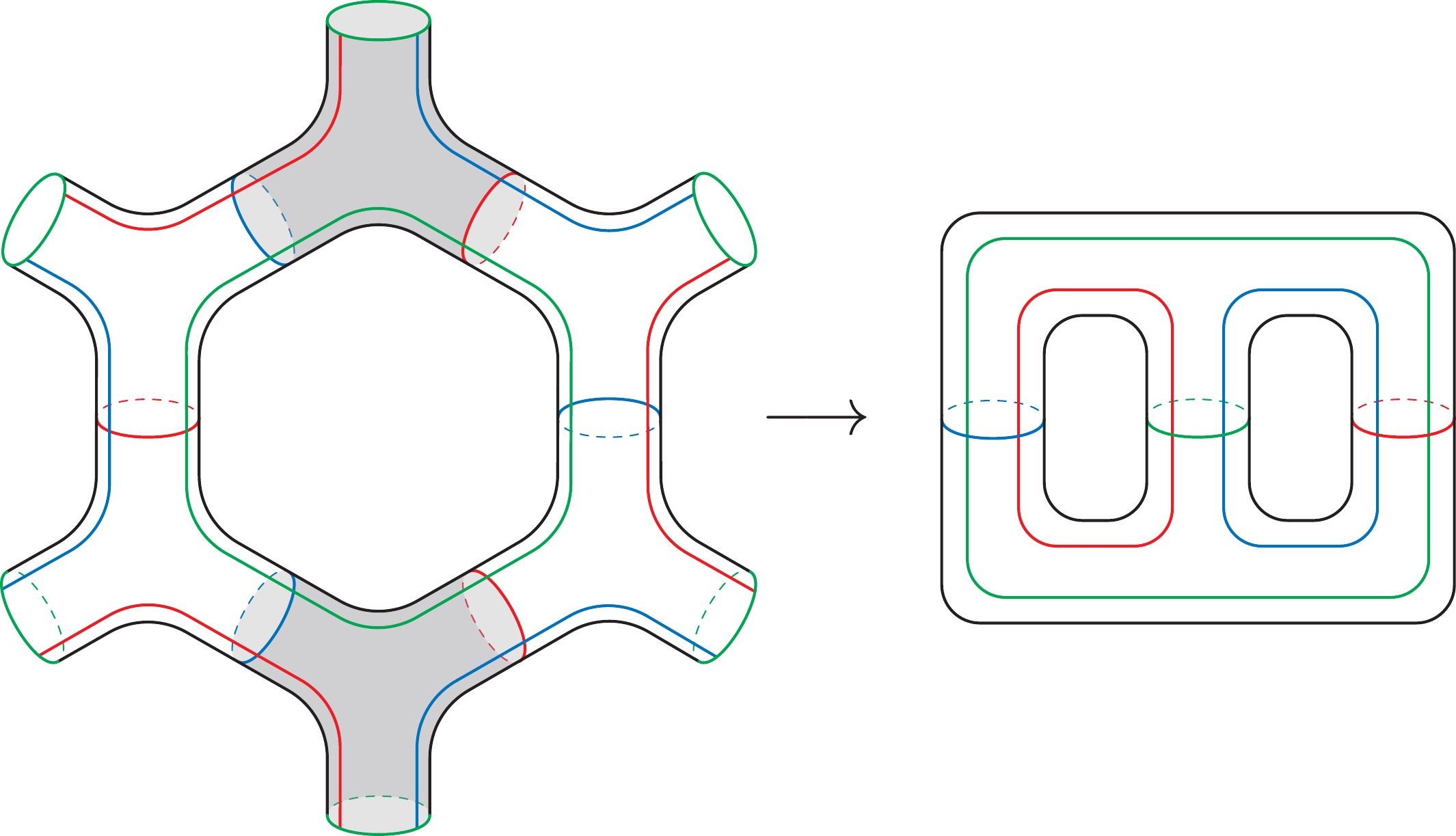}
	\caption{$S^2\times S^2$ as the quotient of a $\Z_3$--action on $T^2\times S^2$}
	\label{fig:t2xs2}
\end{figure}

This section is a proof of concept for the idea that one can infer group actions on manifolds through symmetries of their trisection diagrams and shows how to use our diagrammatic tools to deduce properties of these actions.
The trisections we consider were constructed by Williams in her thesis~\cite{Wil_20_Trisections-of-flat-surface-bundles}.
Consider the trisection diagram $\DD$ of $T^2\times S^2$ given in~\cite[Figure~3.12]{Wil_20_Trisections-of-flat-surface-bundles} and recreated on the left in Figure~\ref{fig:t2xs2}, where opposite boundary circles are identified by reflection.

This trisection diagram is invariant with respect to the $\Z_3$--action by generated by rotation about the center of the figure.
By Proposition~\ref{prop:diag_spine}, this action on the diagram determines a unique $\Z_3$--equivariant trisection $\TT$ on $T^2\times S^2$.

To understand the quotient $(T^2\times S^2)/\Z_3$,
Corollary~\ref{cor:quotient_diagram} indicates that we need only consider the quotient $\DD/\Z_3$ of the diagram $\DD$.
To compute this quotient, consider the fundamental domain $\Ii$ shaded in Figure~\ref{fig:t2xs2}.
The quotient $\Ii/\Z_3$ is obtained by identifying the boundary circles of $\Ii$ according to color.
The result is the standard genus-two trisection diagram for $S^2\times S^2$.

By Corollary~\ref{cor:quotient_diagram}, we have that $(T^2\times S^2)/\Z_3 = S^2\times S^2$.
Since the fixed-point set of the action on the spine is not in bridge position, we cannot apply Corollary~\ref{cor:bridge_in_spine}.
We expect, however, that some ad hoc analysis could be applied to conclude that this $\Z_3$--action on $T^2\times S^2$ is a branched covering action with branching set $\{a, b, c\}\times S^2\subseteq S^2\times S^2$ for some triple $\{a,b,c\}\in S^2$.

\begin{figure}[ht!]
	\centering
	\includegraphics[width = 0.9\linewidth]{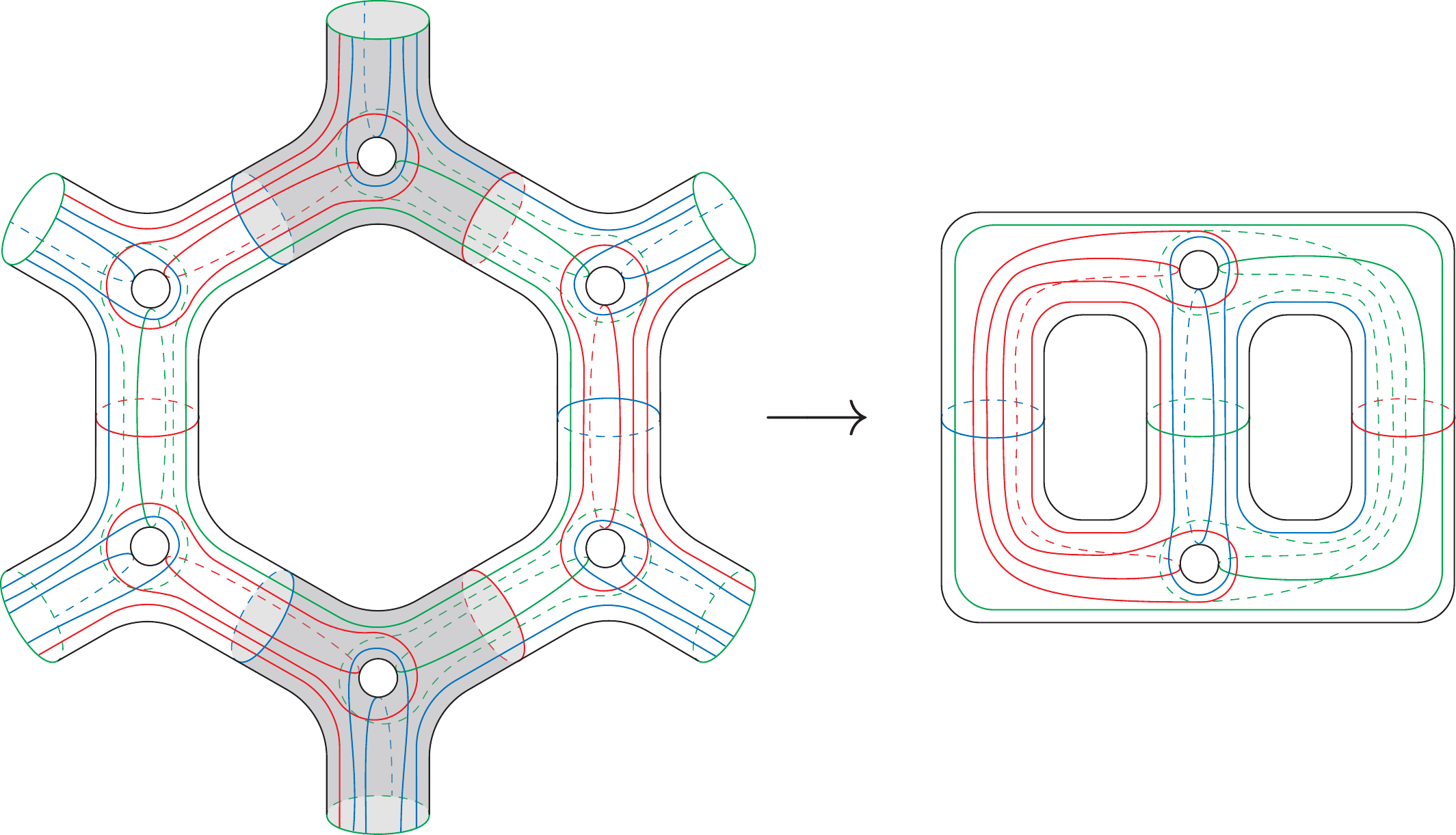}
	\caption{$T^2\times S^2$ as the quotient of a $\Z_3$--action on $T^4$}
	\label{fig:t2xt2}
\end{figure}

Continuing, consider the trisection diagram $\DD$ shown to the left in Figure~\ref{fig:t2xt2}, which is adapted from~\cite[Figure~2.3]{Wil_20_Trisections-of-flat-surface-bundles} and describes $T^4 = T^2\times T^2$.
Again, this trisection diagram is invariant with respect to the $\Z_3$--action by generated by rotation about the center of the figure.
By Proposition~\ref{prop:diag_spine}, this action on the diagram determines a unique $\Z_3$--equivariant trisection $\TT$ on $T^2\times T^2$.
To compute this quotient, consider the fundamental domain $\Ii$ shaded in Figure~\ref{fig:t2xt2}.
The quotient $\Ii/\Z_3$ is obtained by identifying the boundary circles of $\Ii$ according to color.
The result is the genus-four trisection diagram for $T^2\times S^2$ shown to the right in Figure~\ref{fig:t2xt2}.
Williams shows the diagram on the right in Figure~\ref{fig:t2xt2} and the diagram on the left in Figure~\ref{fig:t2xs2} are diagrams for the same trisection of $T^2\times S^2$.

Again, we expect that one could verify that this $\Z_3$--action corresponds with a branched covering over three fibers, this time $T^2\times\{a,b,c\}\subseteq T^2\times S^2$; however, the same caveats as before apply here.

\subsection{$PU(3)$--equivariant trisections of $\CP^2$}
\label{subsec:PU(3)_trisections}

In this subsection, we describe equivariant trisections for a large class of linear actions on $\CP^2$.

The space $\C^3$ comes equipped with a natural action of $U(3)$.
This descends to an action of $PU(3)$ on $\CP^2$ as follows.
For a point $x\in \C^3$, let $[x]$ be the quotient point in $\CP^2$.
Similarly, for a matrix $A\in U(3)$, let $[A]\in PU(3)$.
Then, we define $[A]\cdot [x] = [Ax]$.
Any subgroup of $PU(3)$ then has a natural action on $\CP^2$ by matrix multiplication.

There has been some study of finite group actions on $\CP^2$, with one major conclusion being that any finite group acting on $\CP^2$ must be isomorphic to a finite subgroup of $PU(3)$~\cite{HamLee_88_Finite-group-actions,HamLeeMad_89_Rigidity-of-certain-finite,Wil_87_Group-actions-on-the-complex}.
However, much is still unknown about which actions on $\CP^2$ are conjugate to an action by a subgroup of $PU(3)$.
In particular, it is open whether every (homologically trivial) cyclic action on $\CP^2$ is linear.
In Subsection~\ref{subsec:genus_1_classification}, we show that every action that respects the genus-one trisection of $\CP^2$ is linear.

Presently, we will restrict our attention to finite subgroups of the subgroup $U(2)<PU(3)$, which we identify with $\{[1\oplus A]\,:\,A\in U(2)\}$.
By work of Blichfeld~\cite{MilBliDic_61_Theory-and-applications-of-finite-groups} nearly all finite subgroups of $PU(3)$ have this form, but a few exceptional subgroups do not; see~\cite[Section~1]{HamLee_88_Finite-group-actions}.
An analysis of equivariant trisections of these exceptional subgroups would require new techniques beyond those that we now describe, which suffice for the subgroups of $U(2)$.

We begin by establishing a basic understanding of how $U(2)<PU(3)$ acts on $\CP^2$.
First, $U(2)$ fixes the point $[1:0:0]$ and leaves the line $\CP^1 = \{[0:z_2:z_3]\}$ invariant.
Second, the center $Z(U(2))\cong U(1)$ consists of elements of the form $[1\oplus\Diag(\lambda)]$ for some $\lambda\in U(1)$.
Finally, for any $[1\oplus A]\in U(2)$, the induced action on $\CP^1$ is by $[A]\in PU(2)\cong U(2)/U(1)$.

As a first step towards building an equivariant trisection with respect to the action of some subgroup of $U(2)$, we will give a $U(2)$--invariant decomposition of $\CP^2$.
Consider the 4--ball
$$X_1 = \{[1:z_2:z_3]\,|\, |z_2|^2 + |z_3|^2 \leq 1\},$$
which is clearly $U(2)$--invariant,
and let $N = \overline{\CP^2\setminus X_1}$.
Then, $\CP^2 = X_1\cup_{Y}N$, where $Y = \partial X_1 = \partial N\cong S^3$ is given by
$$Y =  \{[1:z_2:z_3]\,|\, |z_2|^2 + |z_3|^2 = 1\}.$$
The projection $\pi\colon N\to \CP^1$ defined by projecting onto the last two homogeneous co-ordinates gives $N$ the structure of a disk-bundle over $\CP^1$, which can be parametrized as
$$N = \{[z_1,z_2:z_3]\,|\, |z_1|\leq 1, |z_2|^2 + |z_3|^2 = 1\}.$$
Summarizing, we have the following $U(2)$--equivariance of the disk-bundle structure on $N$.

\begin{proposition}
\label{prop:U2_bundle}
	Let $q\colon U(2)\to PU(2)$ be the usual quotient map by $Z(U(2))$.
	Let $\pi\colon N\to\CP^1$ denote projection onto the last two coordinates.
	For $[1\oplus A]\in U(2)$, $x\in N$, we have \[\pi([1\oplus A]x) = q([1\oplus A])\pi(x) = [A]\pi(x)\text{.}\]
	In particular, the $U(2)$--action takes fibers to fibers.
	Furthermore, if $y\in \CP^1$ is a fixed point of $q([1\oplus A])$, then $y$ is an eigenvector of $A$ with some eigenvalue $\lambda$, which is unit magnitude since $A$ is unitary, and the action of $[1\oplus A]$ on the disk $\pi^{-1}(y)$ is by complex multiplication by $\lambda^{-1}$.
\end{proposition}

\begin{proof}
	Writing $[z_1:\vec{v}]$ to mean $[z_1:v_1:v_2]$, let $[z_1:\vec{v}]\in N$, so $|v_1|^2+|v_2|^2 = 1$, $|z_1|\leq 1$.
	Applying the definition of the $PU(3)$ action, we compute \[\pi\left([1\oplus A][z_1:\vec{v}]\right) = \pi\left([z_1:A\vec{v}]\right) = [A\vec{v}]\] where $[A\vec{v}]$ means $[(A\vec{v})_1:(A\vec{v})_2]$.
	By the definition of the $PU(2)$--action on $\CP^1$, we have $[A\vec{v}] = [A][\vec{v}]$ where $A\in U(2)$ and $[A]\in PU(2)$ is its associated projective equivalence class.
	Thus $[A\vec{v}] = [A][\vec{v}] = q([1\oplus A])\pi([z_1:\vec{v}])$ and the equality holds as desired.

	For the second part, let $[\vec{y}]\in \CP^1$ be fixed by $q([1\oplus A]) = [A]$, so we have $[A][\vec{y}] = [A\vec{y}] = [\vec{y}]$ as classes in $\CP^1$, and thus $A\vec{y} = \lambda\vec{y}$ for some $\lambda\in \C$.
	Thus $\vec{y}$ is an eigenvector and $\lambda$ is its eigenvalue, which is unit magnitude since $A$ is unitary.
	The projection $\pi$ projects down from $N$, so the disk $\pi^{-1}([\vec{y}])$ can be parameterized as $\{[z_1:\vec{y}]\,|\, |z_1|\leq 1\}$.
	Thus if $A\vec{y} = \lambda \vec{y}$, we have \[[1\oplus A][z_1:\vec{y}] = [z_1:A\vec{y}] = [z_1:\lambda\vec{y}] = [\lambda^{-1}z_1:\vec{y}]\] as claimed.
\end{proof}

Now, let $\widetilde{G}$ be a finite subgroup of $U(2)$ acting on $\CP^2$.
Let $G = q(\widetilde{G})$, and let $C$ be the kernel of $q|_{\widetilde{G}}$. 
Note that $C = \widetilde G\cap Z(U(2))$, thus $C$ is cyclic.
By Proposition~\ref{prop:U2_bundle}, for any $\widetilde g\in\widetilde G$ and $q(\widetilde g) = g\in G$, we have $\pi(\widetilde g\cdot x) = g\cdot\pi(x)$ for all $x\in N$.

In abuse of notation, we let $\pi\colon Y\to \CP^1$ denote the restriction of the projection introduced above to the 3--sphere $Y$.
This projection gives the Hopf fibration on $Y$: for $x\in\CP^1$, $\pi^{-1}(x)\cong S^1$.
The action of $Z(U(2)) = U(1)$ on $Y$ is by rotation along Hopf fibers, and by Proposition~\ref{prop:U2_bundle}, we see that $\pi$ is the quotient map of this $U(1)$--action.

The next step in trisecting $\CP^2$ equivariantly is finding an appropriate equivariant Heegaard splitting for $Y$.

\begin{proposition}
\label{prop:PU(2)_Heegaard}
	There is a $\widetilde G$--invariant Heegaard splitting $Y = H_1\cup_\Sigma H_2$ with the property that $H_i$ is a $\widetilde G$--invariant neighborhood of an invariant graph $\widetilde\Gamma_i$ such that $\pi(\widetilde\Gamma_i)$ is a $G$--invariant graph $\Gamma_i\subset\CP^1$ that gives $\CP^1$ the structure of a polyhedron.
	The graph $\widetilde\Gamma_i$ contains $|V(\Gamma_i)|$ disjoint Hopf fibers as cycles, and $g(\Sigma) = |O_{E(\Gamma_1)}||\widetilde G|+1$, where $|O_{E(\Gamma_1)}|$ denotes the number of $G$--orbits of the edges of $\Gamma_1$.
\end{proposition}

\begin{proof}
	The action of $G<PU(2)$ on $\CP^1$ can be regarded as a linear action by a subgroup of $SO(3)\cong PU(2)$.
	We can thus imagine $\CP^1$ as a polyhedron with $G$ acting on the polyhedral decomposition vertex--transitively.
	Let $\Gamma_1$ be the 1--skeleton of this polyhedral decomposition of $\CP^1$.
	Figure~\ref{fig:graph_lifts} shows $\CP^1$ with the octahedral 1--skeleton shown in cyan.
	If $\Gamma_1$ has an edge that is inverted by an element of $G$, equivariantly replace each inverted edge with a pair of parallel edges bounding a bigon.
	After this, $\Gamma_1$ has the property that no edge is left invariant by an element of $G$.

	Let $\Gamma_2$ denote the dual graph, which we choose to be $G$--invariant.
	Figure~\ref{fig:graph_lifts} shows the dual graph to the octahedral decomposition in orange.
	Let $M = \Gamma_1\cap\Gamma_2$, which we refer to as the \emph{midpoints} of the edges.
	Let $\nu(\Gamma_1)$ and $\nu(\Gamma_2)$ denote $G$--invariant tubular neighborhoods such that $\CP^1 = \nu(\Gamma_1)\cup\nu(\Gamma_2)$ and $\nu(\Gamma_1)\cap\nu(\Gamma_2) = \partial(\nu(\Gamma_1)) \cup \nu(M) = \partial(\nu(\Gamma_2))\cup\nu(M)$.
	Note that $\nu(M)$ is a disjoint collection of disks in $\CP^1$; one such disk is shown in Figure~\ref{fig:graph_lifts}.

	Now, let $H_1' = \pi^{-1}(\nu(\Gamma_1)\setminus\nu(M))$, let $H_2'=\pi^{-1}(\nu(\Gamma_2)\setminus\nu(M))$, and let $W = \pi^{-1}(\nu(M))$.
	Each of $H_1'$, $H_2'$, and $W$ is a disjoint collection of solid tori that are fibered by Hopf fibers.
	Moreover, $Y = H_1'\cup W \cup H_2'$.

\begin{figure}[ht!]
     \centering
     \includegraphics[width=.8\textwidth]{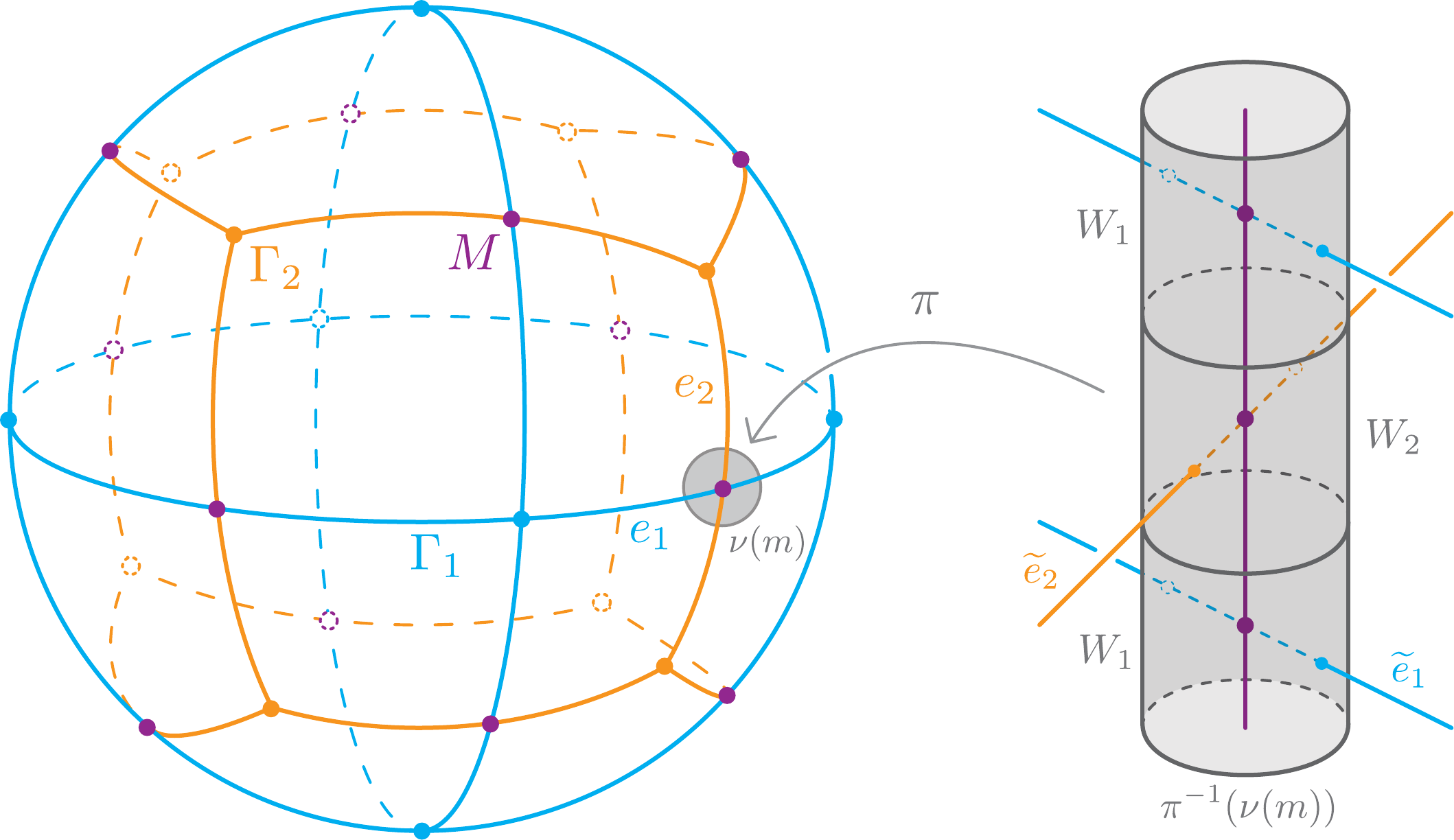}
     \caption{The octahedral graph $\Gamma_1$ on $\CP^1$ (cyan), together with its dual graph $\Gamma_2$ (orange), and a portion of the solid torus fiber $\pi^{-1}(\nu(m))$ for a point $m\in\Gamma_1\cap\Gamma_2$}
     \label{fig:graph_lifts}
\end{figure}

	Define a lift $\widetilde\Gamma_1$ of $\Gamma_1$ to $Y$ as follows.
	First, take the total pre-image $\pi^{-1}(V(\Gamma_1))$ of the vertices $V(\Gamma_1)$; this is a disjoint union of Hopf fibers, one for each vertex.
	Choose a point $p$ on one such Hopf fiber, and define $V(\widetilde\Gamma_1))$ to be the orbit $\widetilde G\cdot p$.
	This is a collection of points spread across the Hopf fibers $\pi^{-1}(V(\Gamma_1))$, with at least one point on each fiber since $G$ is vertex--transitive on $\Gamma_1$.

	Consider a point $q\in V(\widetilde\Gamma_1)$, let $\pi(q) = v$, and let $e\in E(\Gamma_1)$ be an edge of $\Gamma_1$ incident to $v$.
	Consider the annulus $\pi^{-1}(e)$.
	Since the $G$--stabilizer of $e$ is trivial, the induced action of $\widetilde G$ on $\pi^{-1}(e)$ is that of $C = \widetilde{G}\cap Z(U(2))$: cyclic rotation in the direction of the Hopf fibers.
	There is at least one $\Stab_{\widetilde G}(\pi^{-1}(e))$--orbit of vertices of $\widetilde\Gamma_1$ along each boundary component of $\pi^{-1}(e)$.
	Choose a $\Stab_{\widetilde G}(\pi^{-1}(e))$--equivariant arc $\widetilde{e}$ from $q$ to a vertex of $\widetilde\Gamma_1$ on the other circular boundary component of $\pi^{-1}(e)$ such that $\pi\colon \widetilde e\to e$ is a diffeomorphism.
	
	The arc $\widetilde{e}$ is an equivariant lift of $e$ to $Y$.
	Repeat this process for one edge in each $G$--orbit of edges in $\Gamma_1$ to obtain a collection of equivariant lifts.
	Taking the union of the orbits of all of these arcs gives a $\widetilde{G}$--invariant collection of edges which projects down onto the edges of $\Gamma_1$; we refer to these edge lifts as \emph{horizontal}.
	Define $E(\widetilde\Gamma_1)$ to be the union of these horizontal edges with with the (\emph{vertical}) segments of the Hopf fibers $\pi^{-1}(V(\Gamma_1))\setminus V(\widetilde\Gamma_1)$.
	Note that the original Hopf fibers $\pi^{-1}(V(\Gamma_1))$ are contained as disjoint cycles in $\widetilde\Gamma_1$ and that $\pi(\widetilde\Gamma_1) = \Gamma_1$, since there is at least one lifted edge in $\widetilde\Gamma_1$ over each edge of $\Gamma_1$.

	We now calculate $\chi(\widetilde\Gamma_1)$, noting that we will have $g(\Sigma) = g(H_1) = 1 - \chi(\widetilde\Gamma_1)$ once we have defined $H_1$.
	Let $|O_{E(\Gamma_1)}|$ denote the number of $G$--orbits of the edges of $\Gamma_1$.
	We can think of $\widetilde\Gamma_1$ as the union of the $|V(\Gamma_1)|$ Hopf fibers with the horizontal edges constructed above.
	By the hypothesis that no edge of $\Gamma_1$ is left invariant by an element of $G$, we know that no edge of $\widetilde\Gamma_1$ is left invariant by an element of $\widetilde G$.
	It follows that the number of horizontal edges is precisely $|\widetilde G||O_{E(\Gamma_1)}|$, yielding $\chi(\widetilde\Gamma_1) = -|\widetilde G||O_{E(\Gamma_1)}|$, as desired.

	We now define $\widetilde\Gamma_2$ in the same way, using $\Gamma_2$ in place of $\Gamma_1$.
	Since $G$ doesn't necessarily act vertex-transitively on $\Gamma_2$, choose one vertex of $\Gamma_2$ for each orbit of the action of $G$ on $V(\Gamma_2)$, and call this collection $T$.
	Choose one point $p$ on each Hopf fiber of $\pi^{-1}(T)$, and let $V(\widetilde G)$ be be the union of orbits $\widetilde G\cdot p$ across each such point $p$. 
	Now, for some choice(s) of edge lifts in $\widetilde\Gamma_2$, the lifts of edges may intersect similar edges of $\widetilde\Gamma_1$ in exactly one point (which is a point of $\pi^{-1}(M)$).
	Suppose that we have chosen lifts for the edges of $\Gamma_2$ so that each lift $\widetilde e$ of an edge $e$ of $\Gamma_2$ intersects at most one lift $\widetilde f$ of an edge $f$ of $\Gamma_1$, with the intersection occurring if and only if $e$ and $f$ intersect in $m$.
	In other words, the lifts of $e$ and $f$ are in bijection and have Kronecker delta intersection along points of $\pi^{-1}(m)$.
	Now, slightly perturb the edges of $\widetilde\Gamma_2$ vertically near the fiber $\pi^{-1}(m)$ so that $\widetilde\Gamma_1\cap\widetilde\Gamma_2 = \varnothing$, and
	 the edges of $\widetilde\Gamma_1$ and $\widetilde\Gamma_2$ intersect the Hopf fiber $\pi^{-1}(m)$ in an alternating fashion; see Figure~\ref{fig:graph_lifts}.

	For each $m\in M$, consider the solid torus $\pi^{-1}(\nu(m))$.
	Edges of $\widetilde\Gamma_1$ and $\widetilde\Gamma_2$ pass through this solid torus as horizontal arcs (with respect to the Hopf fibration) according to the alternating pattern described above.
	It follows that we can decompose $\pi^{-1}(\nu(m))$ equivariantly into some number of 3--balls $D^2\times I$, where each $D^2\times\{t\}$ is horizontal with respect to the fibration, such that each 3--ball is an equivariant neighborhood of the midpoint of an edge of $\widetilde\Gamma_1$ or $\widetilde\Gamma_2$.
	Let $W_1$ be the collection of $3$--balls which are neighborhoods of edges in $\widetilde \Gamma_1$, and let $W_2$ be the $3$--balls which are neighborhoods of edges in $\widetilde \Gamma_2$.
	Let $H_i = H_i'\cup W_i$, and let $\Sigma = \partial H_i$.
	Each 3--ball of $W_i$ can be seen as a 3--dimensional 1--handle attached to $H_i$, so it follows from the fact that $H_i'$ is a handlebody that $H_i$ is a handlebody; see Figure~\ref{fig:graph_lifts}.
	Since $Y_1 = H_1\cup_\Sigma H_2$, we have the desired equivariant Heegaard splitting.
%
\end{proof}

So far, we have the equivariant decomposition $\CP^2 = X_1\cup_{Y}N$, where $X_1$ is a 4--ball, and $Y = H_1\cup_\Sigma H_2$ is an equivariant Heegaard decomposition.
We now divide $N$ into the sectors $X_2$ and $X_3$ that will complete the equivariant trisection.

By Proposition~\ref{prop:PU(2)_Heegaard}, $H_1$ is a neighborhood of the graph $\widetilde\Gamma_1$, which contains $|V(\Gamma_1)|$ Hopf fibers.
For each such Hopf fiber $f$, let $D_f$ be the associated disk-fiber in $N$ with $\partial D_f = f$.
Let $X_3$ be a 4--dimensional equivariant thickening of the invariant set $H_1\bigcup_fD_f$ into $N$ such that $X_3\cap X_1= H_1$.
Let $X_2 = N\setminus X_3$.

\begin{proposition}
\label{prop:PU(3)_trisection}
	The decomposition $\CP^2 = X_1\cup X_2\cup X_3$ is a $\widetilde G$--equivariant $(g;0,k_2,k_3)$--trisection, where
	\begin{enumerate}
		\item $g = |\widetilde G||O_{E(\Gamma_1)}| + 1$
		\item $k_2 = |V(\Gamma_1)|-1$
		\item $k_3 = g - |V(\Gamma_1)|$
	\end{enumerate}
\end{proposition}

\begin{proof}
	First, we observe that $X_2$ and $X_3$ are 4--dimensional 1--handlebodies.
	The space $X_2$ is obtained from $N$ by first removing a finite number of thickened disk-fibers, then removing a portion of a collar of the boundary.
	(This second removal does not affect the diffeomorphism type of $X_2$.)
	It follows that $X_2$ is a disk-bundle over a punctured sphere, so $X_2\cong \natural^{k_2}(S^1\times B^3)$, where
	\[k_2 = 1 - \chi(\CP^1\setminus\nu(V(\Gamma_1))) = |V(\Gamma_1)| - 1.\]
	Furthermore, $X_3$ is obtained from a thickening of $H_1$ by attaching $|V(\Gamma_1)|$ 4--dimensional 2--handles along core curves of $H_1$.
	It follows that $X_3\cong\natural^{k_3}(S^1\times B^3)$, where $k_3 = g - |V(\Gamma_1)|$.
	Moreover, we have $\partial X_3 = H_1\cup_\Sigma H_3$, where $H_3$ is the result of Dehn surgery on the cores in $H_1$.
	By construction, $X_2\cap X_3 = H_3$, so the decomposition is a trisection, as desired.
	In addition, by construction this trisection is $\widetilde G$--invariant.
	
	The calculation of $g = g(\Sigma)$ follows from Proposition~\ref{prop:PU(2)_Heegaard}.
	
	It remains to show that the pieces of the trisection are linearly parted with respect to the $\widetilde G$--action.
	The $\widetilde G$--action on $X_1$ is simply the cone of the $\widetilde G$--action on the 3--sphere $Y$, so this sector is as desired.
	
	To find a linearly parting ball-system for $X_3$, first choose an equivariant cut-disk $D_m$ for $H_1$ dual to each midpoint $m$ of each horizontal edge of $\widetilde\Gamma_1$.
	These disks cut $H_1$ into the disjoint union of solid tori neighborhoods of the distinguished Hopf fibers.
	Let $\Bb_3 = \cup_m(D_m\times I)$ be the thickening of $D_m$ into $X_3$.
	Then $X_3\setminus\Bb_1$ is a disjoint union of 4--balls, each of which is one of the 4--dimensional 2--handles used in the construction of $X_3$.
	The stabilizer of each such 4--ball is the same as the stabilizer of a disk-fiber in $N$, which is linear.
	
	For $X_2$, let $\pi_2\colon X_2\to \CP^1$ be the restriction of $\pi\colon N\to \CP^1$ to $X_2\subset N$.
	Let $\Bb_2 = \pi_2^{-1}(E(\Gamma_1))$.
	Since $X_2 = \pi_2^{-1}(\CP^1\setminus\nu(V(\Gamma_1))$, we have that $\Bb_2$ cuts $X_2$ into 4--balls, each of which is the pre-image of a face of $\Gamma_1$.
	Each of these is a neighborhood of a disk-fiber in $N$, so the action has linear stabilizer.
\end{proof}

\begin{remark}
\label{rmk:SU2}
	There is a natural subgroup $SU(2)<U(2)$ consisting of the matrices with unit determinant.
	The action of this subgroup on $Y$ is the standard, free $SU(2)$--action on $S^3$.
	Thus if $\widetilde G<SU(2)$, it is clear that $\widetilde G$ acts pseudo-freely on $\mathbb CP^2$. In particular, we see that $\widetilde G$ acts freely on the handlebodies $H_1$ and $H_2$ in the equivariant Heegaard splitting in Proposition~\ref{prop:PU(2)_Heegaard}.

	For a free finite group action of $G$ on a handlebody $V$, the genus of $V$ is $1+|G|(\mu-1)$, where $\mu$ is the genus of $V/G$ \cite[Section~1]{McCMilZim_89_Group-actions-on-handlebodies}. 
	This ``discretizes'' the possible values for the genus of $V$; in particular, if $V$ has genus at least $2$, we see that the genus of $V$ is actually at least $1+|G|$.
	Looking at the computation of the genus of the Heegaard splitting in Proposition~\ref{prop:PU(3)_trisection} above, we see that if we choose a graph $\Gamma_1$ upon which $G$ acts \emph{both} vertex-transitively and edge-transitively, then the Heegaard genus of the $\widetilde G$--equivariant Heegaard splitting $Y = H_1\cup_\Sigma H_2$ is minimal.
	In fact, since $\widetilde G$ acts pseudo-freely on $\mathbb CP^2$, in any equivariant trisection of the action, $\widetilde {G}$ must act freely on each handlebody of the spine, and hence if we choose $\Gamma_1$ as above, the trisection genus of the resulting trisection is minimal.
	This computes the equivariant trisection genus for each finite group action on $\CP^2$ arising from a subgroup $\widetilde G<SU(2)$; specifically, the (equivariant) trisection genus is $1+|\widetilde{G}|$.
	
	Generalizing the above, we have that for any $G$ acting on $X$ freely or pseudo-freely, the equivariant trisection genus must be $0$ or $1+n|G|$, for $n\in \N\cup\{0\}$.
\end{remark}

\begin{example}
	We now examine a particular case of interest to illustrate the construction. The \emph{binary tetrahedral group} $2T$ is group double cover of the tetrahedral group $T = A_4$ and can be seen as a subgroup of the matrix group $SU(2)$.
	The binary tetrahedral group is familiar to topologists as the fundamental group of an exceptional $3$--dimensional spherical space form (a Seifert fibered space over $S^2(2,3,3)$), since the usual action of $SU(2)$ (hence $2T$) on $S^3$ as the unit sphere in $\C^2$ is free and thus the quotient is a closed $3$--manifold covered by $S^3$.
	We can see $SU(2)<U(2)$ acting on $\CP^2$ as in Remark~\ref{rmk:SU2}, and the action of $SU(2)$ on the boundary of $X_1$ is equivalent to usual action of $SU(2)$ on $S^3$ as above.

	The extention from $T$ to $2T$ corresponds to extending $PU(2) \equiv SO(3)$ to $SU(2)$.
	In the notation of the remarks following Proposition~\ref{prop:U2_bundle}, we have $\widetilde{G} = 2T$ and  $G = T$, with $C = \Z_2$, generated by the diagonal matrix $-I$.
	Viewing $2T$ acting on $\CP^2$ via the above linear action, when we restrict to $\CP^1$ we see the usual $T$--action on $\CP^1$, acting as the orientation-preserving symmetries of the tetrahedron.
	
	There is a particularly auspicious lift of the tetrahedral graph on $\CP^1$ to $S^3$, namely the $1$--skeleton of the $24$--cell, a $4$--dimensional regular polytope with $2T$ symmetry group acting freely on it.
	Interestingly (and as part of the unusual symmetry of the $24$--cell), we can obtain this same lift if we choose $\Gamma_1$ to be the $1$--skeleton of the octahedron, with $T$ acting as a subgroup of the full $S_4$ symmetry group of the octahedron.
	This is the same $T$--action on the sphere, which can be seen by inscribing the tetrahedron in the cube, the dual of the octahedron.
	Furthermore, the $24$--cell is self-dual, and also appears as a lift if we pick $\Gamma_1$ to be the $1$--skeleton of the dual cube to that octahedron.

\begin{figure}[ht!]
     \centering
     \includegraphics[width=.9\textwidth]{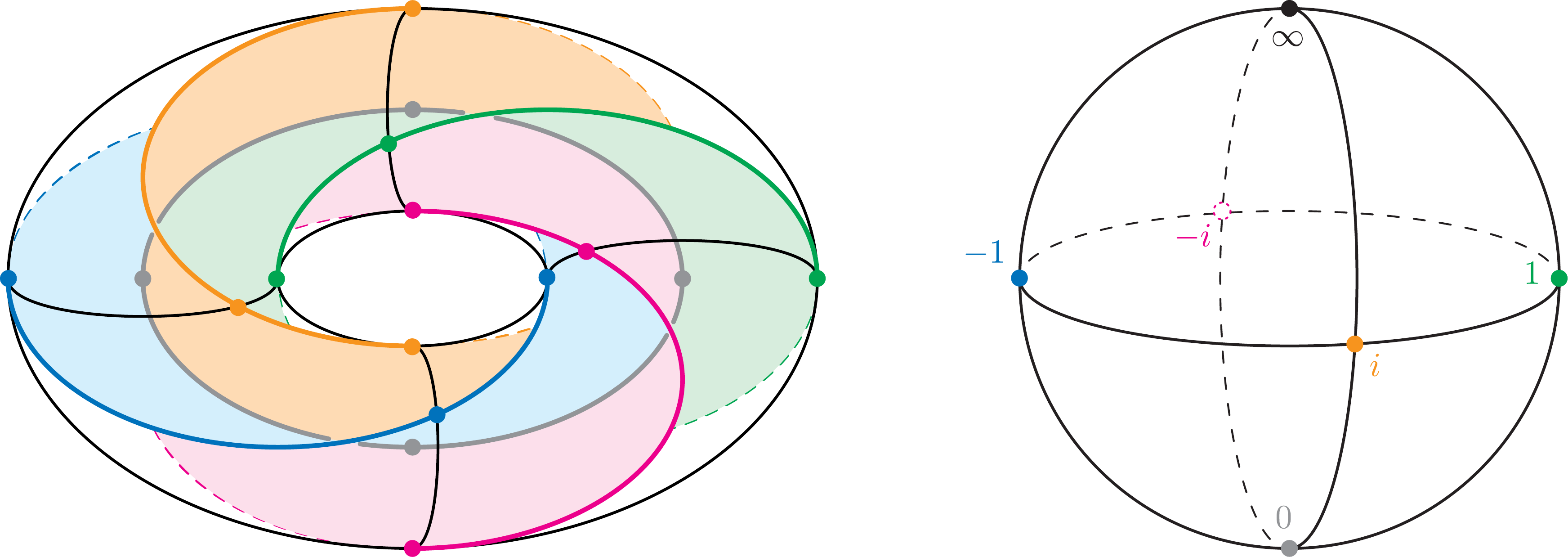}
     \caption{Lifting the southern half of the octahedral graph on $\CP^1$ to the Hopf fibration of $S^3$}
     \label{fig:Hopf_fibration}
\end{figure}
	
	Presently, we'll work with $\Gamma_1$ as the octahedral graph, denoting the six vertices by
	$$\{[1:0], [0:1], [1:1], [1:-1], [1:i], [1:-i]\}\subseteq \CP^1,$$
	since for this choice there is a nice parameterization of $\widetilde\Gamma_1$ that makes the computations simpler.
	See the right side of Figure~\ref{fig:Hopf_fibration}, where the vertices are labeled via the mapping $[a:b]\mapsto b/a$.
	We now describe explicitly the graph $\widetilde\Gamma_1$ constructed in Proposition~\ref{prop:PU(2)_Heegaard}.
	This graph will be closely related to the $24$--cell, but will differ in key ways; for example, no edge of the $24$--cell is vertical with respect to the Hopf fibration, while many edges of $\widetilde\Gamma_1$ have this property.

	First, to describe the $24$--cell and the $2T$ action on $\CP^2$, we keep in mind two perspectives at once: $2T$ as a subgroup of $SU(2)$ and as a subgroup of the unit quaternions.
	The unit quaternions are a subgroup of the quaternions, and have an action on $\C^2$ as follows: $a1+bi+cj+dk$ acts via the complex matrix
$$\begin{bmatrix}
a+bi & c+di \\
-c+di & a-bi 
\end{bmatrix}.$$
	Furthermore, each unit quaternion $a1+bi+cj+dk$ can be seen as a unit-length element of $\C^2$ by sending $a1+bi+cj+dk\mapsto (a+bi, c+di)$.
	These identifications are compatible in the sense that the action of the unit quaternions on themselves is the same as the action of these matrices on these points.
	The structure of the unit quaternions justifies the choice of the octahedral graph for $\Gamma_1$.
	The solid torus on the left side of Figure~\ref{fig:Hopf_fibration} contains the five Hopf fibers that project to the vertices of $\Gamma_1$, excluding $[1:0] = \infty$.
	The Hopf fiber over $\infty$ would be the core of the complementary solid torus in $S^3$ and is not pictured.
	The vertical annuli sitting over the four edges incident to $0$ in $\Gamma_1$ and also shown in the solid torus.

	The binary tetrahedral group $2T$ is the subgroup of the unit quaternions generated by the elements $A = \frac{1}{2}(1+i+j+k)$ and $B = \frac{1}{2}(1+i+j-k)$.
	Equivalently, it is the group of unit quaternions consisting of the 24 points
	$$\left\{\pm1, \pm i, \pm j, \pm k, \frac{\pm1\pm i\pm j\pm k}{2}\right\},$$
	which, as a set of points in $\C^2$, is
	$$\left\{(\pm1, 0), (\pm i, 0), (0, \pm 1), (0, \pm i), \left(\frac{\pm1\pm i}{2},\frac{\pm1\pm i}{2}\right)\right\}.$$
	This latter set of points are the vertices of the $24$--cell, which we include into $\CP^2$ in the usual way, sending $(a,b)$ to $[1:a:b]$.
	Projecting this set of points onto $\CP^1$, we get the vertices above for the octahedron, with four vertices of the $24$--cell projecting down to each of the $6$ vertices of the octahedron.
	For example, the preimage of $[1:0]$ is a single Hopf fiber that contains the four points $[1:\pm1:0]$ and $[1:\pm i:0]$; this is indicated by the gray dot in $\CP^1$ and the four gray dots of the gray core of the solid torus to the left in Figure~\ref{fig:Hopf_fibration}.
	We leave it as an exercise to the reader to label the points on the Hopf fibers in the solid torus with the corresponding unit quaternions.
	
	We now construct $\widetilde\Gamma_1$, beginning by taking the vertex set to be the above vertices of the $24$--cell.
	For each of the six vertices of $\Gamma_1$, we have a Hopf fiber containing four vertices of $\widetilde\Gamma_1$, giving 24 vertices total.
	This also gives 24 vertical edges: the segments of the six Hopf fibers.
	The action of $T$ on $\Gamma_1$ is edge transitive, so $|O_{E(\Gamma_1)}| = 1$.
	It follows that by choosing a single edge $e$ of $\Gamma_1$, lifting it to a horizontal edge $\widetilde e$, and taking the orbit $2T\cdot\widetilde e$, we get 24 horizontal edges in $\widetilde\Gamma_1$.
	It follows that $\widetilde\Gamma_1$ is a 4--regular graph with 24 vertices and 48 edges.
	The annular preimage $\pi^{-1}(e)$ contains two horizontal edges for each $e\in E(\Gamma_1)$
	
	Carrying out the procedure outlined in Proposition~\ref{prop:PU(2)_Heegaard}, the equivariant Heegaard splitting of $S^3$ is given as a tubular neighborhood of $\widetilde\Gamma_1$ and its ``dual'' $\widetilde\Gamma_2$, which yields a $2T$--equivariant $(25; 0, 17, 7)$--trisection $\TT$ of $\CP^2$.
	By Remark~\ref{rmk:SU2}, this equivariant trisection is genus minimizing for this $2T$--action.
		
	We note that different choices in the above exposition (such as choosing $\Gamma_1$ to be the cube or tetrahedron) yield possibly different, but still genus-minimal, trisections.
	We remark that the 24 horizontal edges of $\widetilde\Gamma_1$ represent some of the 96 edges of the 24--cell; to get all 96 edges, one must include eight horizontal edges in in each annulus $\pi^{-1}(e)$, rather than two.
	Moreover, the 24--cell contains no vertical edges, so the present construction will never recover the 24--cell exactly.
\end{example}

\section{Classification in low genera}
\label{sec:classification}

In this section we classify equivariant trisections with low genus.
We recall Definition~\ref{def:trisection_equivalence}: two equivariant trisections are \emph{$G$--diffeomorphic} if they are related by an equivariant diffeomorphism up to permuting the sector labels.

\subsection{Genus-zero equivariant trisections}
\label{subsec:genus_0_classification}

To characterize genus-zero equivariant trisections, we first prove a general proposition about suspending equivariant Heegaard splittings of $S^3$.
This proposition will be used when studying genus-one equivariant trisections of $S^4$ below, as well.

\begin{proposition}
\label{prop:suspension_heegaard}
	Let $G$ act linearly on $S^3$.
	Viewing $S^3$ as the equatorial sphere of $S^4$, let $G$ act on $S^4$ as the suspension of the linear action on $S^3$.
	Let $V_1\cup_{\Sigma} V_2$ be a $G$--equivariant, genus--$g$ Heegaard splitting of the equatorial $S^3$.
	Then, there exists a $G$--equivariant $(g;0,0,g)$--trisection of $S^4$ whose central surface is $\Sigma$.
\end{proposition}

\begin{proof}
	Let $h\colon S^4\to [-1, 1]$ be standard Morse function on $S^4$, so $h^{-1}(0)$ is the equatorial $S^3$.
	Adopting the Morse theory notation from Section~\ref{sec:existence}, define
	\begin{itemize}
		\item $X_3 = (V_1)_{[-\frac{1}{2}, \frac{1}{2}]}$,
		\item $X_1 = h^{-1}([0,1])\setminus \Int(X_3)$, and
		\item $X_2 = h^{-1}([-1, 0])\setminus \Int(X_3)$.
	\end{itemize}
	Verifying that $(X_1,X_2,X_3)$ gives a $(g;0,0,g)$--trisection of $S^4$ is routine, and the fact that this trisection is $G$--invariant follows from the fact that $h$ is $G$--invariant.

	To show that this trisection is equivariant, it suffices then to show that the $X_i$ are linearly parted.
	For $X_1$ and $X_2$, this follows since these are linear actions on a $4$--ball minus a neighborhood of a graph on the boundary, i.e. a linear action on a $4$--ball.
	For $X_3$, note that $V_1$ is linearly parted as a $3$--dimensional $1$--handlebody (see~\cite[Section~3]{MeiSco_LP}), and extending the parting system in the height direction gives rise to a linearly parting ball-system for $X_3$.
\end{proof}

As a corollary to the above proposition, we have a classification of equivariant genus-zero trisections that establishes that all such trisections are ``geometric''.

\begin{theorem}
\label{thm:class_g0}
	The set of genus-zero trisections, considered up to $G$--diffeomorphism, is in bijection with the set of actions on $S^2$ by subgroups $G<SO(3)$, considered up to smooth equivalence.
	Each $G$--diffeomorphism class of genus-zero equivariant trisections of $S^4$ is represented by an action of a subgroup of $SO(5)$.
\end{theorem}

\begin{proof}
	Given an action on $S^2$ by a subgroup $G<SO(3)$, we suspend to obtain a $G$--action on $S^3$ leaving the genus-zero Heegaard splitting invariant.
	By Proposition~\ref{prop:suspension_heegaard}, we can suspend again to obtain a genus-zero $G$--equivariant trisection of $S^4$.
	This shows that every smooth equivalence class of actions on $S^2$ gives rise to a $G$--diffeomorphism class of equivariant genus-zero trisections.
	Since we are viewing $S^4$ as the double suspension of $S^2$, we have a natural identification of $SO(3)$ as a subgroup of $SO(5)$.
	Any equivariant diffeomorphism of $S^2$ relating a given pair of such $G$--actions gives rise to a $G$--diffeomorphism of trisections by Proposition~\ref{prop:diag_spine}(2).
	
	The converse, that every equivariant genus-zero trisection is one of these, is immediate.
\end{proof}

\subsection{Genus-one equivariant trisections}
\label{subsec:genus_1_classification}

Moving on to genus-one equivariant trisections, we'll begin by characterizing the actions of finite groups on solid tori.

\begin{lemma}
\label{lem:genus_1_maximal_handlebody_group}
	Let $G$ be a finite group acting smoothly and orientation-preservingly on a solid torus $H\cong S^1\times D^2$.
	Then $G$ is a finite subgroup of $(\Z_m\times\Z_m)\rtimes \Z_2$ for some $m$, where $(\Z_m\times \Z_m)$ acts factor-wise and $\Z_2$ acts by reflection in each factor.
\end{lemma}

\begin{proof}
	Let $G$ act on a solid torus $H$.
	The~\ref{ELT} gives an equivariant handle decomposition for $H$ whose 0--handles and 1--handles are arranged according to a cyclic graph.
	In fact, with respect to the $G$--action, $H$ is an invariant tubular neighborhood of its core circle; i.e., $G$ acts on $H \cong S^1\times D^2$ preserving the product structure.
	Since $G$ is orientation-preserving, it must be a finite subgroup of the group of orientation-preserving diffeomorphisms of $H$ that preserve the product, which is $(S^1\times S^1)\rtimes \Z_2$, with $\Z_2$ acting as reflection in both factors and the group $S^1\times S^1$ acting on the concentric tori. 
	It now suffices to show that $G$ includes into $(\Z_m\times \Z_m)\rtimes\Z_2$. 
	The group $G$ has an index--2 subgroup $R$, consisting of those elements whose action on each factor of $S^1\times D^2$ is orientation-preserving.
	The group $R$ is a finite subgroup of $S^1\times S^1$, and thus is a finite subgroup of $\Z_t\times \Z_n$ for some $t$ and $n$, where each factor group is the image of $R$ under projection onto each circle-factor.
	In fact, by including $\Z_t, \Z_n\hookrightarrow \Z_m$ in the usual way, where $m$ is the least common multiple of $n$ and $t$, we can see that $R$ is a finite subgroup of $\Z_m\times\Z_m$. 
\end{proof}

We call a group of the form $(\Z_m\times \Z_m)\rtimes \Z_2$ a \emph{maximal finite handlebody group}, which we denote by $M_m$, and we call the above action of $M_m$ on $S^1\times D^2$ its \emph{natural action}.
Let $N_m = \Z_m\times \Z_m \vartriangleleft M_m$.
The action of $M_m$ on $\partial (S^1\times D^2) = S^1\times S^1$ will also be called the natural action of $M_m$ on $S^1\times S^1$.
We call a genus-one, $M_m$--equivariant trisection \emph{natural} if the action of $M_m$ on the central surface of the trisection is the natural one.

Recall that a genus-one trisection diagram $(\alpha_1,\alpha_2,\alpha_3)$ is standard if $|\alpha_i\cap\alpha_{i+1}|\leq 1$ for all $i\in\Z_3$.
A curve on $S^1\times S^1$ is a \emph{slope} if it lifts to a straight line in the universal cover $\R^2$.
An invariant trisection diagram is \emph{standard} if it is the orbit of a standard trisection diagram.

\begin{proposition}
\label{prop:natural}
	Let $(X,\TT,M_m)$ be a natural, genus-one, equivariant trisection.
	Both $\TT$ and $\TT/N_m$ have standard invariant trisection diagrams, and $\TT$ and $\TT/N_m$ are diffeomorphic (as non-equivariant trisections).
	The induced action of $\Z_2 = M_m/N_m$ on $\TT/N_m$ is by the hyperelliptic involution on the central surface $\Sigma/N_m$.
\end{proposition}

\begin{proof}
	Since the action of $M_m$ on $\Sigma\cong S^1\times S^1$ is natural, it acts factor-wise.
	Let $\alpha_i^0$ be a slope on $\Sigma$ bounding a disk in $H_i$.
	Since $\alpha_i^0$ is a slope, it is $M_m$--equivariant.
	By the~\ref{ELT}, $\alpha_i^0$ bounds an equivariant cut-disk.
	Let $\alpha_i$ be the orbit of $\alpha_i^0$.
	Then $(\alpha_1,\alpha_2,\alpha_3)$ is an invariant trisection diagram for $\TT$.
	
	Thinking of $\Sigma$ as the unit square in the usual way and subdividing it into an $m\times m$ grid, we see that the intersection of $(\alpha_1,\alpha_2,\alpha_3)$ with any fundamental domain (a small square in the grid) is just a copy of $(\alpha_1^0,\alpha_2^0,\alpha_3^0)$.
	Thus, $\TT/N_m$ has the same trisection diagram as $\TT$; since the $\alpha_i$ are slopes, this diagram is standard.
	
	The $\Z_2$--factor in $M_m = (\Z_m\times\Z_m)\rtimes\Z_2$ acts as rotation of the unit square through $\pi$ radians about its center.
	(The center, the corner(s), and the midpoints of the edges are the four fixed points.)
	The action of $\Z_2 = M_m/N_m$ is the same rotation, but applied to one of the smaller squares representing the quotient torus; this is the hyperelliptic involution.
\end{proof}

\begin{corollary}
\label{cor:natural}
	For each $m\in\N$, any two natural, genus-one, $M_m$--equivariant trisections of a given 4--manifold $X$ are $M_m$--diffeomorphic.
\end{corollary}

\begin{proof}
	Let $\TT$ and $\TT'$ be natural, genus-one, $M_m$--equivariant trisections,  for some $m\in\N$.
	We analyze the quotient of $\TT/N_m$ by the hyperelliptic involution $\tau$.
	This involution has two fixed arcs in each solid torus $H_i/N_m$, and the quotient $(\TT/N_m)/\langle \tau\rangle$ is a $2$--bridge trisection of a surface--link in $S^4$ by Theorem~\ref{thm:hyperelliptic}.
	
	So, the quotients $\TT/M_m$ and $\TT'/M_m$ are $2$--bridge trisections, which are classified: each is the quotient by $\langle \tau\rangle$ of one of the standard genus-one trisections~\cite[Section~4]{MeiZup_17_Bridge-trisections}.
	The diffeomorphism type of such a 2--bridge trisection is determined up to relabeling by the diffeomorphism type of $X/N_m$.
	Such a diffeomorphism and relabeling lifts to an equivariant diffeomorphism and relabeling of $\TT$ and $\TT'$, as desired.
\end{proof}

Having introduced and analyzed natural, genus-one trisections, we can state the classification:
In essence, every genus-one trisection extends to a natural, geometric one, which is unique up to the choice of $m$ by the previous corollary.

\begin{theorem}
\label{thm:class_genus_one}
	Let $(X,\TT,G)$ be an equivariant, genus-one trisection.
	For some $m\in\N$, there is a linear action of $M_m$ on $X$ extending the action of $G$ and leaving $\TT$ invariant so that $(X,\TT,M_m)$ is natural.
\end{theorem}

\begin{proof}
	Consider the induced action of $G$ on $H_1$.
	Since $G$ acts factor-wise on $H_1 = S^1\times D^2$, as in the proof of Proposition~\ref{prop:natural}, we can find a $G$--equivariant trisection diagram $(\alpha_1,\alpha_2,\alpha_3)$ for $\TT$ consisting of slopes in the coordinates of $H_1$.
	By Lemma~\ref{lem:genus_1_maximal_handlebody_group}, let $M_m$ be the maximal finite handlebody group for which $G<M_m$ and $M_m$ acts via its natural action on $H_1$ extending the action of $G$.
	Because the $\alpha_i$ are slopes and $M_m$ acts factor-wise, the diagram is also $M_m$--equivariant.
	The $M_m$--action on the diagram determines a $M_m$--equivariant trisection $\TT'$, by Proposition~\ref{prop:diag_spine}(1).
	Since the $M_m$--action extends the $G$--action on the diagram, by Proposition~\ref{prop:diag_spine}(2), up to $G$--diffeomorphism, $\TT'$ extends $\TT$.
	In other words, the $\TT$ is $M_m$--equivariant, with the $M_m$--action extending the original $G$--action.
	To complete the proof, we show that the natural $M_m$ actions are linear below.
\end{proof}
	
	We now need to determine the smooth equivalence class of natural $M_m$--actions on the genus-one manifolds $S^4$, $\CP^2$, and $S^1\times S^3$.
	Using Corollary~\ref{cor:natural}, we see that all natural $M_m$--actions on a fixed $X$ are smoothly equivalent, and so to determine the smooth equivalence class of the natural actions, it suffices to first exhibit an $M_m$--action whose smooth equivalence class we understand on each such manifold, then show that this action is natural by demonstrating a genus-one natural equivariant trisection of this $M_m$ action.
	This is done over the following propositions.
	
\begin{proposition}
	Consider the subgroup $G = ((O(2)\times O(2))\rtimes \Z_2)<SO(4)$ acting on $S^3$ as the setwise stabilizer of a genus-one Heegaard handlebody of $S^3$.
	Let $M_m$ act on $S^3$ as the obvious $(\Z_m\times \Z_m)\rtimes \Z_2$ subgroup of $G$.
	Then the suspension of this $M_m$--action on $S^3$ is a natural $M_m$--action on $S^4$.
\end{proposition}
	
\begin{proof}
	The $M_m$--action on $S^3$ admits a genus-one invariant Heegaard splitting by the invariant handlebody and its complement, and the action of $M_m$ on the central surface is natural by construction. The suspension of this action then admits a natural genus-one equivariant trisection by Proposition~\ref{prop:suspension_heegaard}.
\end{proof}
	
\begin{proposition}
\label{prop:CP2_genus_1}
	Recall the action of $U(2)<PU(3)$ on $\CP^2$ described in Subsection~\ref{subsec:PU(3)_trisections}.
	Consider the subgroup $U(1)\times U(1) < U(2)$ given by the diagonal matrices, and let $G = (U(1)\times U(1))\rtimes \Z_2$ be the subgroup generated by $U(1)\times U(1)$ and compex conjugation.
	Let $M_m$ act on $\CP^2$ as the obvious $(\Z_m\times \Z_m)\rtimes \Z_2$ subgroup of $G$.
	This $M_m$--action on $\CP^2$ is natural.
\end{proposition}
	
\begin{proof}
	The moment map $\mu\colon \CP^2\to\R^2$ described in Theorem~\ref{thm:massey-kuiper} is $M_m$--invariant, and hence the standard trisection of $\CP^2$ is $M_m$--invariant.
	Furthermore, by observation, the $M_m$--action on each handlebody is natural.
	So, to show that this trisection is $M_m$--equivariant, it suffices to show that the actions on the sectors are linearly parted.
	As in Theorem~\ref{thm:massey-kuiper}, we see that each sector of $\im(\mu)$ is a cone on an arc with cone point the unique vertex in that sector.
	Since $\mu$ is $M_m$--invariant, the pre-image of this cone structure shows that each $4$--ball sector of the trisection is the equivariant cone on its boundary and hence the trisection is $M_m$--equivariant as desired.
\end{proof}
	
	We remark that one can recover the same trisection via the construction described in Subsection~\ref{subsec:PU(3)_trisections}:
	The natural map from the $U(1)\times U(1)$ subgroup described above into $PU(2)$ is the quotient by the diagonal.
	Thus, the induced $U(1)\times U(1)$ action on $\CP^1$ is a circle group acting by rotation and fixing two antipodal points, and the $\Z_2$--action is by reflection in a great circle passing through both antipodal points.
	We can then recover the genus-one natural $M_m$--equivariant trisection above via Proposition~\ref{prop:PU(2)_Heegaard} by choosing $\Gamma_1$ to be the single vertex graph given by one of the two antipodal fixed points and choosing $\Gamma_2$ to be the single vertex graph corresponding to the other fixed point.
	
\begin{proposition}
	Let $O(2)\rtimes \Z_2$ act on $S^1$ so that the first factor acts by rotation and the second factor acts as reflection fixing an $S^0$.
	Let $O(2)\rtimes \Z_2$ act on $S^3$ as the double suspension of the above $(O(2)\rtimes \Z_2)$--action.
	Let $G = (O(2)\times O(2))\rtimes \Z_2$ act on $S^1\times S^3$ so that $O(2)\times O(2)$ acts as the product of the above two $O(2)$ actions on $S^1$ and $S^3$ and so that $\Z_2$ acts by simultaneous reflection of both $S^1$ and $S^3$ across their equatorial spheres.
	Let $M_m$ act on $S^1\times S^3$ as the obvious $(\Z_m\times \Z_m)\rtimes \Z_2$ subgroup of $G$.
	This $M_m$--action on $S^1\times S^3$ is natural. 
\end{proposition}
	
\begin{proof}
	View $S^3$ as the suspension of an equatorial $S^2$, and view this $S^2$ as the suspension of its equator circle $E$.
	View the $(O(2)\rtimes\Z_2)$--action via this double-suspension.
	This equator $E\subseteq S^3$ is clearly invariant.
	Consider the following three invariant disks bounded by $E$: 
	Let $D_1$ and $D_2$ be the upper and lower hemispheres of the equatorial $S^2$.
	Choose one of the two $3$--ball cones on $S^2$ making up the suspension $S^3$, and let $D_3$ be the cone on $E$ in that $B^3$.
	These three disks divide $S^3$ into three balls $B_1$, $B_2$, and $B_3$.
	
	Now, consider the trisection of $S^1\times S^3$ whose sectors are $X_i = S^1\times B_i$.
	Clearly this trisection is $M_m$--invariant, and by observation the action of $M_m$ on each handlebody $H_i = S^1\times D_i$ is natural.
	To show that this is an equivariant trisection, in each sector $X_i$ consider the parting ball-system given by the orbit of balls $M_m\cdot B_i$ and note that each $4$--ball in $X_i\setminus\nu(M_m\cdot B_i)$ is the equivariant product of $B_i$ and an interval, so has linear stabilizer and hence this is a linearly parting ball-system as desired.
\end{proof}

All the $M_m$ actions constructed here respect the geometry of the underlying manifold; hence, the induced actions by any subgroups do so, as well.
We can thus close this section with the following corollary.

\begin{corollary}
\label{cor:geometric}
	Every group action admitting a genus-one equivariant trisection is geometric.
\end{corollary}

\subsection{Genus-two equivariant trisections}
\label{subsec:genus_two}

In this subsection, we give a partial classification of genus-two equivariant trisections.

In \cite{Zim_96_Genus-actions-of-finite}, Zimmermann introduces a method of studying equivariant Heegaard splittings by focusing on those that are particularly symmetric, i.e. those carrying a maximal group action like the natural actions discussed in the previous subsection.
We will now make further use of this technique to study equivariant genus-two trisections. 

Zimmermann gives the following definitions: the \emph{handlebody-genus} of a finite group $G$ is the smallest genus $g$ for which $G$ acts faithfully on a genus--$g$ ($3$--dimensional) handlebody, and a group action realizing this genus is called a \emph{genus action} of $G$.
The \emph{strong genus} of $G$ is the smallest such genus that is larger than $1$, and the associated actions are \emph{strong genus actions}.
(For most groups $G$ these genera coincide, but since all cyclic and dihedral groups act in a simple way on solid tori,  analysis of those groups benefits from excluding this trivial case; cf. the previous subsection.)
Recall that a finite group acting faithfully on a handlebody of genus $g\geq 2$ can have order at most $12(g-1)$; an action of a group $G$ on a genus--$g$ handlebody realizing this bound is automatically a strong genus action.
We call this special case a \textit{maximally symmetric action}.
We now adapt the main definition in \cite{Zim_96_Genus-actions-of-finite}.

\begin{definition}
\label{def:maximally_symmetric}
	A $G$--equivariant trisection is \emph{strongly minimal} if the action of $G$ on each handlebody of the spine is a strong genus action.
	Moreover, a $G$--equivariant trisection of genus--$g$ is \emph{maximally symmetric} if $G$ attains the order bound $12(g-1)$.
\end{definition}

It appears that the full classification of genus-two, equivariant trisections should be feasible, but will likely require much technical analysis.
For this reason, rather than pursuing a complete classification, we instead classify the strongly minimal and maximally symmetric equivariant trisections of genus-two.

\begin{figure}[ht!]
     \centering
     \includegraphics[width=.9\textwidth]{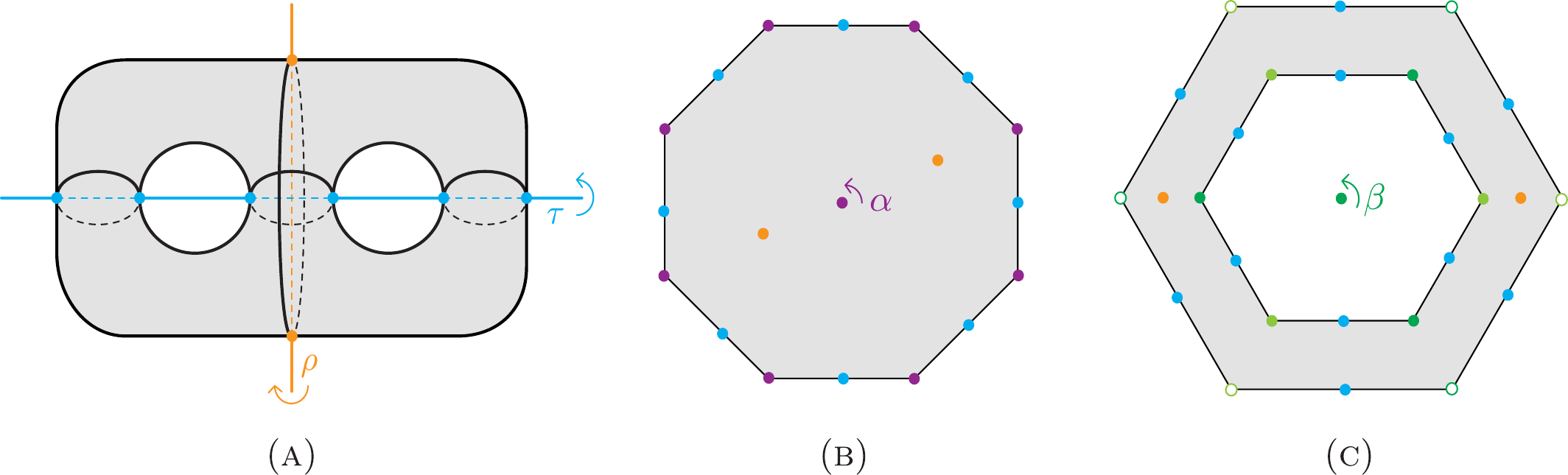}
     \caption{Four symmetries of the genus-two surface}
     \label{fig:genus_two_symmetries}
\end{figure}

Following Zimmerman~\cite{Zim_96_Finite-groups-of-outer,Zim_96_Genus-actions-of-finite}, we know that the only two strong genus actions on a genus-two handlebody are actions of $D_4$ and $D_6$, with $D_6$ being maximally symmetric.
See Figure~\ref{fig:genus_two_symmetries} for a description of these actions on the boundary surface.
The elements $\tau$ and $\rho$ are common to both actions; $\tau$ is the hyperelliptic involution and fixes six points, while $\rho$ is an involution fixing two points.
The $D_4$ action is generated by $\rho$ and $\alpha$, where $\alpha$ is order--$4$ and fixes two points.
The $D_6$ action is generated by $\rho$ and $\beta$, where $\beta$ is order--$6$ and fixed-point free.
Finally, note that $\alpha^2 = \beta^3 = \tau$, so $\rho$ and $\tau$ are contained in both groups as claimed.
In Figure~\ref{fig:genus_two_symmetries}, the fixed points of $\tau$ and $\rho$ are indicated in all three diagrams, with the caveat that the fixed points of $\alpha$ are also fixed points of $\tau$.
In Figures~\ref{fig:genus_two_symmetries}(\textsc{b}) and~(\textsc{c}), the genus-two surface is obtained by identifying each pair of parallel, congruent segments.

The $D_6$--action can be seen as follows: take a $Y$--graph lying in an equatorial disk in a $3$--ball, and consider the double; this gives a $\theta$--graph in the $3$--sphere.
The action of $D_6 = D_3\times \Z_2$ is the product of the obvious $D_3$--action on the $Y$--graph with the orientation preserving involution of $S^3$ that is given by rotation about the equatorial circle of the Heegaard $2$--sphere boundary of the balls.
The neighborhood of this graph is the maximally symmetric $D_6$--action.
Furthermore, note that the genus-two equivariant Heegaard splitting of $S^3$ given by this handlebody and its complement is maximally symmetric in both handlebodies.

\begin{theorem}
\label{thm:genus_two}
	Suppose that $(X,\TT,G)$ is a strongly minimal, genus-two equivariant trisection.
	Then, up $G$--diffeomorphism, $\TT$ is uniquely determined by the pair $(X,G)$, which is one of
	$$\left(\#^2 (S^1\times S^3), D_4\right), \hspace{4mm} \left(\#^2 (S^1\times S^3), D_6\right), \hspace{4mm} \text{ or } \hspace{4mm}\left(S^4, D_6\right).$$
	
	Each such action of $G$ on $\#^2 (S^1\times S^3)$ is the double of the product action on $V\times I$, where the action on $V$ is the strong genus action of $G$ on the genus-two handlebody $V$, and $G$ acts trivially on the interval factor.
	
	The action of $D_6$ on $S^4$ is linear, arising as the suspension of a genus-two $D_6$--equivariant Heegaard splitting of $S^3$.
	These three genus-two equivariant trisections are shown in Figure~\ref{fig:genus_two_diagrams}.
\end{theorem}

\begin{proof}
	\textbf{Case 1:} $G \cong D_4$.
	
	Let $\TT$ be a strongly minimal $G$--equivariant trisection, and let $G' = \langle \alpha \rangle$.
	Referring to the octagon in Figure~\ref{fig:genus_two_symmetries}(\textsc{b}), the singular set of the $G'$--action on the central surface $\Sigma$ of $\TT$ consists of the two fixed points (indicated in purple), which have local degree four, and the two orbits of the four midpoints of the edges of the octagon (indicated in cyan), which have local degree two.
	It is easy to verify that $\Sigma/G'$ is a 2--sphere.
	Since $H_i/G'$ is a handlebody (by the~\ref{ELT}), it follows that $H_i/G'$ is a 3--ball.
	By Corollary~\ref{coro:cyclic_3d}, the quotient $(H_i^*,\Tt_i^*) = (H_i/G',\Sing(G')/G')$ is a trivial, 2--stranded tangle orbifold (also called a solid pillow), with cone-points of order four on one strand and order two on the other strand.
	A disk $D_i^*$ separating the strands of $H_i^*$ lifts to a cut-system $D_i$ for $H_i$.
	
	Let $Y_i = H_i\cup\overline H_{i+1}$, and let $L_i$ be the singular set of the induced $G'$--action on $Y_i$.
	Suppose $Y_i\cong S^3$.
	Note that since $G'$ has only one nontrivial subgroup, the fixed-point set of the $\Z_2$ subgroup of $G'$ is the entire singular set of $G'$.
	By the Smith Conjecture, this fixed-point set is an unknot, contradicting the fact that $L_i$ has two components of different singular orders.
	It follows that $Y_i$ is not $S^3$.

	Since $Y_i\not\cong S^3$, by the Equivariant Sphere Theorem~\cite{MeeSimYau_82_Embedded-minimal}, $L_i$ is either composite or an unlink.
	However, only the latter option is possible for two-bridge links.
	It follows that for each $i\in\Z_3$, the branching set $L_i^*$ is an unlink.
	
	Each $(H_i^*,\Tt_i^*)$ is a rational tangle orbifold given by a rational number $a_i/b_i$ for each $i\in\Z_3$.
	After a diffeomorphism of the common boundaries of the tangles $(H_i^*,\Tt_i^*)$, we can assume $a_1/b_1 = 0/1$.
	Since each $L_i$ is an unlink, it follows that $a_2/b_2 = a_3/b_3 = 0/1$, as well.
	Therefore, the $G$--action on $\Sigma$ extends in the same way across all three handlebodies, and $\TT$ is diffeomorphic to the equivariant trisection described by the diagram shown as Figure~\ref{fig:genus_two_diagrams}(\textsc{a}).
	This shows that the $D_4$ action on $\#^2(S^1\times S^3)$ is the double of $H_1\times I$ as claimed.
	So, there is a unique $D_4$--equivariant, genus-two trisection, as claimed.

\begin{figure}[ht!]
     \centering
     \includegraphics[width=.9\textwidth]{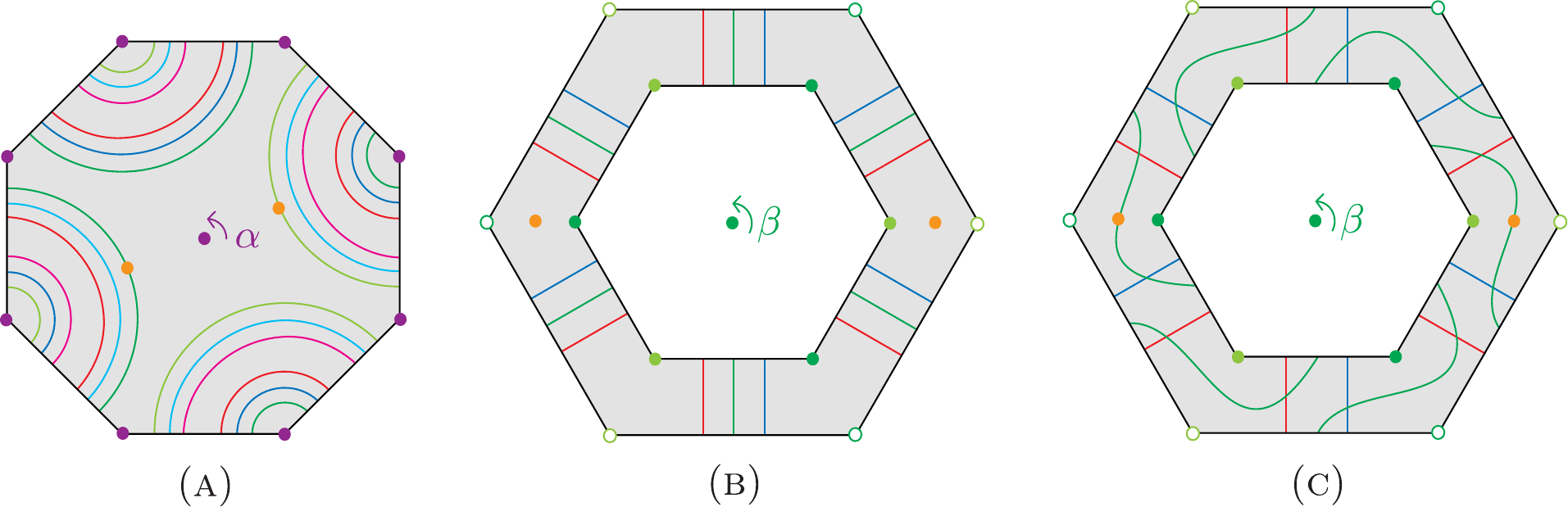}
     \caption{The genus-two equivariant trisections}
     \label{fig:genus_two_diagrams}
\end{figure}
		
	\textbf{Case 2:} $G \cong D_6$.
	
	Again, let $\TT$ be a strongly minimal $G$--equivariant trisection, and let $G'=\langle\beta\rangle$.
	This case is similar to the previous case, except that the local degrees are three and two, so $\beta$ has no fixed points; cf. Figure~\ref{fig:genus_two_diagrams}(\textsc{c}).
	As above, if $Y_i\not\cong S^3$, then by the Equivariant Sphere Theorem, $L_i^*$ is an unlink.
	
	If $Y\cong S^3$, consider the fixed-point set $J$ of $\langle \beta^2\rangle$, which is an unknot by the Smith Conjecture.
	The unknot $J$ is $\beta$--invariant, and is either rotated or reflected by $\beta$; since the singular set of $G'$ is a link, $J$ must be rotated.
	From this, we determine that the singular set $L_i$ is a Hopf link, as is the branch locus $L_i^*$.
	
	After a diffeomorphism of the boundary of the quotient tangles, we can assume that $a_1/b_1 = 0/1$ and $a_2/b_2,a_3/b_3\in\{-2/1,0/1,2/1\}$.
	In order for $L_2^*$ to be a Hopf link or unlink, it is not possible for both $2/1$ and $-2/1$ to appear.
	So, after potentially more diffeomorphisms and permuting the indices, we can assume $a_1/b_1 = a_2/b_2 = 0/1$, and $a_3/b_3\in\{0/1,2/1\}$.
	 
	If $a_3/b_3 = 0/1$, then the action extends in the same way across all three handlebodies, and $\TT$ is diffeomorphic to the equivariant trisection described by the diagram shown as Figure~\ref{fig:genus_two_diagrams}(\textsc{b}), again showing that this action is a double.
	If $a_3/b_3 = 2/1$, then $L_2^*$ is a Hopf link, and $Y_2$ is the cyclic, 6--fold branched cover of $L_2^*$, described above, with branching index three along one component and two along the other.
	Unlike Case 1 above, this is not a contradiction, and the essential disks $D_i^*$ for the two tangles can be lifted to give the diagram shown as Figure~\ref{fig:genus_two_diagrams}(\textsc{c}).
	Note that forgetting the blue curves on this diagram gives a maximally symmetric genus-two $D_6$--equivariant Heegaard splitting of $S^3$ whose suspension is the diagram shown.
\end{proof}

\section{Open questions and further directions}
\label{sec:open_questions_further_directions}

In this section, we discuss several open questions closely related to the present work that we believe to be of broad interest in smooth 4--dimensional topology.
We then discuss some directions for further work that appear fruitful.

\subsection{Open questions}
\label{subsec:open_questions}

We start by applying the perspective of equivariant bridge trisections to the Unknot Linearization Conjecture from~\cite{MeiSco_LP}.
A linear $\Z_p$--action on $S^4$ fixing an $S^2$ admits an equivariant bridge trisection of genus zero, where $\Z_p$ simply acts by rotation on the central surface of the trisection, and so the Unknot Linearization Conjecture is equivalent to the following.

\begin{conjecture}[Unknot Linearization Redux]
\label{conj:unknot_linear_2}
	If $\mathbb Z_n$ acts smoothly on $S^4$ and the fixed-point set $S$ is an unknotted $2$--sphere, then $S$ admits a 1--bridge, genus-zero, equivariant bridge trisection.
\end{conjecture}

This rephrases Proposition~6.2 of~\cite{MeiSco_LP}: a counterexample $(S^4,S)$ must admit an equivariant bridge trisection by Theorem~\ref{thm:exist_main}, but cannot admit an equivariant bridge trisection with trisection genus less than two: if the equivariant trisection had genus one or zero, the action would be linear, by the work of Subsections~\ref{subsec:genus_0_classification} and~\ref{subsec:genus_1_classification}.
Hence, either the quotient $S^4/\rho$ has trisection genus zero, but the fixed-point set $S/\rho$ has bridge number larger than one (so $S^4/\rho \cong S^4$ and $S/\rho$ must be a nontrivial $\Z$--knot), or the quotient has trisection genus greater than one (and $S^4/\rho$ is exotic).

We now turn to a question about $\mathbb{CP}^2$ related to our work in Subsection~\ref{subsec:PU(3)_trisections}.
Recall the natural action of the group $PU(3)$ on $\mathbb{CP}^2$ that we study in that subsection.
This Lie group action is a source of finite group actions on $\mathbb{CP}^2$.
It very clearly does not contain every finite action, or even every cyclic action on $\mathbb{CP}^2$, however, since the action of $PU(3)$ is trivial on the second homology of $\mathbb{CP}^2$.
(The $\mathbb Z_2$--action given by complex conjugation instead acts by negation on $H_2$, for example.)
All odd-order cyclic actions must act trivially on $H_2$, however, and no nonlinear actions are known, motivating the following conjecture.

\begin{conjecture}[$\mathbb{CP}^2$ Homologically Trivial Linearization]
\label{conj:homo_linear}
	Every homologically trivial $\Z_n$--action on $\mathbb{CP}^2$ is smoothly equivalent to a linear action.
\end{conjecture}

Chen has resolved the above conjecture in the affirmative in the symplectic case, and the conjecture is known to be true topologically~\cite[Section~2]{Che_10_Group-actions-on-4-manifolds:}.
In addition, using $\mathbb Z_p$--actions on $S^4$ fixing a nontrivial knot, one can produce $\mathbb Z_p\oplus \mathbb Z_p$--actions which are nonlinear, but in these examples each cyclic subgroup has a linear action~\cite{HamLeeMad_89_Rigidity-of-certain-finite}.
We can again reformulate this conjecture in terms of trisection genus.

\begin{conjecture}[$\mathbb{CP}^2$ Homologically Trivial Linearization]
\label{conj:homo_linear_2}
	Every homologically trivial $\Z_n$--action on $\mathbb{CP}^2$ admits a genus-one equivariant trisection.
\end{conjecture}

This conjecture motivates developing an understanding of how the $G$--action on $H_2(X)$ is determined by an equivariant trisection of $X$.
	
\subsection{Further directions}
\label{subsec:further_work}

In this subsection we collect ideas for further work on this material.
Our first topic is uniqueness.
Any two trisections of a given $4$--manifold become diffeomorphic after sufficiently many stabilizations~\cite{GayKir_16_Trisecting-4-manifolds}.
This statement can be construed as a uniqueness statement: Each $4$--manifold admits a unique trisection up to stabilization.
We thus have the following question.

\begin{question}
\label{quest:stabilization_uniqueness}
	Is there a set of ``stabilization moves'' that relate any two $G$--equivariant trisections of a given $4$--dimensional $G$--manifold $X$?
\end{question}

We remark first that the obvious generalization of stabilization alone is not enough to extract a uniqueness statement.
This is due to discrete singular points: stabilizing along an invariant collection of equivariantly boundary-parallel arcs in any obvious way cannot move isolated singular (or fixed) points between sectors in an equivariant trisection.
By Theorem~\ref{thm:exist_main}, especially Remark~\ref{rmk:sector1}, every $4$--dimensional $G$--manifold $X$ admits an equivariant trisection where all singular points occur in the same sector of the trisection.
On the other hand, most examples we construct in Section~\ref{sec:examples} have singular points distributed between sectors.
Simple equivariant stabilizations could not relate trisections arising in these different ways, so a more complex stabilization move would be needed for any kind of uniqueness statement to hold.
We suggest investigating this in the context of orbifold trisections, where a technique similar to that of Gay-Kirby may be adaptable.

This question is also closely related to the issue mentioned in Remark~\ref{rmk:conjugation}: 
If one can establish uniqueness up to stabilization moves and the stabilization moves utilized preserve the diffeomorphism type of $X/G$, then the diffeomorphism type of $X/G$ would be completely determined by the $G$--action on $X$.

In the non-equivariant case, if two bridge trisections of a given surface agree on the spine, then these two surfaces are \emph{isotopic}.
The result we obtain about determining an equivariant bridge trisection from its spine, Theorem~\ref{thm:spines_bridge}, only establishes uniqueness up to equivariant diffeomorphism, due to a limitation in the proof of~\cite[Theorem~5.9]{MeiSco_LP}; see~\cite[Lemma 5.7]{MeiSco_LP} and discussion thereafter for details.
Again, analysis of quotient orbifolds could be useful in obtaining statements about equivariant isotopy.

For equivariant bridge trisections, we also have an analog for Question~\ref{quest:stabilization_uniqueness}:
Determine a set of trisection-based moves that relate any two equivariant bridge trisections of a given invariant surface-link.
We also lack a key tool in the application of bridge trisections in the non-equivariant setting:
Non-equivariantly, a given surface can always be isotoped to lie in bridge position with respect to any given trisection, which among other things allows us to represent any surface-link in $S^4$ via a triplane diagram.
Our existence result for equivariant bridge trisections (Theorem~\ref{thm:exist_main}) creates a bespoke equivariant bridge trisection for a given surface that may have large trisection genus, so we cannot presently get equivariant triplane diagrams.
This gives the following question.

\begin{question}
\label{ques:tri-plane}
	Let $G$ act on $S^4$ and let $\mathfrak T$ be a genus-zero $G$--equivariant trisection of $S^4$.
	Let $S$ be a $G$--invariant knotted surface in $S^4$.
	Can $S$ be equivariantly isotoped to lie in bridge position with respect to $\mathfrak T$?
\end{question}

Next, in Remark~\ref{rmk:orb_tri}, we outline a potential definition for an orbifold (bridge) trisection of a $4$--dimensional orbifold.
We believe that this idea should be developed further.

A major source of potential future work lies in trisecting specific, interesting group actions.
There are many such potential examples, and here we give a few actions that we are particularly interested in.

The first examples of equivariant trisections of $S^4$ with respect to non-linear cyclic actions occur at genus three:
The 2--fold branched covering of the 3--twist-spun trefoil corresponds to a non-linear $\Z_2$--action on $S^4$~\cite{Gor_74_On-the-higher-dimensional-Smith} that admits a genus-three equivariant trisection via the techniques of Subsection~\ref{subsec:branched_coverings} and~\cite[Section~5.2]{MeiZup_17_Bridge-trisections}.
It would be interesting to give a systematic analysis of the trisections arising in this way.
These coverings are closely related to the classification of circle actions on $S^4$ given by Montgomery and Yang, Fintushel, and Pao~\cite{MonYan_60_Groups-on-Sn-with-principal,Fin_76_Locally-smooth-circle,Pao_78_Nonlinear-circle-actions}; cyclic actions embedding in these circle actions would be an interesting class of actions to investigate.

Recent exciting work of Miyazawa (in light of follow-up work of Hughes, Kim, and Miller) produced an infinite family of pairs $(\CP^2, \iota_i)$, where $\iota_i\colon\CP^2\to\CP^2$ is a smooth involution with quotient $S^4$, such that the $(\CP^2, \iota_i)$ are equivariantly homeomorphic to $(\CP^2,\tau)$, where $\tau$ denotes complex conjugation, but pairwise equivariantly non-diffeomorphic~\cite{HugKimMil_24_Branched-covers-of-twist-roll,Miy_23_A-gauge-theoretic-invariant}.
These actions are branched covering actions, and it would be interesting to apply the techniques of Subsection~\ref{subsec:branched_coverings} to understand the complexity of these actions and their fixed-point sets.

In Subsection~\ref{subsec:PU(3)_trisections}, we investigate a subclass of finite subgroups of $PU(3)$ and trisect their induced actions on $\mathbb{CP}^2$.
A few group actions obtained from finite subgroups of $PU(3)$ cannot arise in this way and thus still remain to be equivariantly trisected, many of which have interesting algebro-geometric properties that may be interesting to some readers.

Consider the $\Z_3$ action on $\CP^2$ arising from cyclic permutation of the homogeneous coordinates, i.e. the generator $\sigma$ acts by $[z_1:z_2:z_3] \mapsto [z_3:z_1:z_2]$. This action behaves nicely with respect to the standard moment map-derived trisection of $\CP^2$, but not as an equivariant trisection:
instead, $\sigma$ cyclically permutes the sectors. 
Many group actions on $4$--manifolds arise in this way, and can be seen via group actions on the trisection diagram which permute the colors of the cut-curve systems; see, for example, the trisections of exotic complex surfaces in~\cite{LamMei_22_Bridge-trisections-in-rational}.
It would be interesting to systematically equivariantly trisect actions arising in this way, and study the manifolds supporting cyclically-permuted trisections.

We can treat the above $\Z_3$--action in specific, however, with the help of some representation theory.
Clearly, $\sigma$ is the quotient of the linear transformation $\sigma'\colon \C^3\to \C^3$ given by cyclically permuting basis elements $\{e_1, e_2, e_3\}$.
Let $\zeta$ be the $3^\text{rd}$ root of unity with positive imaginary part.
In the basis $\{e_1+e_2+e_3, e_1+\zeta e_2 + \zeta^2 e_3, e_1+\zeta^2 e_2 + \zeta e_3\}$, $\sigma'$ is represented by the diagonal matrix $\text{Diag}(1, \zeta, \zeta^2)$.
Via a suitable change-of-basis, which descends to a diffeomorphism of $\CP^2$, we can see therefore that $\sigma$ generates a subgroup of $U(2)<PU(3)$, which we trisect in Subsection~\ref{subsec:PU(3)_trisections}.
In fact, $\langle\sigma\rangle$ admits a genus-one trisection as in Proposition~\ref{prop:CP2_genus_1}.

\bibliographystyle{amsalpha}
\bibliography{EqTri.bib}

\providecommand{\bysame}{\leavevmode\hbox to3em{\hrulefill}\thinspace}
\providecommand{\MR}{\relax\ifhmode\unskip\space\fi MR }
\providecommand{\MRhref}[2]{%
  \href{http://www.ams.org/mathscinet-getitem?mr=#1}{#2}
}
\providecommand{\href}[2]{#2}
\begin{thebibliography}{BCKM24}

\bibitem[Arn88]{Arn_88_The-branched-covering-bf-Crm-P2to}
V.~I. Arnol'd, \emph{The branched covering {${\bf C}{\rm P}^2\to S^4$},
  hyperbolicity and projective topology}, Sibirsk. Mat. Zh. \textbf{29} (1988),
  no.~5, 36--47, 237. \MR{971226}

\bibitem[BCKM24]{BlaCahKju_24_Note-on-three-fold-branched-covers}
Ryan Blair, Patricia Cahn, Alexandra Kjuchukova, and Jeffrey Meier, \emph{Note
  on three-fold branched covers of {$S^4$}}, Ann. Inst. Fourier (Grenoble)
  \textbf{74} (2024), no.~2, 849--866. \MR{4748188}

\bibitem[Ber24]{Ber_24_Star-Neighborhoods-in-Double}
Aleksandr Berdnikov, \emph{Star neighborhoods in double barycentric
  subdivision},
  \url{http://demonstrations.wolfram.com/StarNeighborhoodsInDoubleBarycentricSubdivision/},
  January 2024.

\bibitem[Bre72]{Bre_72_Introduction-to-compact-transformation-groups}
Glen~E. Bredon, \emph{Introduction to compact transformation groups}, Pure and
  Applied Mathematics, vol. Vol. 46, Academic Press, New York-London, 1972.
  \MR{413144}

\bibitem[Can79]{Can_79_Shrinking-cell-like-decompositions}
J.~W. Cannon, \emph{Shrinking cell-like decompositions of manifolds.
  {C}odimension three}, Ann. of Math. (2) \textbf{110} (1979), no.~1, 83--112.
  \MR{541330}

\bibitem[Che10]{Che_10_Group-actions-on-4-manifolds:}
Weimin Chen, \emph{Group actions on 4-manifolds: some recent results and open
  questions}, Proceedings of the {G}\"okova {G}eometry-{T}opology {C}onference
  2009, Int. Press, Somerville, MA, 2010, pp.~1--21. \MR{2655301}

\bibitem[CK17]{CahKju_17_Singular-branched}
Patricia {Cahn} and Alexandra {Kjuchukova}, \emph{{Singular branched covers of
  four-manifolds}}, arXiv e-prints (2017), arXiv:1710.11562.

\bibitem[CMR23]{CahMatRup_23_Algorithms-for-Computing-Invariants-of-Trisected}
Patricia {Cahn}, Gordana {Matic}, and Benjamin {Ruppik}, \emph{{Algorithms for
  Computing Invariants of Trisected Branched Covers}}, arXiv e-prints (2023),
  arXiv:2308.11689.

\bibitem[DL09]{DinLee_09_Equivariant_Ricci_flow}
Jonathan Dinkelbach and Bernhard Leeb, \emph{{Equivariant Ricci flow with
  surgery and applications to finite group actions on geometric
  $3$--manifolds}}, Geometry and Topology \textbf{13} (2009), no.~2, 1129 --
  1173.

\bibitem[Edm86]{Edm_86_A-topological-proof}
Allan~L. Edmonds, \emph{A topological proof of the equivariant {D}ehn lemma},
  Trans. Amer. Math. Soc. \textbf{297} (1986), no.~2, 605--615. \MR{854087}

\bibitem[Edm87]{Edm_87_Construction-of-group-actions}
Allan~L. Edmonds, \emph{Construction of group actions of four-manifolds},
  Transactions of the American Mathematical Society \textbf{299} (1987), no.~1,
  155--170.

\bibitem[Edm18]{Edm_18_A-survey-of-group-actions-on-4-manifolds}
Allan~L. Edmonds, \emph{A survey of group actions on 4-manifolds}, Handbook of
  group actions. {V}ol. {III}, Adv. Lect. Math. (ALM), vol.~40, Int. Press,
  Somerville, MA, 2018, pp.~421--460. \MR{3888625}

\bibitem[Fin76]{Fin_76_Locally-smooth-circle}
Ronald Fintushel, \emph{Locally smooth circle actions on homotopy
  {$4$}-spheres}, Duke Math. J. \textbf{43} (1976), no.~1, 63--70. \MR{394716}

\bibitem[FM12]{FarMar_12_A-primer-on-mapping-class-groups}
Benson Farb and Dan Margalit, \emph{A primer on mapping class groups},
  Princeton Mathematical Series, vol.~49, Princeton University Press,
  Princeton, NJ, 2012. \MR{2850125}

\bibitem[FS97]{FinSte_97_Rational-blowdowns-of-smooth-4-manifolds}
Ronald Fintushel and Ronald~J. Stern, \emph{Rational blowdowns of smooth
  {$4$}-manifolds}, J. Differential Geom. \textbf{46} (1997), no.~2, 181--235.
  \MR{1484044}

\bibitem[Gif66]{Gif_66_The-generalized-Smith-conjecture}
Charles~H. Giffen, \emph{The generalized {S}mith conjecture}, Amer. J. Math.
  \textbf{88} (1966), 187--198. \MR{198462}

\bibitem[GK16]{GayKir_16_Trisecting-4-manifolds}
David Gay and Robion Kirby, \emph{Trisecting 4-manifolds}, Geom. Topol.
  \textbf{20} (2016), no.~6, 3097--3132. \MR{3590351}

\bibitem[GM22]{GayMei_22_Doubly-pointed-trisection-diagrams}
David Gay and Jeffrey Meier, \emph{Doubly pointed trisection diagrams and
  surgery on 2-knots}, Math. Proc. Cambridge Philos. Soc. \textbf{172} (2022),
  no.~1, 163--195. \MR{4354420}

\bibitem[Gor74]{Gor_74_On-the-higher-dimensional-Smith}
C.~McA. Gordon, \emph{On the higher-dimensional {S}mith conjecture}, Proc.
  London Math. Soc. (3) \textbf{29} (1974), 98--110. \MR{0356073 (50 \#8544)}

\bibitem[GS99]{GomSti_99_4-manifolds-and-Kirby}
Robert~E. Gompf and Andr{{\'a}}s~I. Stipsicz, \emph{{$4$}-manifolds and {K}irby
  calculus}, Graduate Studies in Mathematics, vol.~20, American Mathematical
  Society, Providence, RI, 1999. \MR{1707327 (2000h:57038)}

\bibitem[HH11]{HamHau_11_Conjugation-spaces-and-4-manifolds}
Ian Hambleton and Jean-Claude Hausmann, \emph{Conjugation spaces and
  4-manifolds}, Math. Z. \textbf{269} (2011), no.~1-2, 521--541. \MR{2836082}

\bibitem[Hir94]{Hir_94_Differential-topology}
Morris~W. Hirsch, \emph{Differential topology}, Graduate Texts in Mathematics,
  vol.~33, Springer-Verlag, New York, 1994, Corrected reprint of the 1976
  original. \MR{1336822}

\bibitem[HKM24]{HugKimMil_24_Branched-covers-of-twist-roll}
Mark {Hughes}, Seungwon {Kim}, and Maggie {Miller}, \emph{{Branched covers of
  twist-roll spun knots}}, arXiv e-prints (2024), arXiv:2402.11706.

\bibitem[HL88]{HamLee_88_Finite-group-actions}
Ian Hambleton and Ronnie Lee, \emph{Finite group actions on {${\rm P}^2({\bf
  C})$}}, J. Algebra \textbf{116} (1988), no.~1, 227--242. \MR{944157}

\bibitem[HLM89]{HamLeeMad_89_Rigidity-of-certain-finite}
Ian Hambleton, Ronnie Lee, and Ib~Madsen, \emph{Rigidity of certain finite
  group actions on the complex projective plane}, Comment. Math. Helv.
  \textbf{64} (1989), no.~4, 618--638. \MR{1022999}

\bibitem[Ill78]{Ill_78_Smooth-equivariant-triangulations}
S\"{o}ren Illman, \emph{Smooth equivariant triangulations of {$G$}-manifolds
  for {$G$} a finite group}, Math. Ann. \textbf{233} (1978), no.~3, 199--220.
  \MR{500993}

\bibitem[JMMZ22]{JosMeiMil_22_Bridge-trisections-and-classical}
Jason Joseph, Jeffrey Meier, Maggie Miller, and Alexander Zupan, \emph{Bridge
  trisections and classical knotted surface theory}, Pacific J. Math.
  \textbf{319} (2022), no.~2, 343--369. \MR{4482720}

\bibitem[Kep22]{Kep_22_An-algorithm-taking-Kirby-diagrams}
Willi Kepplinger, \emph{An algorithm taking {K}irby diagrams to trisection
  diagrams}, Pacific J. Math. \textbf{318} (2022), no.~1, 109--126.
  \MR{4460230}

\bibitem[Kui74]{Kui_74_The-quotient-space-of-bf-CP2-by-complex-conjugation}
Nicolaas~H. Kuiper, \emph{The quotient space of {${\bf C}P(2)$} by complex
  conjugation is the {$4$}-sphere}, Math. Ann. \textbf{208} (1974), 175--177.
  \MR{346817}

\bibitem[Lan16]{Lan_16_Characterization-of-finite-groups}
Christian Lange, \emph{Characterization of finite groups generated by
  reflections and rotations}, J. Topol. \textbf{9} (2016), no.~4, 1109--1129.
  \MR{3620454}

\bibitem[Lan19]{Lan_19_When-is-the-underlying-space}
\bysame, \emph{When is the underlying space of an orbifold a manifold?}, Trans.
  Amer. Math. Soc. \textbf{372} (2019), no.~4, 2799--2828. \MR{3988594}

\bibitem[LCM22]{LamMei_22_Bridge-trisections-in-rational}
Peter Lambert-Cole and Jeffrey Meier, \emph{Bridge trisections in rational
  surfaces}, J. Topol. Anal. \textbf{14} (2022), no.~3, 655--708. \MR{4493476}

\bibitem[LCMS21]{LamMeiSta_21_Symplectic-4-manifolds-admit-Weinstein}
Peter Lambert-Cole, Jeffrey Meier, and Laura Starkston, \emph{Symplectic
  4-manifolds admit {W}einstein trisections}, J. Topol. \textbf{14} (2021),
  no.~2, 641--673. \MR{4286052}

\bibitem[Liv82]{Liv_82_Surfaces-bounding-the-unlink}
Charles Livingston, \emph{Surfaces bounding the unlink}, Michigan Math. J.
  \textbf{29} (1982), no.~3, 289--298. \MR{674282}

\bibitem[LP72]{LauPoe_72_A-note-on-4-dimensional-handlebodies}
Fran{\c{c}}ois Laudenbach and Valentin Po{{\'e}}naru, \emph{A note on
  {$4$}-dimensional handlebodies}, Bull. Soc. Math. France \textbf{100} (1972),
  337--344. \MR{0317343 (47 \#5890)}

\bibitem[Mas73]{Mas_73_The-quotient-space-of-the-complex-projective}
W.~S. Massey, \emph{The quotient space of the complex projective plane under
  conjugation is a {$4$}-sphere}, Geometriae Dedicata \textbf{2} (1973),
  371--374. \MR{341511}

\bibitem[MB84]{MorBas_84_The-Smith-conjecture}
John~W. Morgan and Hyman Bass (eds.), \emph{The {S}mith conjecture}, Pure and
  Applied Mathematics, vol. 112, Academic Press, Inc., Orlando, FL, 1984,
  Papers presented at the symposium held at Columbia University, New York,
  1979. \MR{758459}

\bibitem[MBD61]{MilBliDic_61_Theory-and-applications-of-finite-groups}
G.~A. Miller, H.~F. Blichfeldt, and L.~E. Dickson, \emph{Theory and
  applications of finite groups}, Dover Publications, Inc., New York, 1961.
  \MR{123600}

\bibitem[Mi78]{Mih_78_Finite-imprimitive-groups}
M.~A. Miha\u~ilova, \emph{Finite imprimitive groups generated by
  pseudoreflections}, Studies in geometry and algebra ({R}ussian), Kirgiz. Gos.
  Univ., Frunze, 1978, pp.~82--93. \MR{608821}

\bibitem[{Miy}23]{Miy_23_A-gauge-theoretic-invariant}
Jin {Miyazawa}, \emph{{A gauge theoretic invariant of embedded surfaces in
  $4$-manifolds and exotic $P^2$-knots}}, arXiv e-prints (2023),
  arXiv:2312.02041.

\bibitem[MMZ89]{McCMilZim_89_Group-actions-on-handlebodies}
Darryl McCullough, Andy Miller, and Bruno Zimmermann, \emph{Group actions on
  handlebodies}, Proc. London Math. Soc. (3) \textbf{59} (1989), no.~2,
  373--416. \MR{1004434}

\bibitem[MS25]{MeiSco_LP}
Jeffrey Meier and Evan Scott, \emph{{An equivariant Laudenbach-Po\'enaru
  theorem}}, arXiv:2501.10524, 2025.

\bibitem[MSY82]{MeeSimYau_82_Embedded-minimal}
William Meeks, III, Leon Simon, and Shing~Tung Yau, \emph{Embedded minimal
  surfaces, exotic spheres, and manifolds with positive {R}icci curvature},
  Ann. of Math. (2) \textbf{116} (1982), no.~3, 621--659. \MR{678484}

\bibitem[MSZ16]{MeiSchZup_16_Classification-of-trisections}
Jeffrey Meier, Trent Schirmer, and Alexander Zupan, \emph{Classification of
  trisections and the {G}eneralized {P}roperty {R} {C}onjecture}, Proc. Amer.
  Math. Soc. \textbf{144} (2016), no.~11, 4983--4997. \MR{3544545}

\bibitem[MY60]{MonYan_60_Groups-on-Sn-with-principal}
D.~Montgomery and C.~T. Yang, \emph{Groups on {$S\sp{n}$} with principal orbits
  of dimension {$n-3$}. {I}, {II}}, Illinois J. Math. \textbf{4} (1960),
  507--517. \MR{125902}

\bibitem[MY79]{MeeYau_79_The-classical-Plateau}
William~H. Meeks, III and Shing~Tung Yau, \emph{The classical {P}lateau problem
  and the topology of {$3$}-manifolds}, Minimal submanifolds and geodesics
  ({P}roc. {J}apan-{U}nited {S}tates {S}em., {T}okyo, 1977), North-Holland,
  Amsterdam-New York, 1979, pp.~101--102. \MR{574258}

\bibitem[MY80]{MeeYau_80_Topology-of-three-dimensional}
\bysame, \emph{Topology of three-dimensional manifolds and the embedding
  problems in minimal surface theory}, Ann. of Math. (2) \textbf{112} (1980),
  no.~3, 441--484. \MR{595203}

\bibitem[MZ17]{MeiZup_17_Bridge-trisections}
Jeffrey Meier and Alexander Zupan, \emph{Bridge trisections of knotted surfaces
  in {$S^4$}}, Trans. Amer. Math. Soc. \textbf{369} (2017), no.~10, 7343--7386.
  \MR{3683111}

\bibitem[MZ18]{MeiZup_18_Bridge-trisections}
\bysame, \emph{Bridge trisections of knotted surfaces in 4-manifolds},
  Proceedings of the National Academy of Sciences \textbf{115} (2018), no.~43,
  10880--10886.

\bibitem[Pao78]{Pao_78_Nonlinear-circle-actions}
Peter~Sie Pao, \emph{Nonlinear circle actions on the {$4$}-sphere and twisting
  spun knots}, Topology \textbf{17} (1978), no.~3, 291--296. \MR{508892}

\bibitem[Per02]{Per_02_The-entropy-formula}
Grisha Perelman, \emph{The entropy formula for the ricci flow and its geometric
  applications}, arXiv:0211159, 2002.

\bibitem[Per03a]{Per_03_Finite-extinction}
\bysame, \emph{Finite extinction time for the solutions to the ricci flow on
  certain three-manifolds}, arXiv:0307245, 2003.

\bibitem[Per03b]{Per_03_Ricci-flow}
\bysame, \emph{Ricci flow with surgery on three-manifolds}, arXiv:0303109,
  2003.

\bibitem[Sch80]{Sch_80_Lifting-smooth-homotopies}
Gerald~W. Schwarz, \emph{Lifting smooth homotopies of orbit spaces}, Inst.
  Hautes \'{E}tudes Sci. Publ. Math. (1980), no.~51, 37--135. \MR{573821}

\bibitem[Was69]{Was_69_Equivariant-differential-topology}
Arthur~G. Wasserman, \emph{Equivariant differential topology}, Topology
  \textbf{8} (1969), 127--150. \MR{250324}

\bibitem[Wil87]{Wil_87_Group-actions-on-the-complex}
Dariusz~M. Wilczy\'nski, \emph{Group actions on the complex projective plane},
  Trans. Amer. Math. Soc. \textbf{303} (1987), no.~2, 707--731. \MR{902793}

\bibitem[Wil20]{Wil_20_Trisections-of-flat-surface-bundles}
Marla Williams, \emph{Trisections of {F}lat {S}urface {B}undles over
  {S}urfaces}, 81, Thesis (Ph.D.)--The University of Nebraska - Lincoln.
  \MR{4144669}

\bibitem[Zim96a]{Zim_96_Finite-groups-of-outer}
Bruno Zimmermann, \emph{Finite groups of outer automorphisms of free groups},
  Glasgow Math. J. \textbf{38} (1996), no.~3, 275--282. \MR{1417356}

\bibitem[Zim96b]{Zim_96_Genus-actions-of-finite}
\bysame, \emph{Genus actions of finite groups on {$3$}-manifolds}, Michigan
  Math. J. \textbf{43} (1996), no.~3, 593--610. \MR{1420594}

\bibitem[Zud08]{Zud_08_Branched-coverings}
Daniele Zuddas, \emph{Branched coverings and 4-manifolds}, Ph.D. thesis, Scuola
  {N}ormale {S}uperiore {P}isa, 2008.

\end{thebibliography}

\end{document}